\ifx\shlhetal\undefinedcontrolsequence\let\shlhetal\relax\fi
\input amstex
\expandafter\ifx\csname mathdefs.tex\endcsname\relax
  \expandafter\gdef\csname mathdefs.tex\endcsname{}
\else \message{Hey!  Apparently you were trying to
  \string twice.   This does not make sense.} 
\errmessage{Please edit your file (probably \jobname.tex) and remove
any duplicate ``\string\input'' lines} \fi




\catcode`\X=12\catcode`\@=11

\def\n@wcount{\alloc@0\count\countdef\insc@unt}
\def\n@wwrite{\alloc@7\write\chardef\sixt@@n}
\def\n@wread{\alloc@6\read\chardef\sixt@@n}
\def\r@s@t{\relax}\def\v@idline{\par}\def\@mputate#1/{#1}
\def\l@c@l#1X{\firstpart.#1}\def\gl@b@l#1X{#1}\def\t@d@l#1X{{}}

\def\crossrefs#1{\ifx\all#1\let\tr@ce=\all\else\def\tr@ce{#1,}\fi
   \n@wwrite\cit@tionsout\openout\cit@tionsout=\jobname.cit 
   \write\cit@tionsout{\tr@ce}\expandafter\setfl@gs\tr@ce,}
\def\setfl@gs#1,{\def\@{#1}\ifx\@\empty\let\next=\relax
   \else\let\next=\setfl@gs\expandafter\xdef
   \csname#1tr@cetrue\endcsname{}\fi\next}
\def\m@ketag#1#2{\expandafter\n@wcount\csname#2tagno\endcsname
     \csname#2tagno\endcsname=0\let\tail=\all\xdef\all{\tail#2,}
   \ifx#1\l@c@l\let\tail=\r@s@t\xdef\r@s@t{\csname#2tagno\endcsname=0\tail}\fi
   \expandafter\gdef\csname#2cite\endcsname##1{\expandafter
     \ifx\csname#2tag##1\endcsname\relax?\else\csname#2tag##1\endcsname\fi
     \expandafter\ifx\csname#2tr@cetrue\endcsname\relax\else
     \write\cit@tionsout{#2tag ##1 cited on page \folio.}\fi}
   \expandafter\gdef\csname#2page\endcsname##1{\expandafter
     \ifx\csname#2page##1\endcsname\relax?\else\csname#2page##1\endcsname\fi
     \expandafter\ifx\csname#2tr@cetrue\endcsname\relax\else
     \write\cit@tionsout{#2tag ##1 cited on page \folio.}\fi}
   \expandafter\gdef\csname#2tag\endcsname##1{\expandafter
      \ifx\csname#2check##1\endcsname\relax
      \expandafter\xdef\csname#2check##1\endcsname{}%
      \else\immediate\write16{Warning: #2tag ##1 used more than once.}\fi
      \multit@g{#1}{#2}##1/X%
      \write\t@gsout{#2tag ##1 assigned number \csname#2tag##1\endcsname\space
      on page \number\count0.}%
   \csname#2tag##1\endcsname}}
\def\multit@g#1#2#3/#4X{\def\t@mp{#4}\ifx\t@mp\empty%
      \global\advance\csname#2tagno\endcsname by 1 
      \expandafter\xdef\csname#2tag#3\endcsname
      {#1\number\csname#2tagno\endcsnameX}%
   \else\expandafter\ifx\csname#2last#3\endcsname\relax
      \expandafter\n@wcount\csname#2last#3\endcsname
      \global\advance\csname#2tagno\endcsname by 1 
      \expandafter\xdef\csname#2tag#3\endcsname
      {#1\number\csname#2tagno\endcsnameX}
      \write\t@gsout{#2tag #3 assigned number \csname#2tag#3\endcsname\space
      on page \number\count0.}\fi
   \global\advance\csname#2last#3\endcsname by 1
   \def\t@mp{\expandafter\xdef\csname#2tag#3/}%
   \expandafter\t@mp\@mputate#4\endcsname
   {\csname#2tag#3\endcsname\lastpart{\csname#2last#3\endcsname}}\fi}
\def\t@gs#1{\def\all{}\m@ketag#1e\m@ketag#1s\m@ketag\t@d@l p
   \m@ketag\gl@b@l r \n@wread\t@gsin
   \openin\t@gsin=\jobname.tgs \re@der \closein\t@gsin
   \n@wwrite\t@gsout\openout\t@gsout=\jobname.tgs }
\outer\def\localtags{\t@gs\l@c@l}
\outer\def\globaltags{\t@gs\gl@b@l}
\outer\def\newlocaltag#1{\m@ketag\l@c@l{#1}}
\outer\def\newglobaltag#1{\m@ketag\gl@b@l{#1}}

\newif\ifpr@ 
\def\m@kecs #1tag #2 assigned number #3 on page #4.%
   {\expandafter\gdef\csname#1tag#2\endcsname{#3}
   \expandafter\gdef\csname#1page#2\endcsname{#4}
   \ifpr@\expandafter\xdef\csname#1check#2\endcsname{}\fi}
\def\re@der{\ifeof\t@gsin\let\next=\relax\else
   \read\t@gsin to\t@gline\ifx\t@gline\v@idline\else
   \expandafter\m@kecs \t@gline\fi\let \next=\re@der\fi\next}
\def\pretags#1{\pr@true\pret@gs#1,,}
\def\pret@gs#1,{\def\@{#1}\ifx\@\empty\let\n@xtfile=\relax
   \else\let\n@xtfile=\pret@gs \openin\t@gsin=#1.tgs \message{#1} \re@der 
   \closein\t@gsin\fi \n@xtfile}

\newcount\sectno\sectno=0\newcount\subsectno\subsectno=0
\newif\ifultr@local \def\ultralocal{\ultr@localtrue}
\def\firstpart{\number\sectno}
\def\lastpart#1{\ifcase#1 \or a\or b\or c\or d\or e\or f\or g\or h\or 
   i\or k\or l\or m\or n\or o\or p\or q\or r\or s\or t\or u\or v\or w\or 
   x\or y\or z \fi}

\def\resetall{\global\advance\sectno by 1\subsectno=0
   \gdef\firstpart{\number\sectno}\r@s@t}
\def\resetsub{\global\advance\subsectno by 1
   \gdef\firstpart{\number\sectno.\number\subsectno}\r@s@t}
\def\newsection#1\par{\resetall\vskip0pt plus.3\vsize\penalty-250
   \vskip0pt plus-.3\vsize\bigskip\bigskip
   \message{#1}\leftline{\bf#1}\nobreak\bigskip}
\def\subsection#1\par{\ifultr@local\resetsub\fi
   \vskip0pt plus.2\vsize\penalty-250\vskip0pt plus-.2\vsize
   \bigskip\smallskip\message{#1}\leftline{\bf#1}\nobreak\medskip}

\def\t@gsoff#1,{\def\@{#1}\ifx\@\empty\let\next=\relax\else\let\next=\t@gsoff
   \def\@@{p}\ifx\@\@@\else
   \expandafter\gdef\csname#1cite\endcsname##1{\zeigen{##1}}
   \expandafter\gdef\csname#1page\endcsname##1{?}
   \expandafter\gdef\csname#1tag\endcsname##1{\zeigen{##1}}\fi\fi\next}
\def\verbatimtags{\ifx\all\relax\else\expandafter\t@gsoff\all,\fi}
\def\zeigen#1{\hbox{$\langle$}#1\hbox{$\rangle$}}

\def\(#1){\edef\dot@g{\ifmmode\ifinner(\hbox{\noexpand\etag{#1}})
   \else\noexpand\eqno(\hbox{\noexpand\etag{#1}})\fi
   \else(\noexpand\ecite{#1})\fi}\dot@g}

\newif\ifbr@ck
\def\eat#1{}
\def\[#1]{\br@cktrue[\br@cket#1'X]}
\def\br@cket#1'#2X{\def\temp{#2}\ifx\temp\empty\let\next\eat
   \else\let\next\br@cket\fi
   \ifbr@ck\br@ckfalse\br@ck@t#1,X\else\br@cktrue#1\fi\next#2X}
\def\br@ck@t#1,#2X{\def\temp{#2}\ifx\temp\empty\let\neext\eat
   \else\let\neext\br@ck@t\def\temp{,}\fi
   \def\teemp{#1}\ifx\teemp\empty\else\rcite{#1}\fi\temp\neext#2X}
\def\resetbr@cket{\gdef\[##1]{[\rtag{##1}]}}
\def\references{\resetbr@cket\newsection References\par}

\newtoks\symb@ls\newtoks\s@mb@ls\newtoks\p@gelist\n@wcount\ftn@mber
    \ftn@mber=1\newif\ifftn@mbers\ftn@mbersfalse\newif\ifbyp@ge\byp@gefalse
\def\defm@rk{\ifftn@mbers\n@mberm@rk\else\symb@lm@rk\fi}
\def\n@mberm@rk{\xdef\m@rk{{\the\ftn@mber}}%
    \global\advance\ftn@mber by 1 }
\def\rot@te#1{\let\temp=#1\global#1=\expandafter\r@t@te\the\temp,X}
\def\r@t@te#1,#2X{{#2#1}\xdef\m@rk{{#1}}}
\def\b@@st#1{{$^{#1}$}}\def\str@p#1{#1}
\def\symb@lm@rk{\ifbyp@ge\rot@te\p@gelist\ifnum\expandafter\str@p\m@rk=1 
    \s@mb@ls=\symb@ls\fi\write\f@nsout{\number\count0}\fi \rot@te\s@mb@ls}
\def\byp@ge{\byp@getrue\n@wwrite\f@nsin\openin\f@nsin=\jobname.fns 
    \n@wcount\currentp@ge\currentp@ge=0\p@gelist={0}
    \re@dfns\closein\f@nsin\rot@te\p@gelist
    \n@wread\f@nsout\openout\f@nsout=\jobname.fns }
\def\m@kelist#1X#2{{#1,#2}}
\def\re@dfns{\ifeof\f@nsin\let\next=\relax\else\read\f@nsin to \f@nline
    \ifx\f@nline\v@idline\else\let\t@mplist=\p@gelist
    \ifnum\currentp@ge=\f@nline
    \global\p@gelist=\expandafter\m@kelist\the\t@mplistX0
    \else\currentp@ge=\f@nline
    \global\p@gelist=\expandafter\m@kelist\the\t@mplistX1\fi\fi
    \let\next=\re@dfns\fi\next}
\def\symbols#1{\symb@ls={#1}\s@mb@ls=\symb@ls} 
\def\bigsymbol{\textstyle}
\symbols{\bigsymbol\ast,\dagger,\ddagger,\sharp,\flat,\natural,\star}
\def\ftnumbers{\ftn@mberstrue} \def\ftsymbols{\ftn@mbersfalse}
\def\paginal{\byp@ge} \def\resetftnumbers{\ftn@mber=1}
\def\ftnote#1{\defm@rk\expandafter\expandafter\expandafter\footnote
    \expandafter\b@@st\m@rk{#1}}

\long\def\jump#1\endjump{}
\def\ssum{\mathop{\lower .1em\hbox{$\textstyle\Sigma$}}\nolimits}

\def\qed{\nobreak\kern 1em \vrule height .5em width .5em depth 0em}
\def\newneq{\hbox{\rlap{\hbox to 1\wd9{\hss$=$\hss}}\raise .1em 
   \hbox to 1\wd9{\hss$\scriptscriptstyle/$\hss}}}
\def\subsetne{\setbox9 = \hbox{$\subset$}\mathrel{\hbox{\rlap
   {\lower .4em \newneq}\raise .13em \hbox{$\subset$}}}}
\def\supsetne{\setbox9 = \hbox{$\subset$}\mathrel{\hbox{\rlap
   {\lower .4em \newneq}\raise .13em \hbox{$\supset$}}}}

\def\vbar{\mathchoice{\vrule height6.3ptdepth-.5ptwidth.8pt\kern-.8pt}
   {\vrule height6.3ptdepth-.5ptwidth.8pt\kern-.8pt}
   {\vrule height4.1ptdepth-.35ptwidth.6pt\kern-.6pt}
   {\vrule height3.1ptdepth-.25ptwidth.5pt\kern-.5pt}}
\def\f@dge{\mathchoice{}{}{\mkern.5mu}{\mkern.8mu}}
\def\b@c#1#2{{\rm \mkern#2mu\vbar\mkern-#2mu#1}}
\def\b@b#1{{\rm I\mkern-3.5mu #1}}
\def\b@a#1#2{{\rm #1\mkern-#2mu\f@dge #1}}
\def\bb#1{{\count4=`#1 \advance\count4by-64 \ifcase\count4\or\b@a A{11.5}\or
   \b@b B\or\b@c C{5}\or\b@b D\or\b@b E\or\b@b F \or\b@c G{5}\or\b@b H\or
   \b@b I\or\b@c J{3}\or\b@b K\or\b@b L \or\b@b M\or\b@b N\or\b@c O{5} \or
   \b@b P\or\b@c Q{5}\or\b@b R\or\b@a S{8}\or\b@a T{10.5}\or\b@c U{5}\or
   \b@a V{12}\or\b@a W{16.5}\or\b@a X{11}\or\b@a Y{11.7}\or\b@a Z{7.5}\fi}}

\catcode`\X=11 \catcode`\@=12

\sectno=-1   
\localtags

\NoBlackBoxes
\def\cite #1{\rm[#1]}
\newbox\noforkbox \newdimen\forklinewidth
\forklinewidth=0.3pt   
\setbox0\hbox{$\textstyle\bigcup$}
\setbox1\hbox to \wd0{\hfil\vrule width \forklinewidth depth \dp0
                        height \ht0 \hfil}
\wd1=0 cm
\setbox\noforkbox\hbox{\box1\box0\relax}
\def\unionstick{\mathop{\copy\noforkbox}\limits}
\def\nonfork#1#2_#3{#1\unionstick_{\textstyle #3}#2}
\def\nonforkin#1#2_#3^#4{#1\unionstick_{\textstyle #3}^{\textstyle #4}#2}     
%
\setbox0\hbox{$\textstyle\bigcup$}
\setbox1\hbox to \wd0{\hfil{\sl /\/}\hfil}
\setbox2\hbox to \wd0{\hfil\vrule height \ht0 depth \dp0 width
                                \forklinewidth\hfil}
\wd1=0cm
\wd2=0cm
\newbox\doesforkbox
\setbox\doesforkbox\hbox{\box1\box0\relax}
\def\nunionstick{\mathop{\copy\doesforkbox}\limits}

\def\fork#1#2_#3{#1\nunionstick_{\textstyle #3}#2}
\def\forkin#1#2_#3^#4{#1\nunionstick_{\textstyle #3}^{\textstyle #4}#2}     
\documentstyle {amsppt}
\topmatter
\title {CHARACTERIZING AN $\aleph_\epsilon$-SATURATED MODEL OF 
SUPERSTABLE NDOP THEORIES BY ITS $L_{\infty,\aleph_\epsilon}$-THEORY \\
Sh401} \endtitle
\rightheadtext{Characterizing Models of NDOP theories}
\author {Saharon Shelah \thanks{\null\newline
Partially supported by the United States-Israel Binational Science 
Foundation and the NSF, and I thank Alice Leonhardt for the beautiful typing.
\newline
Done Oct/89 \newline
Publ No. 401 \newline
Latest Revision - 96/Sept/10} \endthanks} \endauthor
\affil {Institute of Mathematics \\
The Hebrew University \\
Jerusalem, Israel
\medskip
Department of Mathematics \\
Rutgers University \\
New Brunswick, NJ  USA
\medskip
MSRI \\
Berkeley, CA  USA} \endaffil
\endtopmatter
\document  
\expandafter\ifx\csname bib4plain.tex\endcsname\relax
  \expandafter\gdef\csname bib4plain.tex\endcsname{}
\else \message{Hey!  Apparently you were trying to \string twice.   This does not make sense.}
\errmessage{Please edit your file (probably \jobname.tex) and remove
any duplicate ``\string\input'' lines} \fi

\def\renewcommand{\newcommand}	       
\edef\cite{\the\catcode`@}%
\catcode`@ = 11
\let\@oldatcatcode = \cite
\chardef\@letter = 11
\chardef\@other = 12
%
%
%
%
\def\@innerdef#1#2{\edef#1{\expandafter\noexpand\csname #2\endcsname}}%
%
%
\@innerdef\@innernewcount{newcount}%
\@innerdef\@innernewdimen{newdimen}%
\@innerdef\@innernewif{newif}%
\@innerdef\@innernewwrite{newwrite}%
%
%
%
\def\@gobble#1{}%
%
%
%
\ifx\inputlineno\@undefined
   \let\@linenumber = \empty 
\else
   \def\@linenumber{\the\inputlineno:\space}%
\fi
%
%
%
\def\@futurenonspacelet#1{\def\cs{#1}%
   \afterassignment\@stepone\let\@nexttoken=
}%
\begingroup 
\def\\{\global\let\@stoken= }%
\\ 
\endgroup
\def\@stepone{\expandafter\futurelet\cs\@steptwo}%
\def\@steptwo{\expandafter\ifx\cs\@stoken\let\@@next=\@stepthree
   \else\let\@@next=\@nexttoken\fi \@@next}%
\def\@stepthree{\afterassignment\@stepone\let\@@next= }%
%
%
%
\def\@getoptionalarg#1{%
   \let\@optionaltemp = #1%
   \let\@optionalnext = \relax
   \@futurenonspacelet\@optionalnext\@bracketcheck
}%
%
%
\def\@bracketcheck{%
   \ifx [\@optionalnext
      \expandafter\@@getoptionalarg
   \else
      \let\@optionalarg = \empty
      \expandafter\@optionaltemp
   \fi
}%
\def\@@getoptionalarg[#1]{%
   \def\@optionalarg{#1}%
   \@optionaltemp
}%
%
%
%
\def\@nnil{\@nil}%
\def\@fornoop#1\@@#2#3{}%
\def\@for#1:=#2\do#3{%
   \edef\@fortmp{#2}%
   \ifx\@fortmp\empty \else
      \expandafter\@forloop#2,\@nil,\@nil\@@#1{#3}%
   \fi
}%
\def\@forloop#1,#2,#3\@@#4#5{\def#4{#1}\ifx #4\@nnil \else
       #5\def#4{#2}\ifx #4\@nnil \else#5\@iforloop #3\@@#4{#5}\fi\fi
}%
\def\@iforloop#1,#2\@@#3#4{\def#3{#1}\ifx #3\@nnil
       \let\@nextwhile=\@fornoop \else
      #4\relax\let\@nextwhile=\@iforloop\fi\@nextwhile#2\@@#3{#4}%
}%
%
%
%
\@innernewif\if@fileexists
\def\@testfileexistence{\@getoptionalarg\@finishtestfileexistence}%
\def\@finishtestfileexistence#1{%
   \begingroup
      \def\extension{#1}%
      \immediate\openin0 =
         \ifx\@optionalarg\empty\jobname\else\@optionalarg\fi
         \ifx\extension\empty \else .#1\fi
         \space
      \ifeof 0
         \global\@fileexistsfalse
      \else
         \global\@fileexiststrue
      \fi
      \immediate\closein0
   \endgroup
}%
%
%
%
%
\def\bibliographystyle#1{%
   \@readauxfile
   \@writeaux{\string\bibstyle{#1}}%
}%
\let\bibstyle = \@gobble
%
%
\let\bblfilebasename = \jobname
\def\bibliography#1{%
   \@readauxfile
   \@writeaux{\string\bibdata{#1}}%
   \@testfileexistence[\bblfilebasename]{bbl}%
   \if@fileexists
      \nobreak
      \@readbblfile
   \fi
}%
\let\bibdata = \@gobble
%
%
\def\nocite#1{%
   \@readauxfile
   \@writeaux{\string\citation{#1}}%
}%
\@innernewif\if@notfirstcitation
%
%
\def\cite{\@getoptionalarg\@cite}%
%
%
\def\@cite#1{%
   \let\@citenotetext = \@optionalarg
   \printcitestart
   \nocite{#1}%
   \@notfirstcitationfalse
   \@for \@citation :=#1\do
   {%
      \expandafter\@onecitation\@citation\@@
   }%
   \ifx\empty\@citenotetext\else
      \printcitenote{\@citenotetext}%
   \fi
   \printcitefinish
}%
\def\@onecitation#1\@@{%
   \if@notfirstcitation
      \printbetweencitations
   \fi
   \expandafter \ifx \csname\@citelabel{#1}\endcsname \relax
      \if@citewarning
         \message{\@linenumber Undefined citation `#1'.}%
      \fi
      \expandafter\gdef\csname\@citelabel{#1}\endcsname{%
\strut
\vadjust{\vskip-\dp\strutbox
\vbox to 0pt{\vss\parindent0cm \leftskip=\hsize 
\advance\leftskip3mm
\advance\hsize 4cm\strut\openup-4pt 
\rightskip 0cm plus 1cm minus 0.5cm ?  #1 ?\strut}}
         {\tt
            \escapechar = -1
            \nobreak\hskip0pt
            \expandafter\string\csname#1\endcsname
            \nobreak\hskip0pt
         }%
      }%
   \fi
   \csname\@citelabel{#1}\endcsname
   \@notfirstcitationtrue
}%
%
%
\def\@citelabel#1{b@#1}%
%
%
\def\@citedef#1#2{\expandafter\gdef\csname\@citelabel{#1}\endcsname{#2}}%
%
%
%
\def\@readbblfile{%
   \ifx\@itemnum\@undefined
      \@innernewcount\@itemnum
   \fi
   \begingroup
      \def\begin##1##2{%
         \setbox0 = \hbox{\biblabelcontents{##2}}%
         \biblabelwidth = \wd0
      }%
      \def\end##1{}
      %
      %
      \@itemnum = 0
      \def\bibitem{\@getoptionalarg\@bibitem}%
      \def\@bibitem{%
         \ifx\@optionalarg\empty
            \expandafter\@numberedbibitem
         \else
            \expandafter\@alphabibitem
         \fi
      }%
      \def\@alphabibitem##1{%
         \expandafter \xdef\csname\@citelabel{##1}\endcsname {\@optionalarg}%
         \ifx\biblabelprecontents\@undefined
            \let\biblabelprecontents = \relax
         \fi
         \ifx\biblabelpostcontents\@undefined
            \let\biblabelpostcontents = \hss
         \fi
         \@finishbibitem{##1}%
      }%
      \def\@numberedbibitem##1{%
         \advance\@itemnum by 1
         \expandafter \xdef\csname\@citelabel{##1}\endcsname{\number\@itemnum}%
         \ifx\biblabelprecontents\@undefined
            \let\biblabelprecontents = \hss
         \fi
         \ifx\biblabelpostcontents\@undefined
            \let\biblabelpostcontents = \relax
         \fi
         \@finishbibitem{##1}%
      }%
      \def\@finishbibitem##1{%
         \biblabelprint{\csname\@citelabel{##1}\endcsname}%
         \@writeaux{\string\@citedef{##1}{\csname\@citelabel{##1}\endcsname}}%
         \ignorespaces
      }%
      %
      %
      \let\em = \bblem
      \let\newblock = \bblnewblock
      \let\sc = \bblsc
      \frenchspacing
      \clubpenalty = 4000 \widowpenalty = 4000
      \tolerance = 10000 \hfuzz = .5pt
      \everypar = {\hangindent = \biblabelwidth
                      \advance\hangindent by \biblabelextraspace}%
      \bblrm
      \parskip = 1.5ex plus .5ex minus .5ex
      \biblabelextraspace = .5em
      \bblhook
      \input \bblfilebasename.bbl
   \endgroup
}%
%
%
\@innernewdimen\biblabelwidth
\@innernewdimen\biblabelextraspace
%
%
%
\def\biblabelprint#1{%
   \noindent
   \hbox to \biblabelwidth{%
      \biblabelprecontents
      \biblabelcontents{#1}%
      \biblabelpostcontents
   }%
   \kern\biblabelextraspace
}%
%
%
%
\def\biblabelcontents#1{{\bblrm [#1]}}%
%
%
\def\bblrm{\rm}%
%
%
\def\bblem{\it}%
%
%
\def\bblsc{\ifx\@scfont\@undefined
              \font\@scfont = cmcsc10
           \fi
           \@scfont
}%
%
%
\def\bblnewblock{\hskip .11em plus .33em minus .07em }%
%
%
\let\bblhook = \empty
%
%
%
\def\printcitestart{[}
\def\printcitefinish{]}
\def\printbetweencitations{, }
\def\printcitenote#1{, #1}
%
%
%
\let\citation = \@gobble
%
%
%
\@innernewcount\@numparams
%
%
\def\newcommand#1{%
   \def\@commandname{#1}%
   \@getoptionalarg\@continuenewcommand
}%
%
%
\def\@continuenewcommand{%
   \@numparams = \ifx\@optionalarg\empty 0\else\@optionalarg \fi \relax
   \@newcommand
}%
%
%
\def\@newcommand#1{%
   \def\@startdef{\expandafter\edef\@commandname}%
   \ifnum\@numparams=0
      \let\@paramdef = \empty
   \else
      \ifnum\@numparams>9
         \errmessage{\the\@numparams\space is too many parameters}%
      \else
         \ifnum\@numparams<0
            \errmessage{\the\@numparams\space is too few parameters}%
         \else
            \edef\@paramdef{%
               \ifcase\@numparams
                  \empty  No arguments.
               \or ####1%
               \or ####1####2%
               \or ####1####2####3%
               \or ####1####2####3####4%
               \or ####1####2####3####4####5%
               \or ####1####2####3####4####5####6%
               \or ####1####2####3####4####5####6####7%
               \or ####1####2####3####4####5####6####7####8%
               \or ####1####2####3####4####5####6####7####8####9%
               \fi
            }%
         \fi
      \fi
   \fi
   \expandafter\@startdef\@paramdef{#1}%
}%
%
%
%
%
\def\@readauxfile{%
   \if@auxfiledone \else 
      \global\@auxfiledonetrue
      \@testfileexistence{aux}%
      \if@fileexists
         \begingroup
            \endlinechar = -1
            \catcode`@ = 11
            \input \jobname.aux
         \endgroup
      \else
         \message{\@undefinedmessage}%
         \global\@citewarningfalse
      \fi
      \immediate\openout\@auxfile = \jobname.aux
   \fi
}%
%
%
\newif\if@auxfiledone
\ifx\noauxfile\@undefined \else \@auxfiledonetrue\fi
%
%
%
%
\@innernewwrite\@auxfile
\def\@writeaux#1{\ifx\noauxfile\@undefined \write\@auxfile{#1}\fi}%
%
%
%
\ifx\@undefinedmessage\@undefined
   \def\@undefinedmessage{No .aux file; I won't give you warnings about
                          undefined citations.}%
\fi
%
%
\@innernewif\if@citewarning
\ifx\noauxfile\@undefined \@citewarningtrue\fi
%
%
%
\catcode`@ = \@oldatcatcode
 
\newpage

\head {Introduction} \endhead
\resetall
\bigskip

After the main gap theorem was proved (see \cite{Sh:c}), in discussion, 
Harrington expressed a desire for a finer structure - of finitary character 
(when we have a structure theorem at all).  I point out that the logic  
$L_{\infty,\aleph_0}$(d.q.)\,(d.q. stands for dimension quantifier) does not 
suffice:  suppose; e.g. for $T = Th \left( \lambda \times {}^\omega 2,E_n
\right)_{n < \omega}$ where $(\alpha,\eta)E_n(\beta,\nu) =: 
\eta \restriction n = \nu \restriction n$ and for $S \subseteq {}^\omega 2$
define $M_S = M \restriction \{(\alpha,\eta):[\eta \in S \Rightarrow \alpha < 
\omega_1]$ and $[\eta \in {}^\omega 2 \backslash S \Rightarrow \alpha <
\omega]\}$.  Hence, it seems to me we should try  
$L_{\infty,\aleph_\epsilon}$(d.q.)\,(essentially, in  ${\frak C}$ we can 
quantify over sets which are included in the algebraic closure of finite 
sets, see below \scite{1.1}, \scite{1.1A}), 
and Harrington accepts this interpretation.  Here the conjecture is
proved for $\aleph_\epsilon$-saturated models. 
\medskip

I.e., the main theorem is $M \equiv_{{\Cal L}_{\infty,\aleph_\epsilon}(d.q.)}
N \Leftrightarrow  M \cong N$ for $\aleph_\epsilon$-saturated 
models of a superstable countable (first order) theory $T$ without dop.  
For this we analyze further regular types, define a kind of infinitary 
logic (more exactly, a kind 
of type of $\bar a$ in $M)$, ``looking only up" in the definition (when 
thinking of the decomposition theorem).  To carry out the isomorphism 
proof we need: (\scite{1.8}) the type of the sum is the sum of types 
(infinitary types) assuming first order independence.  The main point to
construct an isomorphism between $M_1$ and $M_2$ when 
$M_1 \equiv_{{\Cal L}_{\infty,
\aleph_\epsilon}(d.q.)}M_2,Th(M_\ell) = T$ as above.  
So by \cite[X]{Sh:c} it is enough to construct isomorphic decompositions.
The construction of isomorphic decompositions is by $\omega$ approximations,
in stage  $n$, $\sim n$ levels of the decomposition tree are approximated,
i.e. we have $I^\ell_n \subseteq {}^{n \ge}$Ord and $\bar a^{n,\ell}_\eta 
\in M_\ell$ for $\eta \in I_n,\ell=1,2$ such that
tp$(\bar a^{n,1}_{\eta \restriction 0} \char 94 
\bar a^{n,1}_{\eta \restriction 1} \char 94 \ldots \char 94
\bar a^{n,1}_\eta,\emptyset,M) = \text{ tp}(\bar a^{n,2}_{\eta
\restriction 0} \char 94 \bar a^{n,2}_{\eta \restriction 1} \char 94 \ldots
\char 94 \bar a^{n,2}_\eta,\emptyset,M_2)$ with $\bar a^{\eta,\ell}_\eta$ 
being $\varepsilon$-finite, so in stage $n+1$, choosing $\bar a^{n+1,\ell}
_{\langle \rangle}$ we cannot take care of all types $\bar a^{n+1,\ell}
_{\langle \rangle} \char 94 \bar a^{\eta,\ell}_{\langle \alpha \rangle}$
so addition theorem takes care. 

In the end of \S1 (in \scite{1.12}) we point up that the addition 
theorem holds in fuller 
generalization.  In the second section we deal with a finer type needed for
shallow $T$, in the appendix we discuss how absolute is the isomorphism type. 

Subsequently Hrushovski and Bouscaren proved that even if  $T$  is 
$\aleph_0$-stable, ${\Cal L}_{\infty,\aleph_0}(d.q.)$ is not sufficient.  In
fact they found two examples: in one $T$ is (in addition) $\aleph_0$-stable 
of rank 3 and in the second it is not multidimensional.  Moreover, the 
examples show that even in the case of $\aleph_0$-stable theories,
$L_{\infty,\aleph_\varepsilon}$(d.q.) is not sufficient if one considers the
class of all models.  The first example is $\aleph_0$-stable shallow of
depth 3, and the second one is superstable (non $\aleph_0$-stable), NOTOP,
non-multidimensional. 

If we deal with $\aleph_\epsilon$-saturated models of shallow (superstable 
NDOP) theories $T$,  we can bound the depth of the quantification
$\gamma = DP(T)$; i.e. ${\Cal L}^\gamma_{\infty,\aleph_\epsilon}$ suffice. 

We assume the reader has a reasonable knowledge of \cite[V,\S1,\S2]{Sh:c} and 
mainly \cite[V,\S3]{Sh:c} and of \cite[X]{Sh:c}. 

Here is a slightly more detailed guide to the paper.  In \scite{1.1} we 
define the logic ${\Cal L}_{\infty,\aleph_\epsilon}$ and in \scite{1.1A} 
give a back and forth characterization of equivalence in this logic which 
is the operative definition for this paper. 

The major tools are defined in \scite{1.3}, \scite{1.5}.  In particular, 
the notion of $\text{tp}_\alpha$ defined in 1.5 is a kind of a depth 
$\alpha$ look-ahead which is actually used in the final construction.  
In \scite{1.9} we point out that equivalence 
in the logic ${\Cal L}_{\infty,\aleph_\epsilon}$ implies equivalence with 
respect to $\text{tp}_\alpha$ for all $\alpha$.  Proposition \scite{1.6} 
contains a number of important concrete assertions which are established 
by means of Facts 1.6A-F.  In general these explain the properties of 
decompositions over a pair $\binom BA$.  Claim \scite{1.8} (which follows 
from \scite{1.7}) is a key step in the final induction.  Definition 
\scite{1.11} establishes the framework for the 
proof that two $\aleph_\epsilon$-saturated structures that have the same  
$\text{tp}_\infty$ are isomorphic.  The induction step is carried out in 
\scite{1.11D}. 

I thank Baldwin for reading the typescript pointing needed corrections and 
writing down some explanations.
\bigskip

\noindent
\underbar{Notation}:  The notation is of \cite{Sh:c}, 
with the following additions (or reminders). 

If $\eta = \nu \char94 \langle \alpha \rangle$ then we let $\eta^- = \nu$; 
for $I$ is a set of sequences ordinals we let \newline
$\text{Suc}_T(\eta) = \{\nu:\text{for some }
\alpha,\nu = \eta \char 94 \langle \alpha \rangle \in I\}$. \newline
We work in  ${\frak C}^{eq}$ and for simplicity every first order formula 
is equivalent to a relation:  
\medskip
\roster
 
\item "{{}}"  $\perp$ means orthogonal (so $q$ is $\perp p$  means $q$ is 
orthogonal to $p$), \newline
remember $p \perp A$ means $p$ orthogonal to $A$; i.e. $p \perp q$ for every 
\newline

$\quad q \in S(ac\ell(A))$ (in ${\frak C}^{\text{eq}}$)
\item "{{}}"  $\perp_a$  means almost orthogonal
\item "{{}}"  $\perp_w$ means weakly orthogonal
\item "{{}}"  ${\frac{\bar a} B}$ and $\bar a / B$ means $\text{tp}(\bar a,B)$
\item "{{}}"  ${\frac{A} B}$ or $A/B$ means $\text{tp}_*(A,B)$
\item "{{}}"  $A + B$ means $A \cup B$
\item "{{}}"  $\nonfork{}{}_{A} \{B_i:i < \alpha \}$ means $\{B_i:i< \alpha\}$
is independent over $A$
\item "{{}}"  $\nonfork {A}{C}_{B}$ means $\{A,C\}$ is independent over $B$
\item "{{}}"  $\{C_i:i < \alpha \}$  is independent over $(B,A)$ means    
\newline
$j < \alpha \Rightarrow \text{ tp}_* \left( C_j,\dsize \bigcup_{i \ne j}
C_i \cup  B \right)$  does not fork over $A$
\item "{{}}"  regular type means stationary regular type $p \in S(A)$ for 
some $A$
\item "{{}}"  for $p \in S(A)$ regular and $C$ a set of elements realizing
$p$, dim$(C,p)$ is \newline

$\qquad$ Max$\{|\bold I|:\bold I \subseteq C$ is independent over
$A\}$
\item "{{}}"  $ac\ell(A) = \{c:\text{tp}(c,A)$ is algebraic$\}$
\item "{{}}"  $dc\ell(A) = \{c:\text{tp}(c,A)$ is realized by one and only one
element$\}$
\item "{{}}"  $Dp(p)$ is depth (of a stationary type, see 
\cite[X,Definition 4.3,p.528,Definition 4.4,p.529]{Sh:c}
\item "{{}}"  $Cb(p)$ is the canonical base of a stationary type $p$
(see \cite[III,6.10,p.134]{Sh:c})
\endroster
\newpage

\head {\S1 $\aleph_\epsilon$-saturated models} \endhead
\resetall
\bigskip

We first define our logic, but as said in \S0, we shall only use the condition
from \scite{1.1B}.
\definition{\stag{1.1} Definition}  1) The logic 
${\Cal L}_{\infty,\aleph_\epsilon}$ is slightly stronger than  
${\Cal L}_{\infty,\aleph_0}$, it consists of the set of formulas in  
${\Cal L}_{\infty,\aleph_1}$ such that any subformula of $\psi$ of the form
$(\exists \bar x)\varphi$ it actually has the form \newline
$(\exists \bar x^0,\bar x^1) \left[ \varphi_1(\bar x^1,\bar y) \and
\dsize \bigwedge_{i < \ell g \bar x^1} \left( \theta_i(x^1_i,\bar x^0) \and
(\exists^{< \aleph_0} z) \theta_i(z,\bar x^0) \right) \right],\bar x^0$
finite, $\bar x^1$ not necessarily finite; so $\varphi$ ``says" 
$\bar x^1 \subseteq ac\ell(\bar x^0)$. \newline
2)  ${\Cal L}_{\infty,\aleph_\epsilon}(d.q.)$  is like  
${\Cal L}_{\infty,\aleph_\epsilon}$ but we have cardinality quantifiers and 
moreover dimensional quantifiers (as in \cite[XIII,1.2,p.624]{Sh:c}). \newline
3)  The logic ${\Cal L}^\gamma_{\infty,\aleph_\epsilon}$ consist of the 
formulas of ${\Cal L}_{\infty,\aleph_\epsilon}$  such that $\varphi$  has 
quantifier depth  $< \gamma$ (but we start the inductive definition by
defining the quantifier depth of all first order as zero). \newline
4)  ${\Cal L}^\gamma_{\infty,\aleph_\epsilon}(d.q.)$ is like  
${\Cal L}^\gamma_{\infty,\aleph_\epsilon}$ but we have cardinality 
quantifiers and moreover dimensional quantifiers.  In fact the dimension 
quantifier is used in a very restricted way (see Definition \scite{1.4} and 
Claim \scite{1.9} + Claim \scite{1.11}).
\enddefinition
\bigskip

\demo{\stag{1.1A} Convention}  1) $T$ is a fixed first order complete theory,
${\frak C}$ is as in \cite{Sh:c}, so is $\bar \kappa$-saturated; ${\frak C}
^{\text{eq}}$ is as in \cite[III,6.2,p.131]{Sh:c}.  
We work in ${\frak C}^{\text{eq}}$ so
$M,N$ vary on elementary submodels of ${\frak C}^{\text{eq}}$ of cardinality
$< \bar \kappa$.  We assume $T$ is superstable with 
NDOP (countability is used only in the Proof of \scite{1.2} for bookkeepping).
\newline
Remember $a,b,c,d$ denote members of ${\frak C}^{\text{eq}},\bar a, \bar b,
\bar c,\bar d$ denote finite sequences of members of ${\frak C}^{\text{eq}},
A,B,C$ denote subsets of ${\frak C}^{\text{eq}}$ of cardinality $< \bar 
\kappa$. \newline
Remember $ac \ell(A)$  is the algebraic closure of $A$, i.e. \newline
$\{b:\text{for some first order and } n < \omega,\varphi(x,\bar y)
\text{ and } \bar a \subseteq A$ we have ${\frak C}^{\text{eq}} \models
\varphi[b,\bar a] \and (\exists^{\le n} y)\varphi(y,\bar a)\}$ and $\bar a$
denotes Rang$(\bar a)$ in places where it stands for a set (as in
$ac \ell(\bar a)$.  We write $\bar a \in A$ instead of $\bar a \in
{}^{\omega >}(A)$. \newline
2) $A$ is $\epsilon$-finite, if for some  $\bar a \in {}^{\omega >}A,A 
= ac \ell(\bar a)$.  (So for stable theories a subset of an 
$\epsilon$-finite set is not necessarily $\epsilon$-finite but as $T$ is 
superstable, a subset of an $\epsilon$-finite set is $\epsilon$-finite as if
$B \subseteq ac \ell(\bar a),\bar b \in B$ is such that $\text{tp}(\bar a,B)$
does not fork over $\bar b$, then trivially $ac \ell(\bar b) \subseteq A$
and if $ac\ell(\bar b) \ne B,\text{ tp}_*(B,\bar a \char 94 \bar b)$ forks
over $B$, hence \newline
(\cite[III,0.1]{Sh:c}) $\text{tp}(\bar a,B)$ forks over 
$\bar b$, a contradiction. \newline
So if $ac\ell(A) = ac\ell(B)$, then $A$ is $\epsilon$-finite iff $B$ is 
$\epsilon$-finite). \newline
3) When $T$ is superstable by \cite[IV,Table 1,p.169]{Sh:c} for $\bold F =
\bold F^a_{\aleph_0}$, all the axioms there hold and we write
$\aleph_\varepsilon$ instead of $\bold F$ and may use implicitly the
consequences in \cite[IV,\S3]{Sh:c}. 
\enddemo
\bigskip

\noindent
We may instead Definition \scite{1.1} use 
directly the standard characterization from \scite{1.1B}; 
as actually less is used we state the condition we shall actually use:
\proclaim{\stag{1.1B} Claim}  For models $M_1,M_2$ of $T$ we have $M_1 
\equiv_{{\Cal L}_{\infty,\aleph_\epsilon}(d.q.)}M_2$ \underbar{if} 
\medskip
\roster
\item "{$\bigotimes$}"  there is a family $F$ such that:
{\roster
\itemitem{ (a) }  each  $f \in  F$  is an $(M_1,M_2)$-elementary mapping,
(so Dom$(f) \subseteq M_1$, \newline
Rang$(f) \subseteq M_2$)
\itemitem{ (b) }  for  $f \in F$, Dom$(f)$ is $\epsilon$-finite (see 
\scite{1.1A}(2)) above)
\itemitem{ (c) }  if  $f \in F,\bar a_\ell \in M_\ell \, (\ell = 1,2)$ 
\underbar{then} for some $g \in F$ we have: \newline
$f \subseteq g$ and $ac\ell(\bar a_1) \in \text{ Dom}(f)$ and $ac\ell
(\bar a_2) \in \text{ Rang}(f)$
\itemitem{ (d) }  if  $f \cup \{ \langle a_1,a_2 \rangle\} \in F$ and tp$(a_1,
\text{Dom}(f))$ is stationary and regular \underbar{then}
dim$(\{a^1_1 \in M_1:f \cup \{ \langle a^1_1,a_2 \rangle \} \in F\},M_1)$    
\newline
$\qquad \quad = \text{ dim}(\{a^1_2 \in M_2:f \cup \{ \langle a_1,a^1_2
\rangle \} \in F\},M_2)$.
\endroster}
\endroster
\endproclaim
\bigskip

\noindent
Our main theorem is
\proclaim{\stag{1.2} Theorem}  Suppose $T$ is countable (superstable 
complete first order theory) with NDOP. 

Then
\medskip
\roster
\item   the ${\Cal L}_{\infty,\aleph_\epsilon}(d.q.)$ theory of an 
$\aleph_\epsilon$-saturated model characterizes it up to isomorphism. 
\item  Moreover, if $M_1,M_2$ are $\aleph_\epsilon$-saturated models of $T$
(so $M_\ell \prec {\frak C}^{\text{eq}}$) and $\bigotimes_{M_0,M_1}$ of
\scite{1.1B} holds, \underbar{then} $M_1,M_2$ are isomorphic.
\endroster
\medskip

\noindent
By \scite{1.1B}, it suffices to prove part (2).
\endproclaim
\bigskip

\centerline {$* \qquad * \qquad *$}
\bigskip

\noindent
The proof is broken into a series of claims (some of them do not use
NDOP).
\bigskip

\noindent
\underbar{Discussion}:  Let us motivate 
the notation and Definition below.  Suppose 
$N_0 \prec N_1 \prec N_2$ are $\aleph_\epsilon$-saturated and for $\ell = 
0,1,N_{\ell +1}$ is $\aleph_\epsilon$-prime over $N_\ell + b_\ell$ where  
$b_\ell /N_\ell$ is regular and $(b_1/N_1) \perp N_0$ (this arises naturally 
in the Decomposition theorem). 

We try to extract a finitary description we can choose.  If  $A_\ell  
\subseteq N_\ell$ is such that $b_\ell /N_\ell$ does not fork over  
$A_\ell,b_\ell /A_\ell$ stationary (and $A_\ell$ is $\epsilon$-finite,
of course) then to a large extent $A_0,b_0$ capture the essential part
of extending  $N_0$ to  $N_1$.  However, to see ``the successor" pair  
$(N_1,N_2)\,\,(\le_b$ below) we need to make $A_1$ explicit, for this we want 
to extend $(A_0,b_0)$ to $(A_0,b_0 + A_1)$, but $A_1$ was not arbitrary,  
$A_1 \subseteq N_1,N_1$ is $\aleph_\epsilon$-prime over $N_0 + b_0$, so the
reasonable way is to extend $A_0$ to some $\epsilon$-finite  $A^+_0 
\subseteq N_0$ such that $A_1/(A^+_0 + b_0)$ isolate the type
$A_1/(N_0 + b_0)$, this is $\le_a$ below.
\bigskip

\definition{\stag{1.3} Definition}  1)  $\Gamma = \{(A,B):A \subseteq B$  are 
$\epsilon$-finite$\}$.  Let \newline
$\Gamma(M) = \{(A,B) \in \Gamma:A \subseteq B \subseteq M\}$. \newline
2)  For members of  $\Gamma$  we write also  $\binom BA$, if $A \nsubseteq
B$ we mean $\binom {B \cup A}A$. \newline
3)  $\binom {B_1}{A_1} \le_a \binom {B_2}{A_2}$ 
(usually we omit $a$) if (both are
in $\Gamma$ and) \newline
$A_1 \subseteq A_2,B_1 \subseteq B_2,
\nonfork{B_1}{A_2}_{A_1}$ and ${\frac{B_2}{B_1+A_2}} \perp_a A_2$. \newline
4)  $\binom {B_1}{A_1} \le_b \binom {B_2}{A_2}$ if 
$A_2 = B_1,B_2 \backslash A_2 
= \bar b$ and ${\frac{\bar b}{A_2}}$ is regular orthogonal to $A_1$. \newline
5) $\le^*$ is the transitive closure of $\le_a \cup \le_b$.  (So 
it is a partial order, whereas in general  $\le_a \cup \le_b$ and $\le_b$ 
are not). \newline
6) We can replace $A,B$ by sequences listing them (we do not always strictly
distinguish).
\enddefinition
\bigskip

\definition{\stag{1.4} Definition}  1) We define $\text{tp}_\alpha 
[\binom BA,M]$ (for $A \subseteq B \subseteq M,A$ and $B$ are \newline
$\epsilon$-finite and $\alpha$ is an ordinal) and  
${\Cal S}_\alpha(\binom BA,M),{\Cal S}_\alpha(A,M)$ and 
${\Cal S}^r_\alpha(\binom BA,M),{\Cal S}^r_\alpha(A,M)$ 
by induction on $\alpha$ (simultaneously we mean; of course, we 
use appropriate variables):
\medskip
\roster
\item "{$(a)$}"  $\text{tp}_0[\binom BA,M]$ is the first 
order type of  $A \cup B$ 
\smallskip
\noindent
\item "{$(b)$}"   $\text{tp}_{\alpha +1}[\binom BA,M] = 
\text{ the triple } \langle Y^{1,\alpha}_{A,B,M},Y^{2,\alpha}_{A,B,M},
\text{tp}_\alpha(\binom BA,M)\rangle$ \newline
where: \newline
$Y^{1,\alpha}_{A,B,M} =: \biggl\{ \text{tp}_\alpha[\binom{B'}{A'},M]:
\text{for some } A',B' \text{ we have }
\binom BA \le_a \binom{B'}{A'} \in \Gamma (M) \biggr\}$, \newline
and $Y^{2,\alpha}_{A,B,M} =: \biggl\{
\langle \Upsilon,\lambda^\Upsilon_{M,B} \rangle:\Upsilon \in  
{\Cal S}^r_\alpha(B,M)\biggr\}$ \newline
where $\lambda^\Upsilon_{M,B} = 
\text{ dim} \left[ \biggl\{ d:\text{tp}_\alpha[\binom{B+d}B,M] = 
\Upsilon \biggr\},B \right]$:
\smallskip
\noindent
\item "{$(c)$}"  for  $\delta$  a limit ordinal,  
$\text{tp}_\delta[\binom BA,M]  = 
\langle \text{tp}_\alpha[\binom BA,M]:\alpha < \delta \rangle$
\newline
(this includes $\delta = \infty$, really  $\|M\|^+$ suffice).
\smallskip
\noindent
\item "{$(d)$}"  ${\Cal S}_\alpha(A,M) = 
\biggl\{ \text{tp}_\alpha[\binom BA,M]:\text{for some } B
\text{ such that } B \subseteq M$, \newline

$\qquad \qquad \qquad \qquad \qquad \qquad \text{and } 
\binom BA \in \Gamma(M) \biggr\}$
\smallskip
\noindent
\item "{$(e)$}"  ${\Cal S}^r_\alpha(\binom BA,M) = 
\biggl\{ \text{tp}_\alpha[\binom{B+c}{B},M]:\text{for some }
c \in M \text{ we have}$ \newline

$\qquad \qquad \qquad \qquad \qquad \qquad \qquad \frac cB \perp A$ 
and $\frac cB \text{ regular} \biggr\}$
\smallskip
\noindent
\item "{$(f)$}"  ${\Cal S}^r_\alpha(A,M) = 
\biggl\{ \text{tp}_\alpha[\binom{A+c}{A},M]:c \in M
\text{ and } \frac cA$ regular$\biggr\}$.
\endroster
\medskip

\noindent
2)  We define also  $\text{tp}_\alpha[A,M]$, for $A$ an $\epsilon$-finite 
subset of $M$:
\medskip
\roster
\item "{$(a)$}"   $\text{tp}_0[A,M] =$ first order type of $A$
\smallskip
\noindent
\item "{$(b)$}"  $\text{tp}_{\alpha +1}[A,M]$ is the triple
$\langle Y^{1,\alpha}_{A,M},Y^{2,\alpha}_{A,M},\text{tp}_\alpha[A,M] 
\rangle$ where \newline
$Y^{1,\alpha}_{A,M} =: {\Cal S}_\alpha(A;M)$ and \newline
$Y^{2,\alpha}_{A,M} =: \biggl\{ \langle \Upsilon,\text{dim} \{d \in M:
\text{tp}_\alpha[\binom{A+d}{A},M] = \Upsilon\}\rangle:
\Upsilon \in {\Cal S}^r_\alpha(A,M) \biggl\}$
\smallskip
\noindent
\item "{$(c)$}"   $\text{tp}_\delta[A,M] = \langle \text{tp}_\alpha(A,M):
\alpha < \delta \rangle$
\endroster
\medskip

\noindent
3)  $\text{tp}_\alpha[M] = \text{ tp}_\alpha[\emptyset,M]$.
\enddefinition
\bigskip

\definition{\stag{1.5} Definition}  1) $\langle N_\eta,a_\eta:\eta \in I 
\rangle$ is an $\aleph_\epsilon$-decomposition inside  $M$  above (or for or
over) $\binom BA$ (but we may omit the  $``\aleph_\epsilon-"$) if: 
\medskip
\roster
\item "{$(a)$}"  $I$ a set of finite sequences 
of ordinals closed under initial segments 
\smallskip
\noindent
\item "{$(b)$}"  $\langle \rangle,\langle 0 \rangle \in I,\eta \in I 
\backslash
\{ \langle \rangle \} \Rightarrow  \langle 0 \rangle \trianglelefteq \eta$,
let $I^- = I \backslash \{ \langle \rangle \}$, really $a_{\langle \rangle}$ 
is meaningless
\smallskip
\noindent
\item "{$(c)$}"  $A \subseteq N_{\langle \rangle},B \subseteq N_{\langle 0
\rangle},N_{\langle \rangle} \nonfork{}{}_{A} B$ and 
$dc \ell(\bar a_{\langle 0 \rangle}) = dc\ell(B)$,
\smallskip
\noindent
\item "{$(d)$}"  if $\nu = \eta \char 94 \langle \alpha \rangle \in I$
\underbar{then} $N_\nu$ is $\aleph_\epsilon$-primary over $N_\eta \cup 
\bar a_\nu,N_{\langle \rangle}$ is $\aleph_\epsilon$-prime \newline
over $A$
\smallskip
\noindent
\item "{$(e)$}"  for $\eta \in I$ such that $k = \ell g(\eta) > 1$ the type
$a_\eta / N_{\eta \restriction(k-1)}$ is orthogonal to 
$N_{\eta \restriction(k-2)}$
\smallskip
\noindent
\item "{$(f)$}"  $\eta \triangleleft \nu \Rightarrow N_\eta \prec N_\nu$ 
\smallskip
\noindent
\item "{$(g)$}"  $M$  is $\aleph_\epsilon$-saturated and $N_\eta \prec M$ for
$\eta \in I$
\smallskip
\noindent
\item "{$(h)$}"  if $\eta \in I \backslash \{ \langle \rangle\}$, then 
$\{a_\nu:\nu \in \text{ Suc}_I(\eta)\}$ is (a set of elements realizing 
over $N_\eta$ types orthogonal to $N_{\eta^-}$ and is) an independent set 
over $N_\eta$.
\endroster
\medskip

\noindent
2)  We replace ``inside  $M$"  by ``of  $M$"  if in addition
\medskip
\roster
\item "{$(i)$}"  in clause $(h)$ the set is maximal.
\endroster
\medskip 

\noindent
3)  $\langle N_\eta,a_\eta:\eta \in I \rangle$  is an 
$\aleph_\epsilon$-decomposition inside  $M$ if $(a)$, $(d)$, $(e)$, $(f)$, 
$(g)$, $(h)$ of part (1) holds and in clause $(h)$ we allow $\eta =
\langle \rangle$ (call this $(h)^+$). \newline
4)  $\langle N_\eta,a_\eta:\eta \in I \rangle$ is an 
$\aleph_\epsilon$-decomposition of $M$ if in addition to \scite{1.5}(3) and
we have the stronger version of clause $(i)$ of \scite{1.5}(2) by including
$\eta = \langle \rangle$, i.e. we have:
\medskip
\roster
\item "{$(i)^+$}"  for $\nu \in I$, the set 
$\{a_\eta:\eta \in \text{Suc}_I(\nu)\}$ is a maximal subset of $M$ 
independent over $N_\nu$.
\endroster
\medskip

\noindent
5) If $\langle N_\eta,a_\eta:\eta \in I \rangle$ is an 
$\aleph_\epsilon$-decomposition inside  $M$  we let  \newline
${\Cal P}(\langle N_\eta,a_\eta:\eta \in I \rangle,M) = \biggl\{ p \in 
S(M):p \text{ regular and for some } \eta \in I \backslash 
\{\langle \rangle\} \text{ we have}$ \newline        

$\qquad \qquad \qquad \qquad
\qquad \qquad \qquad \qquad p \text{ is orthogonal to }
N_{\eta^-} \text{ but not to }  N_\eta \biggr\}$.
\enddefinition
\bigskip

\definition{\stag{1.5A} Definition}  1) We say that 
$\langle N_\eta,a_\eta:\eta \in I \rangle$, an 
$\aleph_\epsilon$-decomposition inside $M$, is $J$-regular 
\underbar{if} $J \subseteq I$ and: 
\medskip
\roster
\item "{$(*)$}"  for each  $\eta \in I \backslash J$ there is $c_\eta$ 
such that $a_\eta \in ac\ell(N_\eta + c_\eta)$ \newline
${\frac{c_\eta}{N_\eta}}$ is regular and if $\eta \ne \langle \rangle$ then
${\frac{a_\eta}{N_\eta + c_\eta}} \perp_a N_{(\eta^-)}$.
\endroster
\medskip

\noindent
2)  We say  $``\langle N_\eta,a_\eta:\eta \in I \rangle$  is a regular 
$\aleph_\epsilon$-decomposition inside $M$ [of $M$]"  if it is an 
$\aleph_\epsilon$-decomposition inside $M$ [of $M$] which is 
$\emptyset$-regular. \newline
3)  We say  $``\langle N_\eta,a_\eta:\eta \in I \rangle$  is a regular 
$\aleph_\epsilon$-decomposition inside $M$ [of $M$] over  
$\binom BA$" \underbar{if} it is an $\aleph_\epsilon$-decomposition inside  
$M$  [of $M$] over $\binom BA$ which is $\{ \langle \rangle \}$-regular.
\enddefinition
\bigskip

\proclaim{\stag{1.5B} Claim}  1)  Every $\aleph_\epsilon$-saturated 
model has an $\aleph_\epsilon$-decomposition (i.e. of it). \newline
2)  If $M$ is $\aleph_\epsilon$-saturated, $\langle N_\eta,a_\eta:\eta  
\in I \rangle$ is an $\aleph_\epsilon$-decomposition inside $M$,  then for
some $J$, and $N_\eta,a_\eta$ for $\eta \in J \backslash I$  we have:  
$I \subseteq J$ and $\langle N_\eta,a_\eta:\eta \in J \rangle$ is an 
$\aleph_\varepsilon$-decomposition of $M$ (even an $I$-regular one). \newline
3)  If $M$ is $\aleph_\epsilon$-saturated, $\langle N_\eta,a_\eta:\eta  
\in I \rangle$ is an $\aleph_\epsilon$-decomposition of $M$ \underbar{then}
$M$ is $\aleph_\epsilon$-prime and $\aleph_\epsilon$-minimal over
$\dsize \bigcup_{\eta \in I} N_\eta$; if in addition $\langle N_\eta,a_\eta:
\eta \in \{ \langle \rangle,\langle 0 \rangle \} \rangle$ is an
$\aleph_\epsilon$-decomposition inside $M$ above $\binom BA$, then
$\langle N_\eta,a_\eta:\eta \in I \and (\eta \ne \langle \rangle \rightarrow
\langle 0 \rangle \trianglelefteq \eta) \rangle$ is an $\aleph_\epsilon$-
decomposition of $M$ above $\binom BA$. \newline
4)  If  $\langle N_\eta,a_\eta:\eta \in I \rangle$  is an 
$\aleph_\epsilon$-decomposition inside $M$ above $\binom BA$, \underbar{then}
it is an $\aleph_\epsilon$-decomposition inside $M$. \newline
5)  If $\langle N_\eta,a_\eta:\eta \in I \rangle$ is an 
$\aleph_\epsilon$-decomposition inside $M$ [above $\binom BA$], 
$\eta \in I$, \newline
$[\eta \in I \backslash \{ \langle \rangle \}],\alpha = \text{ Min}
\{\beta:\eta \char 94 
\langle \beta \rangle \notin I\},\nu =: \eta \char 94 \langle \alpha \rangle,
a_\nu \in M,{\frac{a_\nu}{N_\eta}}$  regular, orthogonal to  $\eta^-$ if
$\eta \ne \langle \rangle,N_\nu \prec M$ is 
$\aleph_\epsilon$-primary over $N_\eta + a_\nu$ and
$a_\nu \nonfork{}{}_{N_\eta} \left( \dsize \bigcup_{\rho \in I} N_\rho
\right)$ (enough to demand $\{a_\rho:\rho^- = 
\eta \text{ and } \rho \in I\}$
is independent over $N_\eta$) \underbar{then}   
$\langle N_\rho,a_\rho:\rho \in I \cup \{\nu\} \rangle$ is an 
$\aleph_\epsilon$-decomposition inside $M$ [over $\binom BA$]. \newline
6)  Assume $\langle N_\eta,a_\eta:\eta \in I \rangle$ is an
$\aleph_\epsilon$-decomposition of $M$, \underbar{if} $p$ is regular 
stationary and is not orthogonal to $M$ (e.g. $p \in S(M)$) \underbar{then}
for one and only one $\eta \in I$, there is a regular stationary 
$q \in S(N_\eta)$  not orthogonal to  $p$  such that: if $\eta^-$ is well 
defined (i.e. $\eta \ne \langle \rangle$), then $p \perp N_{\eta^-}$.
\newline
7)  Assume  $I = \dsize \bigcup_{\alpha < \alpha(*)} I_\alpha$, for each
$\alpha$ we have $\langle N_\eta,a_\eta:\eta \in I_\alpha \rangle$ is an 
$\aleph_\epsilon$-decomposition inside $M$ [above $\binom BA$] and 
for each  $\eta \in I$ for every $n < \omega$ and $\nu_\ell = \eta \char 94 
\langle \beta_\ell \rangle$ for $\ell < n$, for some $\alpha$ we have: 
$\{\nu_\ell:\ell < n\} \subseteq I_\alpha$  (e.g. $I_\alpha$ increasing).
\underbar{Then} $\langle N_\eta,a_\eta:\eta \in I \rangle$ is an 
$\aleph_\epsilon$-decomposition inside $M$ [above $\binom BA$]. \newline
8)  In (7), if $\eta \ne \langle \rangle$ and some $\nu_\ell$ is not 
$\triangleleft$-maximal in $I$ and ${\frac{a_{\nu_\ell}}{N_\eta}}$ 
is regular, it is enough:  

$$
\ell_1 < \ell_2 < n \Rightarrow  \dsize \bigvee_{\alpha < \alpha(*)}
[\{\nu_{\ell_1},\nu_{\ell_2}\} \subseteq I_\alpha].
$$
\medskip
\noindent
9)  If $\langle N_\eta,a_\eta:\eta \in I \rangle$  is an 
$\aleph_\epsilon$-decomposition inside  $M,I_1,I_2 \subseteq I$  are 
closed under initial segments and  $I_0 = I_1 \cap I_2$ \underbar{then}
$\left( \dsize \bigcup_{\eta \in I_1} N_\eta \right)
\nonfork{}{}_{\dsize \bigcup_{\eta \in I_0} N_\eta} \left( \dsize \bigcup
_{\eta \in I_2} N_\eta \right)$. \newline
10)  Assume that for $\ell = 1,2$ that $\langle N^\ell_\eta,a^\ell_\eta:\eta  
\in I \rangle$ is an $\aleph_\epsilon$-decomposition inside $M_\ell$, and  
for $\eta \in I$ the function $f_\eta$ is an isomorphism from $N^1_\eta$ 
onto $N^2_\eta$ and $\eta \triangleleft \nu \Rightarrow f_\eta \subseteq 
f_\nu$.  \underbar{Then}  $\dsize \bigcup_{\eta \in I} f_\eta$ is an 
elementary mapping; if in addition $\langle N^\ell_\eta,a^\ell_\eta:
\eta \in I \rangle$ is an
$\aleph_\epsilon$-decomposition of $M_\ell$ (for $\ell = 1,2)$ then
$\dsize \bigcup_{\eta \in I} f_\eta$ can be extended to an isomorphism from
$M_1$ onto $M_2$. \newline
11) If $\langle N_\eta,a_\eta:\eta \in I \rangle$ is an $\aleph_\epsilon$-
decomposition inside $M$ (above $\binom BA$) and $M^- \prec M$ is
$\aleph_\epsilon$-prime over $\dsize \bigcup_{\eta \in I} N_\eta$
\underbar{then} $\langle N_\eta,a_\eta:\eta \in I \rangle$ is an
$\aleph_\epsilon$-decomposition of $M$ \newline
(above $\binom BA$).
\endproclaim
\bigskip

\demo{Proof}  1), 2), 3), 5), 6), 9), 10).  
Repeat the proofs of \cite[X]{Sh:c}.
(Note that here $a_\eta /N_\eta$ is not necessarily regular, a minor change).
\newline
4), 7)  Check. \newline
8)  As  Dp$(p) > 0 \Rightarrow p$ is trivial, by \cite[ChX,7.2,p.551]{Sh:c}
and \cite[ChX,7.3]{Sh:c}. \hfill$\square_{\scite{1.5B}}$
\enddemo
\bigskip

\noindent
We shall prove:
\proclaim{\stag{1.6} Claim}  1) If $M$ is $\aleph_\epsilon$-saturated,  
$\binom BA \in \Gamma(M)$, \underbar{then} there is 
$\langle N_\eta,a_\eta:\eta \in I \rangle$, an $\aleph_\epsilon$-decomposition
of $M$ above $\binom BA$. \newline
2)  Moreover if $\langle N_\eta,a_\eta:\eta \in I \rangle$ satisfies clauses
$(a)-(h)$ of Definition \scite{1.5}(1), we can extend it to satisfy 
clause $(i)$ of \scite{1.5}(2), too.
\newline
3)  If  $\langle N_\eta,a_\eta:\eta \in I \rangle$ is an 
$\aleph_\epsilon$-decomposition of $M$ above $\binom BA$, $M^- \prec M$ is 
$\aleph_\epsilon$-prime over $\dsize \bigcup_{\eta \in I} N_\eta$ 
\underbar{then}: 
\medskip
\roster
\item "{$(a)$}"  $\langle N_\eta:\eta \in I \rangle$ is a 
$\aleph_\epsilon$-decomposition of $M^-$ 
\item "{$(b)$}"  we can find an $\aleph_\epsilon$-decomposition  
$\langle N_\eta,a_\eta:\eta \in J \rangle$ of $M$ such that $J \supseteq I$ 
and
$[\eta \in J \backslash I \Rightarrow (\eta \ne \langle \rangle \and \neg 
\langle 0 \rangle \triangleleft \eta)]$, moreover the last phrase follows from
the previous ones.
\endroster
\medskip

\noindent
4)  If in 3)(b) the set 
$J \backslash I$  is countable  (finite is enough for our 
applications), \underbar{then} necessarily $M,M^-$ are isomorphic, even 
adding all 
members of an $\epsilon$-finite subset of $M^-$ as individual constants.
\newline
5)  If $\langle N_\eta,a_\eta:\eta \in I \rangle$ is an 
$\aleph_\epsilon$-decomposition of $M$ above $\binom BA,I \subseteq J$ and
\newline
$\langle N_\eta,a_\eta:\eta \in J \rangle$ is an 
$\aleph_\epsilon$-decomposition of $M,M^- \prec M$ is
$\aleph_\epsilon$-prime over $\dsize \bigcup_{\eta \in I} N_\eta$ \newline
and
$\binom BA \le^* \binom{B_1}{A_1}$ and $B_1 \subseteq M$ and $c \in M$ and
${\frac{c}{B_1}} \perp A_1$ and ${\frac{c}{B_1}}$ is (stationary and) 
regular \underbar{then}
\medskip
\roster
\item "{$(\alpha)$}"  ${\frac{c}{B_1}} \perp 
{\frac{\cup\{N_\eta:\eta \in J \backslash I\}}{N_{\langle \rangle}}}$
\item "{$(\beta)$}"  ${\frac{c}{B_1}}$ is not orthogonal to some 
$p \in {\Cal P}(\langle N_\eta,a_\eta:\eta \in I \rangle,M^-)$.
\endroster
\medskip

\noindent
6) If $\langle N_\eta,a_\eta:\eta \in I \rangle$ is an 
$\aleph_\epsilon$-decomposition of $M$ above $\binom BA$ and $M^-$ is
$\aleph_\epsilon$-prime over $\dsize \bigcup_{\eta \in I} N_\eta$, 
\underbar{then}
the set ${\Cal P} = {\Cal P}(\langle N_\eta:\eta \in I \rangle,M)$
depends on $\binom BA$ and $M$ only \newline
(and not on $\langle N_\eta:\eta \in I \rangle$ or $M^-$), remember:

$$
\align
{\Cal P} = {\Cal P}(\langle N_\eta:\eta \in I \rangle,M) = \biggl\{ p \in 
S(M):&\,p \text{ regular and for some} \\
  &\,\eta \in I \backslash \{<>\},\text{ we have}: \\
  &\,p \text{ is orthogonal to } N_{\eta^-} \text{ but not to } N_\eta
\biggr\}.
\endalign
$$
\medskip

\noindent
So let ${\Cal P}(\binom BA,M) = {\Cal P}(\langle N_\eta:\eta \in I \rangle,
M)$. \newline
(7)  If $\frac BA$ is regular of depth zero or just $\frac bA \le_a
\frac BA,\frac bA$ regular of depth zero and $M$ is  
$\aleph_\epsilon$-saturated and $B \subseteq M$ \underbar{then} for any 
$\alpha$, we have  
tp$_\alpha(\binom BA,M)$ depend just on tp$_0(\binom BA,M)$. \newline
8)  For $\alpha < \beta$, from tp$_\beta(\binom BA,M)$ we can
compute tp$_\alpha(\binom BA,M)$. \newline
9)  If $f$ is an isomorphism from $M_1$ onto $M_2,A_1 \subseteq B_1$ are
$\varepsilon$-finite subsets of $M_1$ and $f(A_1) = A_2,f(B_1) = B_2$ then

$$
\text{tp}_\alpha(\binom{B_1}{A_1},M_1) = \text{ tp}_\alpha(\binom{B_2}{A_2},
M_2)
$$
\medskip

\noindent
(more pedantically tp$_\alpha(\binom{B_2}{A_2},M_2) = f[\text{tp}_\alpha
(\binom{B_1}{A_1},M_1)]$ or considered the $A_\ell,B_\ell$ as indexed
sets).
\endproclaim
\bigskip

\noindent
We delay the proof (parts (1), (2), (3) are proved after \scite{1.6E}, 
part (4), (6) after \scite{1.6F}, and parts (5), (7), (8) after \scite{1.6G}).
Part (9) is obvious.
\bigskip

\demo{\stag{1.5C} Definition}  1) If $\binom BA \in \Gamma(M),M$ is
$\aleph_\epsilon$-saturated let ${\Cal P}^M_{\binom BA}$ be the set 
${\Cal P}$ from Claim \scite{1.6}(6) above (by \scite{1.6}(6) this is 
well defined as we shall prove below). \newline
2) Let ${\Cal P}_{\binom BA} = \biggl\{ p:p \text{ is (stationary 
regular and) parallel to some } p' \in {\Cal P}^{{\frak C}^{\text{eq}}}
_{\binom BA} \biggr\}$.
\enddemo
\bigskip

\definition{\stag{1.6A} Definition}  If $\langle N^\ell_\eta,a_\eta:\eta \in 
J \rangle$ is a decomposition inside  ${\frak C}$ for $\ell = 1,2$  we say  
that $\langle N^1_\eta,a_\eta:\eta \in J \rangle \le^*_{\text{direct}} 
\langle N^2_\eta,a_\eta:\eta \in J \rangle$ if: 
\medskip
\roster
\item "{$(a)$}"  $N^1_{\langle \rangle} \prec N^2_{\langle \rangle}$ 
\item "{$(b)$}"  $N^2_{\langle \rangle} \nonfork{}{}_{N^1_{\langle \rangle}}
\{a_{\langle \alpha \rangle}:\langle \alpha \rangle \in J\}$ 
\item "{$(c)$}"  for  $\eta \in J \backslash \{\langle \rangle\},N^2_\eta$ is 
$\aleph_\epsilon$-prime over $N^1_\eta \cup N^2_{\eta^-}$.
\endroster
\enddefinition
\bigskip

\proclaim{\stag{1.6A1} Claim}  1) $M$ is $\aleph_\epsilon$-prime over $A$ iff 
$M$ is $\aleph_\epsilon$-primary over $A$ iff $M$ is 
$\aleph_\epsilon$-saturated, $A \subseteq M$, and for every $\bold I 
\subseteq M$ indiscernible over $A$ we have: $\text{dim}(\bold I,M) \le 
\aleph_0$ \underbar{iff} $M$ is $\aleph_\epsilon$-saturated and for every
finite $A \subseteq M$ and regular (stationary) $p \in S(A)$, we have
dim$(p,M) \le \aleph_0$. \newline
2) If $N_1,N_2$ are $\aleph_\epsilon$-prime over $A$, \underbar{then} they
are isomorphic over $A$.
\endproclaim
\bigskip

\demo{Proof}  By \cite[IV,4.18]{Sh:c} (see Definition \cite[IV,4.16]{Sh:c}, 
noting that part (4) disappears when we are speaking on 
$\bold F^a_{\aleph_0}$).
\enddemo
\bigskip

\demo{\stag{1.6B}Fact}  0) If $A$ is countable, 
$N$ is $\aleph_\epsilon$-primary 
over $A$ \underbar{then} $N$ is $\aleph_\epsilon$-primary \newline
over $\emptyset$. \newline
1) If  $N$  is $\aleph_\epsilon$-prime over $\emptyset,A$ countable,  
$N^+$ is $\aleph_\epsilon$-prime over $N \cup A$ \underbar{then} $N^+$ is 
$\aleph_\epsilon$-prime over $\emptyset$. \newline
2) If $\langle N_n:n < \omega \rangle$  is increasing, each $N_n$ is
$\aleph_\epsilon$-prime over $\emptyset$ or just 
$\aleph_\varepsilon$-constructible and $N_\omega$ 
is $\aleph_\epsilon$-prime over 
$\dsize \bigcup_{n < \omega} N_n$ \underbar{then} $N_\omega$ is 
$\aleph_\epsilon$-prime over $\emptyset$, (note $N_\omega = 
\dsize \bigcup_{n < \omega} N_n$). \newline
3)  If  $N_2 \nonfork{}{}_{N_0} N_1$, each $N_\ell$ is $\aleph_\epsilon$
-saturated, $N_2$ is $\aleph_\epsilon$-prime over $N_0 \cup \bar a$, and
$N_3$ is $\aleph_\epsilon$-prime over $N_2 \cup N_1$ \underbar{then}
$N_3$ is $\aleph_\epsilon$-prime over $N_1 \cup \bar a$. \newline
4) If $N_1 \prec N_2$ are $\aleph_\epsilon$-primary over $\emptyset$  
\underbar{then} for some $\aleph_\epsilon$-saturated $N_0 \prec N_1$ 
(necessarily $\aleph_\epsilon$-primary over $\emptyset$) we have: $N_1,N_2$ 
are isomorphic over $N_0$. \newline
5) In part (4), if $A \subseteq N_1$ is $\epsilon$-finite we can demand 
$A \subseteq N_0$. \newline
6)  If $M_0$ is $\aleph_\epsilon$-saturated, $\nonfork{A}{B}_{M_0},M_1$ is 
$\aleph_\epsilon$-primary over $M_0 \cup A$ \underbar{then}
$\nonfork{M_1}{B}_{M_0}$. \newline
7)  Assume $N_0 \prec N_1 \prec N_2$ are $\aleph_\epsilon$-saturated,  
$N_2$ is $\aleph_\epsilon$-primary over $N_1 + a$  and ${\frac{a}{N_1}} \perp 
N_0$ (and $a \notin N_1$).  \underbar{If} $N'_0 \prec N_0,N'_0 \prec
N'_1 \prec N_1,\nonfork{N'_1}{N_0}_{N'_0}$ and $N_1$ is 
$\aleph_\epsilon$-primary over $N_0 \cup N'_1$, \underbar{then} we can find 
$a',N'_2$  such that: $N'_2$ is 
$\aleph_\epsilon$-saturated, $\aleph_\epsilon$-primary over $N'_1 + a',N'_1 
\prec N'_2 \prec N_2,\nonfork{N_1}{N'_2}_{N'_1}$ and $N'_2$ is 
$\aleph_\epsilon$-primary over $N_1 \cup N'_2$. \newline
8) Assume $N'_0 \prec N_0 \prec N_1$ and $a \in N_1$ and $N_1$ is
$\aleph_\epsilon$-prime over $N_0 + a$ and ${\frac{a}{N_0}} \pm N'_0$
\underbar{then} we can find $a',N'_1$ such that $a' \in N',N'_0 \prec N'_1
\prec N_1,\nonfork{N'_1}{N_0}_{N'_0},N'_1$ is $\aleph_\epsilon$-prime over
$N'_0 + a$ and $N_1$ is $\aleph_\epsilon$-prime over $N_0 + N'_1$. \newline
9)  If $N_1$ is $\aleph_\epsilon$-saturated and $A \subseteq B \subseteq
N_1$ and $A,B$ are $\epsilon$-finite, \underbar{then} we can find
$N_0$ such that: $A \subseteq N_0 \prec N_1,N_0$ is $\aleph_\epsilon$-prime
over $\emptyset,A \subseteq N_0,\nonfork{B}{N_0}_{A}$, and 
$N_1$ is $\aleph_\epsilon$-prime over $N_0 \cup B$. \newline
10) If $N_0$ is $\aleph_\epsilon$-prime over $A$ and $A,B \subseteq N_0$ are
$\epsilon$-finite, \underbar{then} $N_0$ is $\aleph_\epsilon$-prime 
over $A \cup B$.
\enddemo
\bigskip

\remark{\stag{1.6B1} Remark}  In the proof of \scite{1.6B}(1)-(6),(10) we do
not use ``$T$ has NDOP".
\endremark
\bigskip

\demo{Proof}  0) There is $\{ a_\alpha:\alpha < \alpha^*\}$ a list of the
members of $N$ such that for $\alpha < \alpha(*)$ we have: tp$(a_\alpha,A
\cup \{a_\beta:\beta < \alpha\})$ is $\bold F^a_{\aleph_0}$-isolated.
\newline
[Why?  by the definition of ``$N$ is $\aleph_\epsilon$-primary over $A$").
Let $\{b_n:n < \omega\}$ list $A$ (if $A = \emptyset$ the conclusion is
trivial so without loss of generality $A \ne \emptyset$, hence we can find 
such a sequence $\langle b_n:n < \omega \rangle$).  Now define $\beta^* =
\omega + \beta$ and $b_{\omega + \alpha} = a_\alpha$ for $\alpha < \alpha^*$.
So $\{ b_\beta:\beta < \beta^* \}$ lists the elements of $N$ (remember
$A \subseteq N$ and check).  We claim that 
tp$(b_\beta,\{b_\gamma:\gamma < \beta\})$ is
$\bold F^a_{\aleph_0}$-isolated for $\beta < \beta^*$. \newline
\medskip

\noindent
(Why? if $\beta \ge \omega$, 
let $\beta' = \beta - \omega$ (so $\beta < \alpha^*$), now the statement 
above means \newline
tp$(a_{\beta'},A \cup \{ a_\gamma:\gamma < \beta
\})$ is $\bold F^a_{\aleph_0}$-isolated which we know; if $\beta < \omega$
this statement is trivial)].  
By the definition, $N$ is $\bold F^a_{\aleph_0}$-
primary over $\emptyset$. \newline
1) Note
\medskip
\roster
\item "{$(*)_1$}"  if $N$ is $\aleph_\epsilon$-primary over $\emptyset$ and
$A \subseteq N$ is finite then $N$ is $\aleph_\epsilon$-primary over $A$
\newline
[why?  see \cite[IV,3.12(3),p.180]{Sh:c} (of course, using 
\cite[IV,Table 1,p.169]{Sh:c} for $\bold F^a_{\aleph_0}$]
\item "{$(*)_2$}"  if $N$ is $\aleph_\epsilon$-primary over $\emptyset,A
\subseteq N$ is finite and $p \in S^m(M)$ does not fork over $A$ and $p
\restriction A$ is stationary then for some $\{ \bar a_\ell:\ell < \omega\}$
we have: $\bar a_\ell \in M$ realize $p,\{ \bar a_\ell:\ell < \omega\}$ is
independent over $A$ and $p \restriction (A \cup \dsize 
\bigcup_{\ell < \omega} \bar a_\ell) \vdash p$ \newline
[why?  \cite[IV,proof of 4.18]{Sh:c} (i.e. by it and \cite[4.9(3),4.11]{Sh:c}) 
or let $N'$ be $\aleph_\epsilon$-primary over
$A \cup \dsize \bigcup_{\ell < \omega} \bar a_\ell$ and note: $N'$ is
$\aleph_\epsilon$-primary over $A$
(proof like the one of \scite{1.6B}(0)) but also $N$ is 
$\aleph_\epsilon$-primary over $A$ so by uniqueness of the 
$\aleph_\epsilon$-primary model $N'$ is isomorphic to $N$ over $A$, so
without loss of generality $N' = N$; and easily $N'$ is as required].
\endroster
\medskip

Now we can prove \scite{1.6B}(1), for any $\bar c \in {}^{\omega >} A$, we
can find a finite $B^1_{\bar a} \subseteq N$ such that tp$(\bar c,N)$ does
not fork over $B^1_{\bar c}$, let $\bar b_{\bar c} \in {}^{\omega >}N$
realize stp$(\bar a,B^1_{\bar a})$ and let $B_{\bar c} = B^1_{\bar c} \cup
\bar b_{\bar c}$, so tp$(\bar c,N)$ does not fork over $B_{\bar c}$ and
tp$(\bar c,B_{\bar c})$ is stationary, hence we can find $\langle
\bar a^{\bar c}_\ell:\ell < \omega\}$ as in $(*)_2$ (for tp$(\bar c,
B_{\bar c})$).  Let \newline
$A' = \cup \{B_{\bar c}:\bar c \in {}^{\omega >} A\}
\cup \{ \bar a^{\bar c}_\ell:\bar c \in {}^{\omega >} A$ and $\ell <
\omega\}$, so $A'$ is a countable subset of $N$ and tp$_*(A,A') \vdash
\text{ tp}(A,N) = \text{ stp}(A,N)$.  
As $N$ is $\aleph_\epsilon$-primary over $\emptyset$ we
can find a sequence $\langle d_\alpha:\alpha < \alpha^* \rangle$ and
$\langle w_\alpha:\alpha < \alpha^* \rangle$ such that
$N = \{ d_\alpha:\alpha < \alpha^* \}$ and $w_\alpha \subseteq \alpha$ is
finite and stp$(d_\alpha,\{d_\beta:\beta \in w_\alpha\}) \vdash
\text{ stp}(d_\alpha,\{d_\beta:\beta < \alpha\})$ and $\beta < \alpha
\Rightarrow d_\beta \ne d_\alpha$.

We can find a countable set $W \subseteq \alpha^*$ such that $A' \subseteq
\{d_\alpha:\alpha \in W\}$ and $\alpha \in W \Rightarrow w_\alpha \subseteq
W$.  By \cite[IV,\S2,\S3]{Sh:c} without loss of generality $W$ is an initial
segment of $\alpha^*$.  \newline
Easily

$$
\alpha < \alpha^* \and \alpha \notin W \Rightarrow \text{ stp}(d_\alpha,
\{d_\beta:\beta \in W_\alpha) \vdash \text{ stp}(d_\alpha,A \cup \{d_\beta:
\beta < \alpha\}).
$$
\medskip

\noindent
As $N^+$ is $\aleph_\epsilon$-primary over $N \cup A$ we can find a list
$\{d_\alpha:\alpha \in [\alpha^*,\alpha^{**})\}$ of $N^+ \backslash (N \cup
A)$ such that tp$(d_\alpha,N \cup A \cup \{d_\beta:\beta \in [\alpha^*,
\alpha^{**})\})$ is $\aleph_\epsilon$-isolated.  So \newline
$\langle d_\alpha:\alpha \notin W,\alpha < \alpha^{**} \rangle$ 
exemplifies that $N^+$ is $\aleph_\epsilon$-primary over $A \cup A'$, 
hence by \scite{1.6B}(1) we know
that $N^+$ is $\aleph_\epsilon$-primary over $\emptyset$. \newline
2) By the characterization of ``$N$ is $\bold F^a_{\aleph_0}$-prime over
$A$" in \scite{1.6A1}). \newline

In more detail we use the last condition in \scite{1.6A1}(1).  Clearly
$N_\omega$ is $\aleph_\epsilon$-saturated (as it is $\aleph_\epsilon$-prime 
over $\dsize \bigcup_{n < \omega} N_n$.  Suppose $A \subseteq N_\omega$ 
is finite and $p \in S(A)$ is stationary and regular.
\bigskip

\noindent
\underbar{Case 1}:  $p$ not orthogonal to $\dsize \bigcup_{n < \omega} N_n$.

So for some $n < \omega,p$ is not orthogonal to $N_n$, hence 
there is a regular $p_1 \in S(N_n)$ such that $p,p_1$ are not 
orthogonal.  Let $A_1 \subseteq N_n$ be finite such that $p_1$ does not 
fork over $A$ and $p_1 \restriction A_1$ is stationary.  So by 
\cite[V,\S2]{Sh:c} we know dim$(p,N_\omega) = \text{ dim}
(p_1 \restriction A_1,N_\omega)$, hence it suffices to prove that the latter 
is $\aleph_0$.  Now use \cite[V,1.16(3),p.237]{Sh:c}.
\bigskip

\noindent
\underbar{Case 2}:  $p$ is orthogonal to $\dsize \bigcup_{n < \omega} N_n$.

Note that if each $N_n$ is $\aleph_\varepsilon$-prime then
$\dsize \bigcup_{n < \omega} N_n$ is $\aleph_\varepsilon$-saturated hence
$N = \dsize \bigcup_{n < \omega} N_n$ hence this case does not arise.

Let $B \subseteq \dsize \bigcup_{n < \omega}N_n$ be a finite set such that
tp$_*(A,\dsize \bigcup_{n < \omega}N_n)$ does not fork over $B$ and
stp$_*(A,B) \vdash \text{ stp}_*(A,\dsize \bigcup_{n < \omega}N_n)$.  Let
$q \in S(N_\omega)$ extend $p$ and does not fork over $A$ and let $q_1 = q
\restriction (B \cup A)$ and $q_2 = q \restriction (\dsize \bigcup_{n<\omega}
N_n \cup A)$.  Now by the assumption of our case $q_1$ is orthogonal to 
tp$_*(\dsize \bigcup_{n < \omega}N_m,B)$
hence (see \cite[V,\S3]{Sh:c}) $q_1 \Vdash q_2$ and let $\{a_\alpha:\alpha <
\alpha^*\} \subseteq q_1(N)$ be maximal set independent over
$\dsize \bigcup_{n<\omega}N_b \cup A$, so $|\alpha^*| \le \aleph_0$ hence
without loss of generality $\alpha^* \le \omega$ and $q \restriction
(\dsize \bigcup_{n < \omega} N_n \cup A \cup \{a_\alpha:\alpha < \alpha^*\})
\vdash q$.  Also clearly tp$_*(\{a_\alpha:\alpha < \alpha^*\},A \cup B) \vdash
\text{ tp}_*(\{a_\alpha:\alpha < \alpha^*\},A \cup \dsize \bigcup_{n<\omega}
N_n)$.
Together dim$(q_1,N) \le \aleph_0$ and clearly dim$(p,N) \le \aleph_0 +
\text{ dim}(q_1,N)$, so we are done.
\medskip

We can use a different proof, note:
\medskip
\roster
\item "{$\otimes_1$}"  if $\kappa = \text{ cf}(\kappa) \ge \kappa_r(T)$ and
$B_\alpha$ is $\bold F^a_\kappa$-constructible over $A$ for $\alpha <
\delta,\delta \le \kappa$ and $\alpha < \beta < \delta \Rightarrow B_\alpha
\subseteq B_\beta$ then $\dsize \bigcup_{\alpha < \delta} B_\alpha$ is
$\bold F^a_\kappa$-constructible over $A$ \newline
[why? see \cite[IV,\S3]{Sh:c}, \cite[IV,5.6,p.207]{Sh:c} for such arguments, 
assume ${\Cal A}_\alpha = \left< A,\langle a^\alpha_i:i < i_\alpha \rangle,
\langle B^\alpha_i:i < i_\alpha \rangle \right>$ is an 
$\bold F^a_\kappa$-construction.  Without loss of 
generality $i < j < i_\alpha \Rightarrow a^\alpha_i \ne a^\alpha_i$, and 
choose by induction on $\zeta,\langle u^\alpha_\zeta:\alpha < \delta 
\rangle$ such that: $u^\alpha_\zeta \subseteq i_\alpha,u^\alpha_\zeta$ 
increasing continuous in $i,u^\alpha_0 = \emptyset,|u^\alpha_{\zeta + 1} 
\backslash u^\alpha_\zeta| \le \kappa,u^\alpha_\zeta$ is 
${\Cal A}_\alpha$-closed and $\alpha < \beta < \delta$ implies 
$\{ a^\alpha_j:j \in u^\alpha_\zeta\} \subseteq \{ a^\beta
_j:j \in u^\beta_\zeta\}$ and tp$_*(\{a^\beta_i:i \in u^\beta_\zeta\},A \cup
\{a^\alpha_i:i < i_\alpha\})$ does not fork over $A \cup \{a^\alpha_i:i \in
u^\alpha_\zeta\}$].
\endroster
\medskip

\noindent
We use $\otimes_1$ for $\kappa = \aleph_0$. \newline
So each $N_n$ is $\aleph_\epsilon$-constructible over $\emptyset$ hence
$\dsize \bigcup_{n < \omega}N_n$ is $\aleph_\epsilon$-constructible over
$\emptyset$ and also $N_\omega$ is $\aleph_\epsilon$-constructible over
$\dsize \bigcup_{n < \omega}N_n$ hence $N_\omega$ is 
$\aleph_\epsilon$-constructible over $\emptyset$.  But $N_\omega$ is 
$\aleph_\epsilon$-saturated hence $N_\omega$ is
$\aleph_\epsilon$-primary over $\emptyset$. \newline
3) Suppose $N'_3$ is $\aleph_\epsilon$-saturated and $N_1 + \bar a
\subseteq N'_3$.  As $N_2$ is $\aleph_\epsilon$-prime over $N_0 + \bar a$ and
$N_0 + \bar a \subseteq N_1 + \bar a \subseteq N'_3$ we can find an
elementary embedding $f_0$ of $N_2$ into $N'_3$ extending 
$\text{id}_{N_0 + \bar a}$.  By \cite[V,3.3]{Sh:c}, the function 
$f_1 = f_0 \cup \text{id}_{N_1}$ is an elementary mapping so 
Dom$(f_1) = N_1 \cup N_2$.  As
$N_3$ is $\aleph_\epsilon$-prime over $N_1 \cup N_2$ and $f_1$ is an
elementary mapping from $N_1 \cup N_2$ into $N'_3$ which is an
$\aleph_\varepsilon$-saturated model there is an elementary
embedding $f_3$ of $N_3$ into $N'_3$ extending $f_2$.  So as for any such
$N'_3$ there is such $f_3$ clearly $N_3$ is $\aleph_\epsilon$-prime over
$N_1 + \bar a$, as required. \newline
4) Essentially by \cite{Sh:225}.

Let $N_0$ be $\aleph_0$-prime over $\emptyset$ and let 
$\{ p_i:i < \alpha\} \subseteq S(N_0)$ be a maximal family of 
pairwise orthogonal regular types.  Let
$\bold I_i = \{ \bar a^i_n:n < \omega\} \subseteq {\frak C}$ be an 
independent set of elements realizing $p_i$ and let 
$\bold I = \dsize \bigcup_{i < \alpha} \bold I_i$ and $N'_1$ be
$\bold F^a_{\aleph_0}$-prime over $N_0 \cup \bold I$.  Now
\medskip
\roster
\item "{$(*)$}"  if $\bar a,\bar b \subseteq N'_1$ and $\bar a/\bar b$ is
regular, then dim$(\bar a/\bar b,N'_1) \le \aleph_0$.
\endroster
\medskip

\noindent
[Why?  If $\bar a/\bar b \perp N_0$ then dim$(\bar a/\bar b,N'_1) \le
\aleph_0$ by \cite[1.2]{Sh:225}.  If $\bar b/\bar a \pm N_0$, 
then for some $\bar b' \char 94 \bar a' \subseteq N_0$ realizing stp$(\bar b 
\char 94 \bar a,\emptyset)$, we have $\bar a'/\bar b' \pm \bar a/b$ hence
dim$(\bar a/b,N'_1) = \text{dim}(\bar a'/\bar b',N'_1)$, so without loss
of generality $\bar b \char 94 \bar a \subseteq N_0$, similarly without loss
of generality there is $i(*) < \alpha$ such that $\bar b/\bar a \subseteq
p_{i(*)}$ and easily dim$(\bar b/\bar a,N'_1) = \text{dim}(\bar b/\bar a,
N_0) + \text{ dim}(p_{i(*)},N_0) \le \aleph_0 + \aleph_0 = \aleph_0$ (see
\cite[V,1.6(3)]{Sh:c}).  So we have proved $(*)$].
\medskip

\noindent
Now use \scite{1.6A1}(1) to deduce: $N'_1$ is 
$\bold F^a_{\aleph_\epsilon}$-prime over $\emptyset$ hence (by uniqueness 
of $\aleph_\epsilon$-prime model, \scite{1.6A1}(2)) $N'_1 \cong N_1$.

By renaming without a loss of generality $N'_1 = N_1$.  Now
\medskip
\roster
\item "{$(**)$}"  $(N_1,c)_{c \in N_0},(N_2,c)_{c \in N_0}$ are
$\aleph_\epsilon$-saturated and for if: $\bar a \in {\frak C},\bar b \in 
N_\ell,\bar a/\bar b$ a regular type and 
$\nonfork{\bar a}{(N_0 + \bar b)}_{\bar b}$ (for $\ell = 1$ or $\ell = 2$), 
then \newline
dim$(\bar a/(\bar b \cup N_0),N_\ell) = \aleph_0$.
\endroster
\medskip

\noindent
[Why?  Remember that we work in $({\frak C}^{\text{eq}},c)_{c \in N_0}$.  The
``$\aleph_\epsilon$-saturated" follows from the second statement. \newline
Note: dim$(\bar a/(\bar b \cup N_0),N_\ell) \le \text{ dim}(\bar a/\bar b,
N_\ell) \le \aleph_0$ (first inequality by monotonicity, second inequality
by \scite{1.6A1}(1) and the assumption ``$N_\ell$ is $\aleph_\epsilon$-prime 
over $\emptyset$").  If $\bar a/\bar b$ is not orthogonal to $N_0$ then for 
some $i < \alpha$ we have $p_i \pm (\bar a/\bar b)$ so easily (using
``$N_\ell$ is $\aleph_\varepsilon$-saturated") we have 
dim$(\bar a/(\bar b \cup N_0),N_0) = \text{ dim}(p_i,N_\ell) \ge \| 
\bold I_i\| = \aleph_0$; so together with the
previous sentence we get equality.  Lastly, if $\bar a/\bar b \perp N_0$ 
by \cite[1.1,p.280]{Sh:225}, we have
dim$(\bar a/(\bar b \cup N_0),N_\ell) <
\aleph_0 \Rightarrow \text{ dim}(\bar a/\bar b,N_\ell) < \aleph_0$ which
contradicts the assumption ``$N_\ell$ is $\aleph_\epsilon$-saturated".
So we have proved $(**)$ hence by \scite{1.6A1}(1) we get ``$N_1,N_2$ are 
isomorphic over $N_0$ as required. \newline
5)  This is proved similarly. \newline
6)  By \cite[V,3.2]{Sh:c}. \newline
7)  As  $N_1$ is $\aleph_\epsilon$-prime over $N_0 \cup N'_1$ and as $T$  
has NDOP (i.e. does not have DOP) we know (by \cite[X,2.1,2.2,p.512]{Sh:c})
that $N_1$ is $\aleph_\epsilon$-minimal over $N_0 \cup N'_1$ and 
${\frac{a}{N_1}}$ is not orthogonal to $N_0$ or to $N'_1$.  But
$a/N_1 \perp N_0$ by an assumption, so $a/N_1$ is not orthogonal to $N'_1$
hence there is a regular  $p' \in S(N'_1)$  not orthogonal to
${\frac{a}{N_1}}$  hence (by \cite[V,1.12,p.236]{Sh:c}) $p'$ is realized 
say by  $a' \in N_2$.  By \cite[V,3.3]{Sh:c}, we know that $N_2$ is 
$\aleph_\epsilon$-prime over $N_1 + a'$.  We can find $N'_2$ which is
$\aleph_\epsilon$-prime over $N'_1 + a'$  such that $N_2$ is 
$\aleph_\epsilon$-prime over $N_1 \cup N'_2$ (by uniqueness, 
i.e. \scite{1.6A1}(1)), so we have finished. \newline
8)  Similar easier proof.  \newline
9)  Let $N'_0$ be $\aleph_\epsilon$-prime over $A$ such that
$\nonfork{B}{N'_0}_{A}$, and let $N'_1$ be $\aleph_\epsilon$-prime over
$N'_0 \cup B$.  By \scite{1.6B}(1), we know that $N'_1$ is 
$\aleph_\epsilon$-prime over $\emptyset$, and by \scite{1.6B}(10) below 
$N'_1$ is
$\aleph_\epsilon$-prime over $A \cup B$, hence by \scite{1.6A1}(2) we know
that $N'_1,N_1$ are isomorphic over $A \cup B$ hence without loss of
generality $N'_1 = N_1$ and so $N_0 = N'_0$ is as required. \newline
10) By \cite[IV,3.12(3),p.180]{Sh:c}.   \hfill$\square_{\scite{1.6B}}$
\enddemo
\bigskip

\demo{\stag{1.6C} Fact}  Assume $\langle N^1_\eta,a_\eta:\eta \in I \rangle 
\le^*_{\text{direct}} \langle N^2_\eta,a_\eta:\eta \in I \rangle$ (see
Definition \scite{1.6A}) and $A \subseteq B \subseteq N^1_{<0>}$ and
$\dsize \bigwedge_{\eta \in I} N^2_\eta \prec M$.
\medskip
\roster
\item  If $\nu = \eta \char94 \langle \alpha \rangle \in I$, then
$\nonfork{N^2_\eta}{N^1_\nu}_{N^1_\eta}$ and even $N^2_\eta
\nonfork{}{}_{N^1_\eta} \, \left(\dsize \bigcup \Sb \rho \in I \\
\neg \eta \triangleleft \rho \endSb N^1_\rho \right)$; and \newline
$\eta \triangleleft \nu \in I$ implies $N^2_\nu \nonfork{}{}_{N^1_\eta} \,
\left( \dsize \bigcup \Sb \rho \in T \\
\neg \eta \triangleleft \rho \endSb N^1_\rho \right)$.
\item   $\langle N^2_\eta,a_\eta:\eta \in I \rangle$ is an 
$\aleph_\epsilon$-decomposition inside $M$ above $\binom BA$ \underbar{iff}
\newline 
$\langle N^1_\eta,a_\eta:\eta \in I \rangle$ is an 
$\aleph_\epsilon$-decomposition inside $M$ above $\binom BA$.
\item   Similarly replacing $``\aleph_\epsilon$-decomposition inside 
$M$ above $\binom BA$" by \newline
$``\aleph_\epsilon$-decomposition of $M$ above $\binom BA$".
\endroster
\enddemo
\bigskip

\demo{Proof}  1)  We prove the first statement by induction on $\ell g(\eta)$.
If $\eta = \,<>$ this is clause (b) by the Definition \scite{1.6A} and clause
(d) of Definition \scite{1.5}(1) (and \cite[V,3.2]{Sh:c}).
If $\eta \ne <>$, then ${\frac{a_\nu}{N_\eta}} \perp N^1_{(\eta^-)}$ (by 
condition (e) of Definition \scite{1.5}(1)).  By the induction hypothesis  
$N^2_{(\eta^-)} \nonfork{}{}_{N^1_{(\eta^-)}} N^1_\eta$ and we know $N^2_\eta$
is $\aleph_\epsilon$-primary over $N^2_{(\eta^-)} \cup N^1_\eta$; we 
know this implies no $p \in S(N^1_\eta)$ orthogonal to $N^1_{\eta^-}$ is
realized in $N^2_\eta$ hence ${\frac{a_\nu}{N^1_\eta}} \perp 
{\frac{N^2_\eta}{N^1_\eta}}$, so ${\frac{a_\nu}{N^1_\eta}} \vdash 
{\frac{a_\nu}{N^2_\eta}}$  hence  ${\frac{N^1_\nu}{N^1_\eta}} \perp 
{\frac{N^2_\eta}{N^1_\eta}}$ hence 
$\nonfork{N^1_\nu}{N^2_\eta}_{N^1_\eta}$ as 
required.  The other statements hold by the non-forking calculus (remember
if $\eta = \nu \char 94 \langle \alpha \rangle \in I$ then use
tp$(\cup\{N^1_\rho:\eta \trianglelefteq \rho \in I\},N^1_\eta)$ is orthogonal
to $N^1_\nu$ or see details in the proof of \scite{1.6D}(1)$(\alpha)$).
\newline
2)  By Definition \scite{1.6A}, for $\ell = 1,2$ we have: $\langle
N^\ell_\eta,a_\eta:\eta \in I \rangle$ is a decomposition inside ${\frak C}$
and by assumption $\dsize \bigwedge_{\eta \in I} N^1_\eta \prec N^2_\eta \prec
M$.  So for $\ell = 1,2$ we have to prove $``\langle N^\ell_\eta,a_\eta:\eta
\in I \rangle$ is an $\aleph_\epsilon$-decomposition inside $M$ for
$\binom BA$" assuming this holds for $1-\ell$.  We have to check Definition
\scite{1.5}(1). \newline
Clauses 1.5(1)(a),(b) for $\ell$ holds because they hold for $1 - \ell$.
\newline
Clause 1.5(1)(c) holds as by the assumptions $A \subseteq B \subseteq N^1
_{<0>} \prec N^2_{<0>}$ and $\nonfork{N^1_{<0>}}{N^2_{<>}}_{N^1_{<>}}$.
\newline
Clauses 1.5(1)(d),(e),(f),(h) holds as $\langle N^\ell_\eta,a_\eta:\eta \in
I \rangle$ is a decomposition inside ${\frak C}$ (for $\ell = 1$ given, for
$\ell = 2$ easily checked). \newline
Clause 1.5(1)(g) holds as $\dsize \bigwedge_\eta N^1_\eta \prec N^2_\eta
\prec M$ is given and $M$ is $\aleph_\epsilon$-saturated. \newline
3)  First we do the ``only if" direction; i.e. prove the maximality of  
$\langle N^1_\eta,a_\eta:\eta \in I \rangle$  as an 
$\aleph_\epsilon$-decomposition inside $M$ for $\binom BA$ (i.e. 
condition (i) from \scite{1.5}(2)), assuming it holds for  
$\langle N^2_\eta,a_\eta:\eta \in I \rangle$.  If this fails, then for some
$\eta \in I \backslash \{<>\}$ and \newline
$a \in M,\{a_{\eta \char94 <\alpha >}:
\eta \char94 <\alpha > \in  I\} \cup  \{a\}$  is 
independent over  $N^1_\eta$ and \newline
$a \notin \{a_{\eta \char94 \langle \alpha \rangle}:
\eta \char94 \langle \alpha \rangle \in I\}$ and 
${\frac{a}{N_\eta}} \perp N^1_{\eta^-}$.
Hence, if  $\eta \char94 \langle \alpha_\ell \rangle \in I$ for 
$\ell < k$ then
$\bar a = \langle a \rangle \char 94 \langle a_{\eta \char94 <\ell>}
:\ell < k \rangle$  realizes over $N^1_\eta$ a type orthogonal to  
$N^1_{\eta^-}$, but \newline
$N^1_{\eta^-} \prec N^1_\eta,N^1_{\eta^-} \prec N^2_{\eta^-}$ and 
$\nonfork{N^1_\eta}{N^2_\eta}_{N^1_{\eta^-}}$ (see \scite{1.6C}(1), hence 
\cite[V,2.8]{Sh:c}) \newline
tp$(\bar a,N^2_\eta) \perp N^2_{\eta^-}$ hence $\{a\} \cup 
\{a_{\eta \char 94 \langle \ell \rangle}:\ell < k\}$ is independent over 
$N^2_\eta$ but $k,\eta \char 94 \langle \alpha_\ell \rangle$ were arbitrary 
so $\{a\} \cup \{a_{\eta \char 94 \langle \alpha 
\rangle}:\eta \char 94 \langle \alpha \rangle \in I\}$ is independent over
$N^2_\eta$ contradicting condition (i) 
from Definition \scite{1.5}(2) for $\langle N^2_\eta,a_\eta:\eta \in I 
\rangle$. \newline
For the other direction use: if the conclusion fails, then for some 
$\eta \in I \backslash \{<>\}$ and \newline
$a \in M$ the set $\{a_{\eta \char 94 \langle \alpha \rangle}:
\eta \char 94 \langle \alpha \rangle \in I\} \cup \{a\}$ is independent
of $N^1_\eta$. \hfill$\square_{\scite{1.6C}}$
\enddemo
\bigskip

\demo{\stag{1.6D} Fact}  Assume $\langle N^1_\eta,a^1_\eta:\eta \in I \rangle$
is an $\aleph_\epsilon$-decomposition inside $M$. \newline
1)  If $N^1_{<>} \prec N^2_{<>} \prec  M,N^2_\eta$ is 
$\aleph_\epsilon$-prime over $\emptyset$ and \newline
$N^2_{<>} \nonfork{}{}_{N^1_{<>}} \{a^1_{<\alpha >}:<\alpha > \in I\}$
\underbar{then}
\medskip
\roster
\item "{$(\alpha)$}"  $\left[ N^2_{<>} \nonfork{}{}_{N^1_{<>}} \,\,
\dsize \bigcup_{\eta \in I} N^1_\eta \right]$ and 
\smallskip
\noindent
\item "{$(\beta)$}"  we can find $N^2_\eta \, (\eta \in I \backslash \{<>\})$
such that $N^2_\eta \prec M$, and \newline
$\langle N^1_\eta,a^1_\eta:\eta \in I \rangle \le^*_{\text{direct}} \langle 
N^2_\eta,a^1_\eta:\eta \in I \rangle$.
\endroster
\medskip 

\noindent
2)  If  $Cb\,{\frac{a^1_{<\alpha >}}{N^1_{<>}}} \subseteq N^0_{<>} 
\prec N^1_{<>}$ or at least $N^0_{<>} \prec N^1_{<>}$ and
${\frac{a^1_{<\alpha>}}{N^1_{<>}}} \pm N^0_{<>}$
whenever $<\alpha > \in I$ \underbar{then} we can find $N^0_\eta \prec M$ 
and  $a^0_\eta \in N_\eta$ (for $\eta \in I\backslash \{<>\})$
such that  $\langle N^0_\eta,a^0_\eta:\eta \in I \rangle  
\le^*_{\text{direct}} \langle N^1_\eta,a^1_\eta:\eta \in I\rangle$.
\enddemo
\bigskip

\demo{Proof}  1) For proving $(\alpha)$ let $\{\eta_i:i < i^*\}$ list the set
$I$ such that $\eta_i \triangleleft \eta_j \Rightarrow i < j$, so 
$\eta_0 = <>$ and without loss of generality for some $\alpha^*$ we have
\newline
$\eta_i \in\{<\alpha>:<\alpha> \in I\} \Leftrightarrow i \in [1,\alpha^*)$.  
Now we prove by induction on $\beta \in [1,\alpha^*)$ \newline
that $N^2_{<>} \nonfork{}{}_{N^1_{<>}} \cup \{N^1_{\eta_i}:i < \beta\}$.
For $\beta = 1$ this is assumed.  For $\beta$ limit use the local character
of non-forking.

If $\beta = \gamma +1 \in [1,\alpha^*)$, then by repeated use of
\cite[V,3.2]{Sh:c} (as $\{a_{\eta_0}:j \in [1,\beta)\}$ is independent over
$(N^1_{<>},N^2_{<>})$ and $N^1_{<>}$ is $\aleph_\epsilon$-saturated and
$N'_{\eta_j}(j \in [1,\gamma))$ is $\aleph_\epsilon$-prime over 
$N^1_{<>} + a_{\eta_j}$) we know that tp$(a_{\eta_\gamma},N^2_{<>} \cup 
\dsize \bigcup_{i < \gamma} N^1_{\eta_i})$ does not fork over $N^1_{<>}$.  
Again by
\cite[V,3.2]{Sh:c}, tp$_*(N^1_{\eta_\gamma},N^2_{<>} \cup \dsize \bigcup
_{i < \gamma} N^1_{\eta_i})$ does not fork over $N^1_{<>}$ hence
$\dsize \bigcup_{i < \beta} N^1_{\eta_i} \nonfork{}{}_{N^1_{<>}} N^2_{<>}$
and use symmetry.
\medskip

Lastly, if $\beta \in \gamma +1 \in [\alpha^*,i^*),
\text{tp}(a_{\eta^-_\gamma},N_{\eta_\gamma})$ is orthogonal to 
$N^1_{<>}$ and even to $N^1_{(\eta^-_\gamma)^-}$ so again by non-forking and
\cite[V,3.2]{Sh:c} we can do it, so clause $(\alpha)$ holds.

For clause $(\beta)$, we choose $N^2_{\eta_i}$ for $i \in [1,i^*)$ by
induction on $i < i^*$ such that $N^2_{\eta_i} \prec M$ is
$\aleph_\epsilon$-prime over $N^2_{\eta^-_i} \cup N^1_{\eta_i}$.  By the
non-forking calculus we can check Definition \scite{1.3}. \newline
2)  We let $\{\eta_i:i < i^*\}$ be as above, now we choose by induction on
$i \in [1,i^*) \, N^0_{\eta_i},a^0_{\eta_i}$ such that:
\medskip
\roster
\item "{$(*)$}"  $N^0_{\eta_i} \prec N^1_{\eta_i}$ and $\dsize 
\bigcup_{N^0_{\eta^-_i}} N^1_{\eta_i}$ and $N^1_{\eta_i}$ is
$\aleph_\epsilon$-prime over $N^0_{\eta_i} \cup N^1_{\eta^-_i}$
\item "{$(**)$}"  $a^0_{\eta_i} \in N^0_{\eta_i}$ and $N^0_{\eta_i}$ is
$\aleph_\epsilon$-prime over $N^0_{\eta^-_i} \pm a^0_{\eta_i}$.
\endroster
\medskip

\noindent
The induction step is done: if $\ell g(\eta_i) > 1$ by \scite{1.6B}(7) and if
$\ell g(\eta_i) = 1$ by \newline
\scite{1.6B}(8). \hfill$\square_{\scite{1.6D}}$
\enddemo
\bigskip

\demo{\stag{1.6E} Fact}  1)  \underbar{If} $\langle N^1_\eta,a_\eta:\eta \in I
\rangle \le^*_{\text{direct}} \langle N^2_\eta,a_\eta:\eta \in I\rangle$ both 
$\aleph_\epsilon$-decompositions inside $M$ \underbar{then}

$$
{\Cal P}(\langle N^1_\eta,a^1_\eta:\eta  \in  I\rangle,M) = 
{\Cal P}(\langle N^2_\eta,a^2_\eta:\eta  \in  I\rangle,M).
$$
\enddemo
\bigskip

\demo{Proof}  By Defintion \scite{1.5}(5) it suffices to prove, for each
$\eta \in I \backslash \{<>\}$ that
\medskip
\roster
\item "{$(*)$}"  for regular $p \in S(M)$ we have \newline
$p \perp N^1_{\eta^-} \and
p \pm N^1_\eta \Leftrightarrow p \pm N^2_{\eta^-} \and p \pm N^2_\eta$.
\endroster
\medskip

\noindent
Now for regular $p \in S(M)$: first assume $p \perp N^1_{\eta^-} \and
p \pm N^1_\eta$ so $p \pm N^2_\eta$ (as $N^1_\eta \prec N^2_\eta$ and
$p \pm N^1_\eta$) and we can find a regular $q \in S(N^1_\eta)$ such that
$q \pm p$; so $q \perp N^1_{\eta^-}$, now $q \perp N^2_{\eta^-}$ (as
$N^1_\eta \nonfork{}{}_{N^1_{\eta^-}} N^2_{\eta^-}$ and $q \perp N^1_\eta$
see \cite[V,2.8]{Sh:c}) hence $p \perp N^2_{\eta^-}$.

Second, assume $p \perp N^2_{\eta^-} \and p \pm N^2_\eta$, remember
$N^1_{\eta^-},N^1_\eta,N^2_\eta,N^3_\eta$ are \newline
$\aleph_\epsilon$-saturated, 
$N^1_\eta \nonfork{}{}_{N^1_{\eta^-}} N^2_{\eta^-}$ and $N^2_\eta$ is
$\aleph_\epsilon$-prime over $N^1_\eta \cup N^2_{\eta^-}$ and $T$ does not
have DOP.  Clearly $N^2_\eta$ is $\aleph_\epsilon$-minimal over $N^1_\eta
\cup N^2_{\eta^-}$ and every regular $q \in S(N^2_\eta)$ is not orthogonal
to $N^1_\eta$ or to $N^2_{\eta^-}$.  Also as $p \pm N^2_\eta$ there is a
regular $q \in S(N^2_\eta)$ not orthogonal to $p$, so as $p \perp
N^2_{\eta^-}$ also $q \perp N^2_{\eta^-}$; hence by the previous sentence 
$q \pm N^1_\eta$ hence $p \pm N^1_\eta$.  Lastly, as $p \perp 
N^2_{\eta^-}$ and $N^1_{\eta^-} \prec N^2_{\eta^-}$ clearly $p \perp 
N^1_{\eta^-}$, as required.  \hfill$\square_{\scite{1.6E}}$
\enddemo
\bigskip

\noindent
At last we start proving \scite{1.6}.
\demo{Proof of \scite{1.6}}  1) Let $N^0 \prec {\frak C}$ be 
$\aleph_\epsilon$-primary over $A$, without loss of generality  
$\nonfork{N^0}{B}_{A}$
(but not necessarily $N^0 \prec M$),  and let  $N^1$ be 
$\aleph_\epsilon$-primary over $N^0 \cup B$.  Now by \scite{1.6B}(0) the model
$N^0$ is $\aleph_\epsilon$-primary over $\emptyset$ and by \scite{1.6B}(1) 
the model $N^1$ is $\aleph_\epsilon$-primary over $\emptyset$ hence (by 
\scite{1.6B}(10)) is
$\aleph_\epsilon$-primary over  $B$, hence without loss of generality $N^1 
\prec M$.  Let  $N_{<>} =: N^0,N_{<0>} = N^1,I = \{<>,<0>\}$ and  
$a_{<0>} = B$.  More exactly $a_\eta$ is such that dcl$(\{a_\eta\}) = 
\text{dcl}(B)$.  Clearly $\langle N_\eta,a_\eta:\eta \in I \rangle$ is an
$\aleph_\epsilon$-decomposition inside  $M$  above  $\binom BA$.  
Now apply part (2) of \scite{1.6} proved below. \newline
2) By \scite{1.5B}(4) we know $\langle N_\eta,a_\eta:\eta \in I \rangle$ is an
$\aleph_\epsilon$-decomposition inside $M$, by \scite{1.6B}(2) we then find
$J \supseteq I$ and $N_\eta,a_\eta$ for $\eta \in J \backslash I$ such that
$\langle N_\eta,a_\eta:\eta \in I \rangle$ is an 
$\aleph_\epsilon$-decomposition of $M$.  By \scite{1.6B}(3), 
$\langle N_\eta,a_\eta:\eta \in J'_1 \rangle$ is an 
$\aleph_\epsilon$-decomposition of $M$ for $\binom BA$ where
$J' =: \{ \eta \in J:\eta = <> \text{ or } \langle 0 \rangle \trianglelefteq
\eta \in J\}$. \newline
3) Part (a) holds by \scite{1.5B}(2),(3).  As for part (b) by \scite{1.5B}(2) 
there is  $\langle N_\eta,a_\eta:\eta  \in  J \rangle$,  an 
$\aleph_\epsilon$-decomposition of  $M$  with  $I \subseteq  J$;  easily  
$[\langle 0 \rangle \trianglelefteq \eta \in J \Rightarrow \eta \in I]$.
\newline
${}{}$  \hfill $\square_{\scite{1.6}(1),(2),(3)}$
\enddemo
\bigskip

\demo{\stag{1.6F} Fact}  If $\langle N^\ell_\eta,a^\ell_\eta:\eta \in 
I^\ell \rangle$ are $\aleph_\epsilon$-decompositions of  $M$  for  
$\binom BA$,  for  $\ell  = 1,2$  and  $N^1_{<>} = N^2_{<>}$ \underbar{then}  
${\Cal P}(\langle N^1_\eta,a^1_\eta:\eta \in I^1 \rangle,M) = 
{\Cal P}(\langle N^2_\eta,a^2_\eta:\eta \in I^2 \rangle,M)$.
\enddemo
\bigskip

\demo{Proof}  By \scite{1.6}(3)(b) we can find $J^1 \supseteq I^1$ and 
$N^1_\eta,a^1_\eta$ for $\eta \in J^1 \backslash I^1$ such that  
$\langle N^1_\eta,a^1_\eta:\eta \in J^1\rangle$ is an 
$\aleph_\epsilon$-decomposition of $M$.  So $\eta \in J^1\backslash I^1 
\Leftrightarrow \eta
\ne \langle \rangle \and \neg (\langle 0 \rangle \triangleleft \eta)$.  Let
$J^2 = I^2 \cup  (J^1\backslash I^1)$  and for  $\eta  
\in  J^2\backslash I^2$ let  $a^2_\eta =: a^1_\eta $,  $N^2_\eta  =: 
N^1_\eta $.  Easily  $\langle N^2_\eta,a^2_\eta:\eta  \in  J^2\rangle $  is 
an $\aleph _\epsilon $-decomposition of  $M$.  By \scite{1.5B}(6) we know 
that for every regular $p \in S(M)$ there is (for $\ell = 1,2)$ a unique  
$\eta(p,\ell) \in  J^\ell$ such that $p \perp N_{\eta(p,\ell)} \and p 
\perp N_{\eta(p,\ell)^-}$ (note $\langle \rangle^-$ -- meaningless).  
By the uniqueness of $\eta(p,\ell)$, if $\eta(p,1) \in  J^1\backslash I^1$ 
then as it can serve as $\eta(p,2)$ clearly it is $\eta(p,2)$ so 
$\eta(p,2) = \eta(p,1) \in  
J^1\backslash I^1 = J^2\backslash I^2$;  similarly  $\eta(p,2) \in J^2 
\backslash I^2 \Rightarrow  \eta (p,1) \in  J^1\backslash I^1$ and
$\eta(p,1) = \langle \rangle \Leftrightarrow \eta(p,2) = \langle \rangle$.  
So 
\medskip
\roster
\item "{$(*)$}"   $\eta (p,1) \in  I^1 \backslash \{ \langle \rangle \}
\Leftrightarrow \eta(p,2) \in I^2 \backslash \{ \langle \rangle\}$.
\endroster
\medskip

\noindent
But 
\roster
\item "{$(**)$}"  $\eta(p,\ell) \in I^\ell \backslash \{ \langle \rangle \}
\Leftrightarrow  p \in {\Cal P}(\langle N^\ell_\eta,a^\ell_\eta:\eta \in 
I^\ell \rangle,M)$.
\endroster
\medskip

\noindent
Together we finish. \hfill$\square_{\scite{1.6F}}$
\enddemo
\bigskip

\noindent
We continue proving \scite{1.6}.
\demo{Proof of \scite{1.6}(4)}  Let $N^1_\eta = N_\eta$ and 
$a^1_\eta = a_\eta$ for $\eta \in I$ and without loss of generality 
$J \ne I$ hence $J \backslash I \ne \emptyset$. 

Let  $N^2_{<>} \prec  M$  be $\aleph_\epsilon$-prime over $\dsize
\bigcup_{\nu \in J \backslash I} N_\nu$; letting $J \backslash I = \{
\eta_i:i < i^\ast \}$  be such that $[\eta_i \triangleleft \eta_j \Rightarrow
i < j]$  we can find $N^2_{<>,i}$ (for $i \le i^*)$ increasing continuous, 
$N^2_{<>,0} = N_{<>}$ and $N^2_{<>,i+1}$ is $\aleph_\epsilon$-prime over  
$N^2_{<>,i} \cup N_{\eta_i}$ hence over $N^2_{<>,i} \cup a_{\eta_i}$.
Lastly, let $N^2_{<>,i^*} = N^2_{<>}$. 

By \scite{1.6B}(1),(2) we know $N^2_{<>}$ is $\aleph_\epsilon$-primary over  
$\emptyset$ and (using repeatedly \scite{1.6B}(6) + finite character of 
forking) we have $N^2_{<>} \nonfork{}{}_{N^1_{<>}} a_{<0>}$.  
By \scite{1.6B}(4) \newline
(with 
$N^1_{<>},N^2_{<>},Cb\,(a_{<>}/N^1_{<>})$  here standing for $N_1,N_2,A$
there) we can find a model $N^0_{<>}$ such that $a_{<0>} 
\nonfork{}{}_{N^0_{<>}} N^1_{<>}$ and 
$Cb\,(a_{<>}/N^1_{<>}) \subseteq N^0_{<>},
N^0_{<>} \prec  N^1_{<>},N^0_{<>}$ is $\aleph_\epsilon$-primary over
$\emptyset$  and  $N^1_{<>},N^2_{<>}$ are isomorphic over $N^0_{<>}$.  
By \scite{1.6D}(1) we can choose  $N^2_\eta \prec M$  
with  $N^1_\eta \prec N^2_\eta$ and $\langle N^1_\eta,a^1_\eta:\eta \in 
I\rangle \le^*_{\text{direct}} \langle N^2_\eta ,a^1_\eta :\eta  \in  
I\rangle$.  Similarly, by \scite{1.6D}(2) (here Suc$_I(<>) = 
\{ \langle 0 \rangle\}$) we can choose an 
$\aleph_\epsilon$-decomposition $\langle N^0_\eta,a^0_\eta:\eta \in  
I \rangle$ with $\langle N^0_\eta,a^0_\eta:\eta \in I \rangle 
\le^*_{\text{direct}} \langle N^1_\eta,a^1_\eta:\eta \in I \rangle$.  By 
\scite{1.6}(3) we know that $\langle N^1_\eta,a^1_\eta:\eta \in I \rangle$ is 
an $\aleph_\epsilon$-decomposition of $M^-$ and easily  
$\langle N^2_\eta,a^1_\eta:\eta \in I \rangle$  is an 
$\aleph_\epsilon$-decomposition of $M$.  Now choose by induction on  $\eta  
\in  I$  an isomorphism  $f_\eta$ from  $N^1_\eta $ onto  $N^2_\eta $ over  
$N^0_\eta $ such that  $\nu \triangleleft \eta \Rightarrow f_\nu  \subseteq  
f_\eta$.  For  $\eta = <>$ we have chosen  $N^0_\eta $ such that  
$N^1_\eta,N^2_\eta$ are isomorphic over  $N^0_\eta $.  For the induction 
step note that  $f_{(\eta^-)} \cup \text{ id}_{N^0_\eta}$ is an elementary 
mapping as $N^2_{(\eta^-)} \nonfork{}{}_{N^0_{(\eta^-)}} N^0_\eta$ and 
$f_{(\eta^-)} \cup \text{ id}_{N^0_\eta}$ can 
be extended to an isomorphism  $f_\eta$ from $N^1_\eta$ onto $N^2_\eta$ 
as  $N^\ell_\eta$ is $\aleph_\epsilon$-primary (in fact even 
$\aleph_\epsilon $-minimal) over  $N^\ell_{(\eta^-)} \cup  N^0_\eta $ 
for $\ell = 1,2$  (which holds easily).  Now by \scite{1.5B}(3) the model 
$M^-$ is $\aleph_\epsilon$-saturated and $\aleph_\epsilon$-primary and 
$\aleph_\epsilon$-minimal over  $\dsize \bigcup_{\eta \in J} N_\eta = 
\dsize \bigcup_{\eta \in I} N^1_\eta$; similarly $M$ is $\aleph_\epsilon$
-primary over $\dsize \bigcup_{\eta \in I} N^2_\eta$.  Now $\dsize 
\bigcup_\eta f_\eta$ is an elementary mapping from  $\dsize \bigcup_{\eta \in
I} N^1_\eta$ onto $\dsize \bigcup_{\eta \in I} N^2_\eta$ hence 
can be extended to an isomorphism $f$ from $M^-$ into $M$.  By 
the $\aleph _\epsilon $-minimality of $M$ over $\dsize \bigcup_{\eta \in I}
N_\eta$ (see \scite{1.5B}(3)), $f$ is onto $M$, so $f$ is as required.        
\hfill$\square_{\scite{1.6}(4)}$
\enddemo
\bigskip

\noindent
We delay the proof of \scite{1.6}(5).
\demo{Proof of \scite{1.6}(6)}  Let  $\langle N^\ell_\eta,a_\eta :\eta  \in  
I^\ell \rangle $  for  $\ell = 1,2$, be $\aleph_\epsilon$-decompositions of
$M$ above $\binom BA$, so dcl$(a^\ell_{<>}) = \text{ dcl}(B)$.  Let $p \in 
S(M)$, and assume that $p \in {\Cal P}(\langle N^1_\eta,a^1_\eta :\eta \in  
I^1\rangle,M)$, i.e. for some  $\eta \in I^1\backslash \{<>\}$,  
$(p_\eta \perp N_{\eta^-}) \and p_\eta \pm N_\eta$.  We shall prove that the
situation is similar for $\ell = 2$; i.e. $p \in {\Cal P}(\langle N^2_\eta,
a^2_\eta:\eta \in I^2 \rangle,M)$; 
by symmetry this suffices. 

Let  $n = \ell g(\eta)$,  choose  $\langle B_\ell:\ell \le n\rangle$  
such that:
\medskip
\roster
\item "{$(\alpha)$}"  $A \subseteq B_0$,
\item "{$(\beta)$}"  $B \subseteq B_1$,
\item "{$(\gamma)$}"  $a_{\eta \restriction \ell} \subseteq B_\ell \subseteq  
N^1_{\eta \restriction \ell}$, for  $\ell \le n$
\item  "{$(\delta)$}"  $B_{\ell +1} \nonfork{}{}_{B_\ell} 
N^1_{\eta \restriction \ell}$,
\item "{$(\epsilon)$}"  ${\frac{B_{\ell +1}}{B_\ell + a^1_{\eta \restriction
(\ell +1)}}} \vdash
{\frac{B_{\ell +1}}{N^1_{\eta \restriction \ell} + a^1_{\eta \restriction
(\ell +1)}}}$,
\item "{$(\zeta)$}"  for some  $d \in B_n,{\frac{d}{B_n\backslash \{d\}}}$
is regular $\pm p$, (hence $\perp B_{n-1}$)
\item "{$(\eta)$}"  $B_\ell$ is $\epsilon$-finite. 
\endroster
\medskip

[Why?  $\langle B_\ell:\ell \le n \rangle$ exists?  Prove by induction on  
$n$ that for any $\eta \in I$ of length $n$ and $\epsilon$-finite
$B' \subseteq N_\eta$, there is $\langle B_\ell:\ell \le n\rangle$  
satisfying $(\alpha)-(\epsilon)$, $(\eta)$ such that $B' \subseteq B_n$.  
Now there is  $p' \in  S(N^1_\eta)$  regular not orthogonal to $p$,  let  
$B^1 \subseteq  N^1_\eta $ be an $\epsilon$-finite set extending $Cb(p')$.  
Applying the previous sentence to  $\eta$, $B^1$ we get $\langle B_\ell:
\ell \le n\rangle$,  let  $d \in  N_\eta $ realize $p' \restriction B_n$. 

Now if $n > 0$, tp$(d,B_n) \perp N_{\eta^-}$ hence tp$(d,B_n) \perp B_{n-1}$,
hence \newline
tp$(d,B_n) \perp tp_*(N_{\eta^-},B_n)$, hence as
tp$(d,B_n)$ is stationary by \cite[V,1.2(2),p.231]{Sh:c}, the types
tp$(d,B_n)$,tp$_*(N_{\eta^-},B_n)$ are weakly orthogonal so \newline
tp$(d,B_n) \vdash \text{ tp}(d,N_{\eta^-} \cup B_n)$ hence
${\frac{B_n + d}{B_{n+1}+a^1_\eta}} \vdash 
{\frac{B_n+d}{N^1_{\eta^-} + a^1_\eta}}$.

Now replace $B_n$ by $B_n \cup \{d\}$ and we finish].

Note that necessarily 
\medskip
\roster
\item "{$(\delta)^+$}"   $B_n \nonfork{}{}_{B_m} N^1_{\eta \restriction m}$ 
for $m \le n$. \newline
[Why?  By the nonforking calculus].
\medskip
\noindent
\item "{$(\epsilon)^+$}"  ${\frac{B_n}{B_m+a^1_{\eta \restriction(m+1)}}} 
\perp_a B_m$ for $m < n$. \newline
[Why?  As $N^1_{\eta \restriction m}$ is $\aleph_\epsilon$-saturated].
\endroster
\medskip 

\noindent
Choose  $D^* \subseteq  N^2_{<>}$ finite such that 
${\frac{B_n}{N^2_{<>}+B}}$ does not fork over $D^* + B$. \newline
[Note: we really mean  $D^* \subseteq N^2_{<>}$,  not $D^* \subseteq
N^1_{<>}$].

We can find $N^3_{<>},\aleph_\varepsilon$-prime over $\emptyset$ such 
that $A \subseteq N^3_{<>} \prec N^2_{<>}$ and
$D^* \nonfork{}{}_{A} N^3_{<>}$ and $N^2_{<>}$ is 
$\aleph_\varepsilon$-prime over $N^3_{<>} \cup D^*$ (by \scite{1.6B}(9)).  
Hence $B_n \nonfork{}{}_{B}N^3_{<>}$ (by non-forking calculus).  
As tp$_*(B,N^2_{<0>})$ does not fork over $A \subseteq N^3_{<>} \subseteq
N^2_{<>}$ by \scite{1.6D}(2) we can find 
$N^3_\eta,a^3_\eta$ (for $\eta \in I^2 \backslash \{<>\})$, such that 
$\langle N^3_\eta,a^3_\eta:\eta \in I^2 \rangle 
\le^*_{\text{direct}} \langle N^2_\eta,a^2_\eta:\eta \in I^2 \rangle$ and  
$a^3_{<0>} = a^2_{<0>}$ (remember dcl$(a^2_{<0>}) =
\text{ dcl}(B)$).  By \scite{1.6C}(2) we know $\langle N^3_\eta,a^3_\eta:
\eta \in I^2 \rangle$ is an $\aleph_\epsilon$-decomposition of $M$ for 
$\binom BA$. \newline
By \scite{1.6E} it is enough to show  $p \in {\Cal P}(\langle N^3_\eta,
a^3_\eta:\eta \in I^2 \rangle,M)$.  
Let $N^4_{<>} \prec N^2_{<>} \prec M$ be $\aleph_\epsilon$-prime over 
$N^3_{<>} \cup B_0$.  Now by the non-forking calculus $B \nonfork{}{}_{A}
(N^3_{<>} \cup B_0)$ [why? because
\medskip
\roster
\item "{$(a)$}"   as said above $B_n \nonfork{}{}_{B} N^3_{<>}$ but
$B_0 \subseteq B_n$ so $B_0 \nonfork{}{}_{B} N^3_{<>}$, and
\smallskip
\noindent
\item "{$(b)$}"   as $B \nonfork{}{}_{A} N^1_{<>}$ and $B_0 \subseteq
N^1_{<>}$ we have $\nonfork{B}{B_0}_{A}$ so $\nonfork{B_0}{B}_{A}$
\endroster
\medskip

\noindent
hence (by (a) + (b) as $A \subseteq B$)
\medskip
\roster
\item "{$(c)$}"   ${\frac{B_0}{N^3_{<>}+B}}$ does not fork over $A$, also
\medskip
\noindent
\item "{$(d)$}"   $B \nonfork{}{}_{A} N^3_{<>}$ (as $A \subseteq N^3_{<>} 
\subseteq N^2_{<>}$ and tp$(B,N^3_{<>})$ does not fork over $A$) \newline
putting (c) and (d) together we get
\medskip
\noindent
\item "{$(e)$}"  $\nonfork{}{}_{A} \{B_0,B,N^3_{<>}\}$ \newline
hence the conclusion].
\endroster
\medskip

\noindent
Hence
$B \nonfork{}{}_{N^3_{<>}} B_0$ so $B \nonfork{}{}_{N^3_{<>}} N^4_{<>}$ 
(by \scite{1.6B}(6)) and so (as  $A \subseteq  N^3_{<>} \subseteq  
N^4_{<>})$  we have  $N^4_{<>} \nonfork{}{}_{N^3_{<>}} N^3_{<0>}$ and 
by \scite{1.6D}(1) we can define  $N^4_\eta \prec M$ (for $\eta  
\in I^2 \backslash \{<>\})$, such that $\langle N^4_\eta,a^3_\eta:\eta  
\in I^2 \rangle \ge_{\text{direct}} \langle N^3_\eta,a^3_\eta:\eta \in  
I^2 \rangle$.  So by \scite{1.6C}(1) \newline
$\langle N^4_\eta,a^3_\eta:\eta \in I^2 \rangle$ is an 
$\aleph_\epsilon$-decomposition of $M$ for $\binom BA$ and by \scite{1.6E} 
it is enough to prove  $p \in  
{\Cal P}\left( \langle N^4_\eta,a^3_\eta:\eta \in I^2\rangle,M\right)$.  
Now as said above $B \nonfork{}{}_{N^3_{<>}} N^4_{<>}$ and
$B \nonfork{}{}_{A} N^3_{<>}$ so together $B \nonfork{}{}_{A} N^4_{<>}$, 
also we have $A \subseteq B_0 \subseteq N^4_{<>}$,  hence
$B \nonfork{}{}_{B_0} N^4_{<>}$ and \newline
${\frac{B_n}{B_0+B}} = {\frac{B_n}{B_0+a^3_{<0>}}} \perp_a B_0$ 
(by $(\epsilon)^+$ above) but $\nonfork{a^3_{<>}}{N^4_{<>}}_{B_0}$ hence
${\frac{B_n}{N^4_{<>}+a^3_{<0>}}}$ is \newline
$\aleph_\epsilon$-isolated.
Also letting $B'_n = B_n \backslash \{d\}$ we have
${\frac{d}{B'_n}} \perp B_0$ (by clause $(\zeta)$), and clearly
$\nonfork{d}{(N^4_{<>} \cup B'_n)}_{B_n}$ so $\frac{d}{B'_n} \perp
N^4_{<>}$.  Hence we can find $\langle N^5_\eta,a^5_\eta:\eta \in I^5 
\rangle$ an $\aleph_\epsilon$-decomposition of $M$ for $\binom BA$ such 
that $N^5_{<>} = N^4_{<>}$, dcl$(B) = \text{ dcl}(a^5_{<0>})$, \newline
$B_n \backslash \{d\} \subseteq N^5_{<0>}$ and $d = a^5_{<0,0>}$ (on $d$ 
see clause $(\zeta)$ above) so $\nonfork{d}{N^5_{<0>}}_{B_n}$.

By \scite{1.6F} it is enough to show $p \in {\Cal P}(\langle N^5_\eta,
a^5_\eta:\eta \in I^5 \rangle,M)$ which holds trivially as 
tp$(d,B_n \backslash \{d\})$ witness. \hfill$\square_{\scite{1.6}(6)}$
\enddemo
\bigskip

\demo{\stag{1.6G} Observation}  If ${\binom{B_1}{A_1}} \le^* \binom{B_2}{A_2}$
\underbar{then} we can find $\langle B'_\ell:\ell \le n\rangle$ and \newline  
$\langle c_\ell:1 \le \ell < n\rangle$ for some $n \ge 1$, satisfying 
$\binom{B_1}{A_1} \le_b \binom{B'_1}{B'_0},c_\ell \in 
B'_{\ell+1},{\frac{c_\ell}{B'_\ell}}$
regular, ${\frac{B'_{\ell+1}}{c_\ell +B'_\ell}} \perp_a B'_\ell,A_2 = 
B'_{n-1},B_2 = B'_n$.
\enddemo
\bigskip

\remark{Remark}  This means that in a sense $\le_b$ is enough, $\le_a$ 
inessential.
\endremark
\bigskip

\demo{Proof}  By the definition of $\le^*$ there are $k < \omega$ and
$\binom{B^\ell}{A^\ell}$ for $\ell \le k$ such that:
$\binom{B^\ell}{A^\ell} \le_{x(\ell)} \binom{B^{\ell+1}}{A^{\ell+1}}$
for $\ell \le k$ and $x(\ell) \in \{a,b\}$ and 
$\binom{B^0}{A^0} = \binom{B_1}{A_1},\binom{B^k}{A^k} = \binom{B_2}{A_2}$ and
without loss of generality $x(2 \ell) = a,x(2 \ell + 1) = b$.
Let $N_0 \prec {\frak C}$  be $\aleph_\epsilon $-prime over  
$\emptyset$ such that $A^0 \subseteq N_0,B_0 \nonfork{}{}_{A^0} N_0$ and
$f_0 = id_{A_0}$.  We choose by induction on  $\ell \le k,N_{\ell+1},
f_{\ell +1}$ such that: 
\medskip
\roster
\item "{$(a)$}"  Dom$(f_{\ell +1}) = B^\ell$
\item "{$(b)$}"  $N_\ell \prec N_{\ell +1}$ 
\item "{$(c)$}"  if $x(\ell) = b$ then $f_{\ell +1}$ is an extension of 
$f_\ell$ (which necessarily has domain $A_\ell$, check) with domain 
$B^\ell$ such that $f_\ell(B^\ell) \nonfork{}{}_{f_\ell(A^\ell)} N_\ell$ 
and $N_{\ell +1}$ is $\aleph_\epsilon$-prime over $N_\ell \cup f_\ell(B^\ell)$
\item "{$(d)$}"  if $x(\ell) = a$, then $f_{\ell +1}$ maps  
$A^\ell$ into $N_{\ell-1},B^\ell$ into $N_\ell$ and $N_{\ell +1} = N_\ell$.
\endroster
\medskip

\noindent
This is straightforward.  Now on $\langle N_\ell:\ell \le k + 1\rangle$  
we repeat the argument (of choosing $\langle B_\ell:\ell \le n \rangle$)
in the proof of \scite{1.6}(6) above.
\hfill$\square_{\scite{1.6G}}$
\enddemo
\bigskip

\demo{Proof of \scite{1.6}(5)}  By \scite{1.6G}, there are 
$\langle B_\ell:\ell \le n \rangle,\langle c_\ell:1 \le \ell < n \rangle$ as 
there.  We can choose $N^1_{<>}$ such that $B_0 \subseteq N^1_{<>},
N^1_{<1>} \nonfork{}{}_{B_0} B_n,N^1_{<>}$ is 
$\aleph_\epsilon$-primary over  $\emptyset$.  Then we choose  
$\left< N^1_\eta,a^1_\eta:\eta \in \biggl\{ <>,<0>,<0,0>,\dotsc,
\underbrace{<0,\dotsc,0>}_n \biggr\} \right>$,
(where \newline 
$N^1_{\underbrace{<0,\dotsc,0>}_n} \prec M$ and $\ell > 0 \Rightarrow
a^1_{\underbrace{<0,\dotsc,0>}_\ell} = c_\ell$ and we choose $N^1_\eta$ 
by induction on $\ell g(\eta)$ being $\aleph_\varepsilon$-prime over
$N^1_{\eta^-} \cup a^1_\eta$ which is an 
$\aleph_\epsilon$-decomposition of $M$ for $\binom{B_1}{A_1}$.  
Now apply \scite{1.6}(6).  \hfill$\square_{\scite{1.6}(5)}$
\enddemo
\bigskip

\demo{Proof of \scite{1.6}(8)}  Read Definition \scite{1.4}.
\enddemo
\bigskip

\demo{Proof of \scite{1.6}(7)}  Should be easy.  Note that
\medskip
\roster
\item "{$(*)_1$}"   for no $\binom{B'}{A'}$ do we have $\binom BA \le_b
\binom{B'}{A'}$ \newline
Why?  By the definition of Depth zero. 
\smallskip
\noindent
\item "{$(*)_2$}"  if $\binom BA <_a \binom{B'}{A'}$, then also
$\binom{B'}{A'}$ satisfies the assumption.
\endroster
\medskip

\noindent
Hence
\medskip
\roster
\item "{$(**)$}"  for no $\binom{B_1}{A_1},\binom{B_2}{A_2}$ do we have
$$
\binom BA \le_a \binom{B_1}{A_1} <_b \binom{B_2}{A_2}.
$$
\noindent
[Why?  As also $\binom{B_1}{A_1}$ satisfies the assumption].
\endroster
\medskip

\noindent
Now we can prove the statement by induction on $\alpha$ for all pairs
$\binom BA$ satisfying the assumption.  For $\alpha = 0$ the statement is
a tautology.  For $\alpha$ limit ordinal reread clause (c) of Definition
\scite{1.4}(1).  For $\alpha = \beta + 1$, reread clause (b) of Definition 
\scite{1.4}(1):
on $\text{tp}_\beta(\binom BA,M)$ use the induction hypothesis also for
computing $Y^{1,B}_{A,B,M}$ (and reread the definition of $\text{tp}_0$,
in Definition \scite{1.4}(1), clause (a)).  Lastly $Y^{2,\beta}_{A,B,M}$ is
empty by $(*)$ above. \hfill$\square_{\scite{1.6}(5),(7),(8)}$
\enddemo
\bigskip

\demo{Discussion}  In particular, the following Claim \scite{1.7} implies 
that if $\langle N_\eta,a_\eta:\eta \in I \rangle$ is an 
$\aleph_\epsilon$-decomposition of $M$ over $\binom BA$ and  
$M^-$ is $\aleph_\epsilon$-prime over $\cup \{N_\eta:\eta \in I\}$ 
then $\binom BA$ has the same $\text{tp}_\alpha$ in $M$ and $M^-$.
\enddemo
\bigskip

\proclaim{\stag{1.19A}Claim}  1) Assume that $M_1 \prec M_2$ are
$\aleph_\varepsilon$-saturated, $\binom BA \in \Gamma(M_1)$.  Then the
following are equivalent:
\medskip
\roster
\item "{$(a)$}"  if $p \in {\Cal P} \bigl( \binom BA,M_1 \bigr)$ \newline
(see \scite{1.6}(6) for definition; so $p \in S(M_1)$ is regular) then $p$
is not realized in $M_2$  
\item "{$(b)$}"  there is an $\aleph_\varepsilon$-decomposition of $M_2$
above $\binom BA$, which is also an $\aleph_\varepsilon$-decomposition of
$M_2$ above $\binom BA$
\item "{$(c)$}"  every $\aleph_\varepsilon$-decomposition of $M_1$ above
$\binom BA$ is also an $\aleph_\varepsilon$-decomposition of $M_2$ above
$\binom BA$.
\endroster
\medskip

\noindent
2) If $M$ is $\aleph_\varepsilon$-saturated, $\binom{B_1}{A_1} \le^*
\binom{B_2}{A_2}$ are both in $\Gamma(M)$ then ${\Cal P} \bigl( \binom{B_2}
{A_2},M \bigr) \subseteq {\Cal P} \bigl( \binom{B_1}{A_1},M \bigr)$. \newline
3) The conditions in \scite{1.19A}(1) above implies
\medskip
\roster
\item "{$(d)$}"  $p \in {\Cal P} \bigl( \binom BA,M_2 \bigr) \Rightarrow
p \pm M_1$.
\endroster 
\endproclaim
\bigskip

\demo{Proof}  1) \underbar{$(c) \Rightarrow (b)$}.

By \scite{1.6}(1) there is an $\aleph_\varepsilon$-decomposition of $M_1$
above $\binom BA$.  By clause (c) it is also an 
$\aleph_\varepsilon$-decomposition of $M_2$ above $\binom BA$, just as needed
for clause (c).
\medskip

\noindent
\underbar{$(b) \Rightarrow (a)$}

Let $\langle N_\eta,a_\eta:\eta \in I \rangle$ be as said in clause (b).  By
\scite{1.6}(3)(b) we can find $J_1,I \subseteq J_1$ and $N_\eta,a_\eta$ (for
$\eta \in J_1 \backslash I$) such that $\langle N_\eta,a_\eta:\eta \in J_1
\rangle$ is an $\aleph_\varepsilon$-decomposition of $M_1$ and $\nu \in J_1
\backslash I \Rightarrow \nu(0) > 0$.  Then we can find $J_2,J_1 \subseteq
J_2$ and $N_\eta,a_\eta$ (for $\eta \in J_2 \backslash J_1$) such that
$\langle N_\eta:\eta \in J_2 \rangle$ is an $\aleph_\varepsilon$-decomposition
of $M_2$ (by \scite{1.6}(2)).  By \scite{1.6}(3)(b), $\nu \in J_2 \backslash
I \Rightarrow\nu(0) > 0$.  So $\eta \in I \backslash \{ \langle \rangle \}
\Rightarrow \text{ Suc}_{J_2}(\eta) = \text{ Suc}_I(\eta)$, hence
\medskip
\roster
\item "{$(*)$}"  $\eta \in I \backslash \{ \langle \rangle\},q \in S(N_\eta)$
regular orthogonal to $N_{\eta^-}$ then the stationarization of $q$ in
$S(M_1)$ is not realized in $M_2$.
\endroster
\medskip

\noindent
Now if $p \in {\Cal P} \bigl( \binom BA,M_1 \bigr)$ then $p \in S(M_1)$ is
regular and (see \scite{1.6}(1), \scite{1.5}(5)) for some $\eta \in I
\backslash \{ \langle \rangle \},p \perp N_{\eta^-},p \pm N_\eta$, so there is
a regular $q \in S(N_\eta)$ not orthogonal to $p$.  Now no $c \in M_2$
realizes the stationarization of $q$ over $M_1$ (by $(*)$ above), hence this
applies to $p$, too.
\medskip

\noindent
\underbar{$(a) \Rightarrow (c)$}

Let $\langle N_\eta,a_\eta:\eta \in I \rangle$ be an 
$\aleph_\varepsilon$-decomposition of $M_1$ above $\binom BA$.  We can find
$\langle N^*_\eta,a^*_\eta:\eta \in J \rangle$ an 
$\aleph_\varepsilon$-decomposition of $M_1$ such that $I \subseteq J$ and
$\nu \in J \backslash I \Rightarrow \nu(0) > 0$ (by \scite{1.6}(3)(b)).  We
should check that it is also an $\aleph_\varepsilon$-decomposition of $M_2$
above $\binom BA$, i.e. Definition \scite{1.5}(1),(2).  Now in \scite{1.5}(1),
clauses (a)-(h) are immediate, so let us check clause (i) (in \scite{1.5}(2)).
Let $\eta \in I \backslash \{ \langle \rangle\}$, now is 
$\{ a_{\eta \char 94 \langle \alpha \rangle}:\eta \char 94 \langle \alpha
\rangle \in I\}$ really maximal (among independent over $N_\eta$ sets of
elements of $M_2$ realizing a type from ${\Cal P}_\eta = \{ p \in S(N_\eta):p
\text{ regular orthogonal to } N_{\eta^-}\}$.  This should be clear from
Claim (a) (and basic properties of dependencies for regular types). \newline
2) By \scite{1.6}(5). \newline
3) Left to the reader. \hfill$\square_{\scite{1.19A}}$
\enddemo
\bigskip

\noindent
\underbar{\stag{1.19B}Conclusion}:  Assume $M_1 \prec M_2$ are
$\aleph_\varepsilon$-saturated and $\binom{B_1}{A_1} \le^* \binom{B_2}{A_2}$
both in $\Gamma(M_1)$.  If clause $(a)$ (equivalently $(b)$ or $(c)$) of
\scite{1.19A} holds for $\binom{B_1}{A_1},M_1,M_2$ then they hold for
$\binom{B_2}{A_2},M_1,M_2$.
\bigskip

\demo{Proof}  By \scite{1.19A}(1), clause (a) for $\binom{B_1}{A_1},M_1,M_2$
implies clause (a) for $\binom{B_2}{A_2},M_1,M_2$. \newline
${}$ \hfill$\square_{\scite{1.19B}}$
\enddemo   
\bigskip

\proclaim{\stag{1.7} Claim}  If $\binom{B_1}{A_1} \in \Gamma(M)$ and $\langle 
N_\eta,a_\eta:\eta \in I \rangle$ is an $\aleph_\epsilon$-decomposition 
of $M$ above $\binom{B_1}{A_1}$ and $M^- \subseteq M$ is 
$\aleph_\epsilon$-saturated and $\dsize \bigcup_{\eta \in I} N_\eta 
\subseteq M^-$ \underbar{then} 
\smallskip
$$
\text{tp}_\alpha \bigl[ \binom{B_1}{A_1},M \bigr] = 
\text{ tp}_\alpha \bigl[ \binom{B_1}{A_1},M^- \bigr] 
$$
\endproclaim
\bigskip

\demo{Proof}  We prove by induction on $\alpha$ (for all $B,A,\langle 
N_\eta,a_\eta:\eta \in I \rangle,I,M$ and $M^-$ as above).  We can find 
an $\aleph_\epsilon $-decomposition  
$\langle N_\eta,a_\eta:\eta \in J\rangle$  of $M$ with $I \subseteq J$
(by \scite{1.5B}(4) + \scite{1.5B}(2)) such that $\eta \in J \backslash I 
\Leftrightarrow \eta \ne \langle \rangle \and \neg \langle 0 \rangle 
\trianglelefteq \eta$ and $M$ is $\aleph_\epsilon$-prime over $\dsize 
\bigcup_{\eta \in J} N_\eta$.
\enddemo 
\bigskip

\noindent
\underbar{Case 0}:  $\alpha = 0$. 

Trivial.
\bigskip

\noindent
\underbar{Case 1}:  $\alpha$  is a limit ordinal. 

Trivial by induction hypothesis (and the definition of tp$_\alpha$).
\bigskip

\noindent
\underbar{Case 2}:  $\alpha = \beta + 1$.

We can find $M^* \prec M^-$ which is $\aleph_\epsilon$-prime over
$\dsize \bigcup_{\eta \in I} N_\eta$, so as equality is transitive it is
enough to prove

$$
\text{tp}_\alpha \left( \binom{B_1}{A_1},M^* \right) =
\text{ tp} \left( \binom{B_1}{A_1},M^- \right)
$$
\medskip

\noindent
and

$$
\text{tp} \left( \binom{B_1}{A_1},M^* \right) = \text{ tp}
\left( \binom{B_1}{A_1},M \right).
$$
\medskip

\noindent
Reflecting, this means that it is enough to prove the statement when $M^-$
is \newline
$\aleph_\epsilon$-prime over $\dsize \bigcup_{\eta \in I} N_\eta$.

Looking at the definition of  $\text{tp}_{\beta +1}$ and remembering 
the induction hypothesis our problems are as follows: \newline
\underbar{First component of tp$_\alpha$}: \newline
\smallskip
\noindent
given  $\binom{B_1}{A_1} \le_a \binom{B_2}{A_2},B_2 \subseteq M$, and
\underbar{it suffices to find}  $\binom{B_3}{A_3}$ such that: 
\medskip
\roster
\item "{$(*)$}"  there is  $f \in \text{ AUT}({\frak C})$ such that:
$f \restriction B_1 = \text{ id}_{B_1},f(A_2) = A_3$, \newline  
$f(B_2) = B_3$ and  $B_3 \subseteq  M^-$ and $\text{tp}_\beta
[\binom{B_2}{A_2},M] = \text{ tp}_\beta[\binom{B_3}{A_3},M^-]$ \newline
(pedantically we should replace $B_\ell,A_\ell$ by indexed sets).
\endroster
\medskip

\noindent
We can find $J',M'$ such that: 
\medskip
\roster
%
\item "{$(i)$}"  $I \subseteq J' \subseteq J,|J'\backslash I| < \aleph_0,J'$
closed under initial segments,
\item "{$(ii)$}"  $M' \prec  M$  is $\aleph_\epsilon$-prime over $M^- \cup  
\cup \{N_\eta:\eta \in  J'\backslash I\}$
\item "{$(iii)$}"  $B_2 \subseteq M'$.
\endroster
\medskip

\noindent
The induction hypothesis for $\beta$ applies, and gives

$$
\text{tp}_\beta \bigl[ \binom{B_2}{A_2},M \bigr] = \text{ tp}_\beta
\bigl[ \binom{B_2}{A_2},M' \bigr].
$$
\medskip

\noindent
By \scite{1.6}(4) there is $g$, an isomorphism from $M'$ onto $M^-$ 
such that $g \restriction B_1 = \text{ id}$.  So clearly $g(B_2) \subseteq 
M^-$ hence

$$
\text{tp}_\beta \bigl[ \binom{B_2}{A_2},M' \bigr] = 
\text{ tp}_\beta \bigl[\binom{g(B_2)}{g(A_2)},M^- \bigr].
$$
\medskip

\noindent
So  $\binom{B_3}{A_3} =: g\binom{A_2}{B_2}$
is as required.
\medskip

\noindent
\underbar{Second component for tp$_\alpha$}: \newline
\smallskip
\noindent
So we are given $\Upsilon$, a tp$_\beta$ type,  (and we assign the lower 
part as $B$)  and we have to prove that the dimension in  $M$  and in  $M^-$ 
are the same, i.e. \newline
Dim$(\bold I,M) = \text{ dim}(\bold I^-,M)$, where:
$\bold I = \{c \in M:\Upsilon = \text{ tp}_\beta\bigl(\binom{c}{B_1},M
\bigr)$ and \newline
$\bold I^- = \bigl\{ c \in  M^-:\Upsilon = \text{ tp}_\beta \bigl( 
\binom{c}{B_1},M^- \bigr) \bigr\}$.  \newline
Let  
$p$ be such that: tp$_\beta\left(\binom{c}{B_1},M \right) = \Upsilon
\Rightarrow  p = {\frac{c}{B_1}}$.  Necessarily $p \perp A_1$ and $p$ is 
regular (and stationary).

Clearly $\bold I^- \subseteq \bold I$, now
\medskip
\roster
\item "{$(*)$}"  for every $c \in \bold I$ for some $k < \omega,c'_\ell \in
M^-$ realizing $p$ for $\ell < k$ we have $c$ depends on
$\{c'_0,c'_1,\dotsc,c'_{k-1}\}$ over $B_1$. \newline
[Why?  Clearly $p \perp N_{<>}$ (as $\nonfork{B_1}{N_{<>}}_{A_1}$ and
$p \perp A_1$) hence \newline
tp$_*( \dsize \bigcup_{\eta \in J
\backslash I} N_\eta,N_{<>}) \perp p$ hence  \newline
tp$_*( \dsize \bigcup_{\eta \in J \backslash I} N_\eta,M^-) \perp p$, but
$M$ is $\aleph_\epsilon$-prime over $M^- \cup \dsize \bigcup_{\eta \in J
\backslash I} N_\eta$ hence \newline
by \cite[V,3.2,p.250]{Sh:c} for no $c \in M \backslash M^-$ is
tp$(c,M^-)$ a stationarization of $p$ hence by \cite[V,1.16(3)]{Sh:c} clearly 
$(*)$ follows].
\endroster
\medskip

\noindent
If the type $p$ has depth zero, then by (\scite{1.6}(7)) 

$$
\gather
\bold I = \{c \in M:\text{tp}(c,B) = p\} \text{ and} \\
\bold I^- = \{c \in M^-:\text{tp}(c,B) = p\}. \endgather
$$
\medskip

\noindent
Now we have to prove dim$(\bold I,A) = \text{ dim}(\bold I^-,A)$, as $A$ is
$\varepsilon$-finite and $M,M^-$ are $\aleph_\epsilon$-saturated and
$\bold I^- \subseteq \bold I$ clearly $\aleph_0 \le \text{ dim}(\bold I^-,A)
\le \text{ dim}(\bold I,A)$.  Now the equality follows by $(*)$ above.

So we can assume ``$p$ has depth zero", hence (by \cite[X,7.2]{Sh:c}) that 
the type $p$ is trivial; hence, see \cite[X7.3]{Sh:c} in $(*)$ without 
loss of generality $k=1$ and
dependency is an equivalence relation, so for ``same dimension" 
it suffices to prove that every equivalence class (in $M$ i.e. in $\bold I$)
is representable in $M^-$ i.e. in $\bold I^-$.  By the remark on $(*)$ in
the previous sentence 
$(\forall d_1 \in \bold I)(\exists d_2 \in \bold I^-) \bigl[ \neg d_1
\nonfork{}{}_{B_1} d_2 \bigr]$.  So it is enough to prove that:
\medskip
\roster
\item "{$\bigotimes$}"   if $d_1,d_2 \in M$ realize same type over $B_1$,
which is (stationary and) regular, and are dependent over $B_1$ and
$d_1 \in  M^-$ \underbar{then} there is  $d'_2 \in  M^-$
such that ${\frac{d'_2}{B_1+d_1}} = {\frac{d_2}{B_1+d_1}}$ and   
tp$_\beta \bigl[\binom{B_1+d_2}{B_1},M \bigr] = \text{ tp}_\beta \bigl[ 
\binom{B_1+d'_2}{B_1},M^- \bigr]$.
\endroster
\medskip

\noindent
Let $M_0 = N_{\langle \rangle}$.  There are $J',M_1,M^+_1$ such that
\medskip
\roster
%
\item "{$(*)_1(i)$}"   $J' \subseteq J$ is finite closed under initial 
segments
\item "{${}(ii)$}"     $\langle \rangle \in J',\langle 0 \rangle \notin J'$
\item "{${}(iii)$}"    $M_1 \prec M$ is $\aleph_\varepsilon$-prime over
$\cup\{N_\eta:\eta \in J'\}$
\item "{${}(iv)$}"     $M^+_1 \prec M$ is $\aleph_\varepsilon$-prime over
$M_1 \cup M^-$ (and $\nonfork{M_1}{M^-}_{M_0}$)
\item "{${}(v)$}"      $d_2 \in M^+_1$.
\endroster
\medskip

\noindent
Now the triple $\binom{B_1+d_2}{B_1},M_1,M$ satisfies the demand on
$\binom{B_1}{A_1},M^-,M$ (because \newline
$\binom{B_1}{A_1} \le^* \binom{B_1+d_2}{B_1}$,by \scite{1.19B}.  
Hence by the induction hypothesis we know that

$$
\text{tp} \bigl[ \binom{B_1+d_2}{B_1},M \bigr] = \text{ tp}_\beta \bigl[
\binom{B_1+d_2}{B_1},M^+_1 \bigr].
$$
\medskip

\noindent
By \scite{1.10}(4) there is an isomorphism $f$ from $M^+_1$ onto $M^-$ which
is the identity on $B_1 + d_1$; let $d'_2 = f(d_2)$ so:

$$
\text{tp}_\beta \bigl[ \binom{B_1+d_2}{B_1},M^+_1 \bigr] = \text{ tp}_\beta 
\bigl[ \binom{B_1+d'_2}{B_1},M^- \bigr].
$$
\medskip

\noindent
Together

$$
\text{tp}_\beta \bigl[ \binom{B_1+d_2}{B_1},M \bigr] = \text{ tp}_\beta \bigl[
\binom{B_1+d'_2}{B_1},M^- \bigr].
$$
\medskip

\noindent
As $\{d_1,d_2\}$ is not independent over $B_1$, also $\{d_1,f(d_2)\}$ is not
independent over $B_1$, hence, as $p$ is regular
\medskip
\roster
\item "{$(*)$}"  $\{d_2,f(d_2)\}$ is not independent over $B_2$.
\endroster
\medskip

\noindent
Together we have proved $\bigoplus$, hence finished proving the equality
of the second component.
\bigskip

\noindent
\underbar{Third component}:  Trivial.

So we have finished the induction step, hence the proof.
\hfill$\square_{\scite{1.7}}$
\bigskip

\proclaim{\stag{1.8} Claim}  1) Suppose $M$ is $\aleph_\epsilon$-saturated, 
$B \subseteq M,\binom BA \in \Gamma,\dsize \bigwedge_{\ell = 1}^{2}
[A \subseteq A_\ell \subseteq M]$, \newline
$A = ac\ell(A),A_\ell$ are 
$\epsilon$-finite, ${\frac{A_1}{A}} = {\frac{A_2}{A}},\nonfork{B}{A_1}_{A}$
and $\nonfork{B}{A_2}_{A}$. 

\underbar{Then} tp$_\alpha\bigl[\binom{A_1 \cup B}{A_1},M\bigr]  = 
\text{ tp}_\alpha\bigl[\binom{A_2 \cup B}{A_2},M\bigr]$ for any ordinal
$\alpha$. \newline
2) Suppose  $M$  is $\aleph_\epsilon $-saturated,  $B \subseteq  M$,  
$\binom BA \in \Gamma,\dsize \bigwedge_{\ell =1}^{2}[A \subseteq A_\ell
\subseteq M],A = ac\ell(A)$, \newline
$B = ac\ell(B),A_\ell = ac\ell(A_\ell),A_\ell$ 
is $\epsilon$-finite, ${\frac{A_1}{A}} = {\frac{A_2}{A}},B \nonfork{}{}_{A}
A_1,B \nonfork{}{}_{A} A_2,f:A_1 \overset\text{onto} {}\to \longrightarrow
A_2$ an elementary mapping, $f \restriction A = \text{ id},g \supseteq f 
\cup \text{ id}_B,g$  elementary mapping from \newline
$B_1 = ac\ell(B \cup A_1)$ onto $B_2 = ac\ell(B \cup A_2)$. 
\smallskip
\noindent
\underbar{Then} $g \left(\text{tp}_\alpha\bigl[\binom{B_1}{A_1},M\bigr]\right)
 = \text{ tp}_\alpha\bigl[\binom{B_2}{A_2},M\bigr]$ for any ordinal $\alpha$.
\endproclaim
\bigskip

\demo{Proof}  1)  Follows from part (2). \newline
2)  We can find  $A_3 \subseteq  M$  such that: 
\medskip
\roster
\item "{$(i)$}"  ${\frac{A_3}{A}} = {\frac{A_1}{A}}$
\smallskip
\noindent
\item "{$(ii)$}"  $A_3 \nonfork{}{}_{A} (B \cup  A_1 \cup  A_2)$.
\endroster
\medskip

\noindent
Hence without loss of generality  $A_1 \nonfork{}{}_{B} A_2$ and even 
$\nonfork{}{}_{A} \{B,A_1,A_2\}$.  Now we can find  $N_{<>}$, an 
$\aleph_\epsilon$-prime model over $\emptyset,N_{<>} \prec M,A \subseteq 
N_{<>}$ and $(B \cup A_1 \cup A_2) \nonfork{}{}_{A} N_{<>}$ (e.g. choose 
$\{A^\alpha_1 \cup A^\alpha_i \cup B^\alpha:\alpha \le \omega \} \subseteq M$
indiscernible over $A,A^\omega_1 = A_1,A^\omega_2 = A_2,B^\omega = B$ and 
let $N_{<>} \prec  M$  be $\aleph_\epsilon $-primary over $\dsize \bigcup_{n
< \omega} (A^n_1 \cup  A^n_2 \cup  B^n \cup  A))$. 

Now find  $\langle N_\eta,a_\eta:\eta \in J\rangle$  an 
$\aleph_\epsilon$-decomposition of $M$ with \newline
dcl$(a_{<0>}) = \text{ dcl}(B),
\text{dcl}(a_{<1>}) = \text{ dcl}(A_1),\text{dcl}(a_{<2>}) = 
\text{ dcl}(A_2)$. 

Let $I = \{\eta \in J:\eta = <> \text{ or } <0> \trianglelefteq \eta \}$ 
and  $J' = I \cup  \{<1>,<2>\}$.  Let  
$N^2_{<>} \prec  M^*$ be $\aleph_\epsilon$-prime over  
$N_{<1>} \cup  N_{<2>}$.  By \scite{1.6D} there is $\langle N^2_\eta,a_\eta:
\eta \in I \rangle$ an $\aleph_\epsilon$-decomposition of $M$ over 
$\binom BA$ such that  
$\langle N_\eta,a_\eta:\eta \in I\rangle \le_{\text{direct}} 
\langle N^2_\eta,a_\eta:\eta \in I\rangle$.  Let $M' \prec M$ be 
$\aleph_\epsilon$-prime over $\dsize \bigcup_{\eta \in I} N^2_\eta$ and
$M^- \prec M'$ be $\aleph_\epsilon$-prime over $\dsize 
\bigcup_{\eta \in I} N_\eta$.  So $M^- \prec M' \prec M$ and $M'$ is 
$\aleph_\epsilon$-prime over $M^- \cup N_{<1>} \cup N_{<2>}$.

Now by \scite{1.7} we have tp$_\alpha\bigl[\binom{B}{A_\ell},M\bigr]  = 
\text{ tp}_\alpha\bigl[\binom{B}{A_\ell},M'\bigr]$ for $\ell = 1,2$ hence
it suffices to find an automorphism of $M'$ extending $g$.  Let $B^+ =
ac \ell(N_{<>} \cup B),A^*_\ell = ac \ell(B \cup A_\ell)$, let
$\bold{\bar a}_\ell$ list $A^*_\ell$ be such that $\bold{\bar a}_2 =
g(\bold{\bar a}_1)$.  Clearly tp$(\bold{\bar a}_\ell,B^+)$ does not fork
over $B$ and $ac \ell(B) = B$ and so stp$(\bold{\bar a}_1,B^+) =
\text{ stp}(\bold{\bar a}_2,B^+)$.  Also tp$_*(A_2,B^+ \cup A_1)$ does not
fork over $A$ hence tp$(\bold{\bar a}_2,B^+ \cup \bold{\bar a}_1)$ does not
fork over $A \subseteq B^+$ hence $\{\bold{\bar a}_1,\bold{\bar a}_2\}$ is
independent over $B^+$ hence there is an elementary mapping $g^+$ from
$ac \ell(B^+ \cup \bold{\bar a}_1)$ onto $ac \ell(B^+ \cup \bold{\bar a}_2),
g^+ \supseteq \text{ id}_{B^+} \cup g$ and even $g' = g^+ \cup (g^+)^{-1}$ is
an elementary embedding.

Let $\bold{\bar a}'_1$ lists $ac \ell(N_{<>} \cup A_1)$ so clearly
$\bold{\bar a}'_2 =: g^+(\bold{\bar a}'_1)$ list $ac \ell(N_{<>} \cup A_2)$.
Clearly
$g' \restriction (\bold{\bar a}'_1 \cup \bold a'_2)$ is an elementary mapping
from $\bold{\bar a}'_1 \cup \bold{\bar a}'_2$ onto itself.  Now $N^2_{<>}$
is $\aleph_\epsilon$-primary over $N_{<>} \cup A_1 \cup A_2$ and
$N_{<>} \cup A_1 \cup A_2 \subseteq \bold{\bar a}'_1 \cup \bold{\bar a}'_2
\subseteq ac \ell(N_{<>} \cup A_1 \cup A_2)$ so easily $N^2_{<>}$ is
$\aleph_\epsilon$-primary over $\bold{\bar a}'_1 \cup \bold{\bar a}'_2$ hence
we can extend $g' \restriction (\bold{\bar a}'_1 \cup \bold{\bar a}'_2)$ to
an automorphism $h_{<>}$ of $N^2_{<>}$ so clearly $h_{<>} \restriction
N_{<>} = \text{ id}_{N_{<>}}$.  Let $\bold{\bar a}^+_1$ list $ac \ell(B^+ \cup
A_1)$ and $\bold{\bar a}^+_2 = g^+(\bold{\bar a}^+_1)$.  So
tp$(\bold{\bar a}^+_\ell,N^2_{<>})$ does not fork over $\bold{\bar a}'_1
(\subseteq N^2_{<>})$ and $ac \ell(\bold{\bar a}'_1) = \text{ Rang}
(\bold{\bar a}'_1)$ ($= ac \ell(N_{<>} \cup A_1)$) and $h_{<>} \restriction
\bold{\bar a}'_1 = g^+ \restriction \bold{\bar a}'_1$ hence $h_{<>} \cup g^+$
is an elementary embedding.  Remember $g^+$ is the identity on $B^+ =
ac \ell(N_{<>} \cup B)$, and tp$_*(N_{<0>},N^2_{<>})$ does not fork over
$N_{<>}$ hence tp$_*(N_{<0>},B^+ \cup N^2_{<>})$ does not fork over $B^+$, so
as $ac \ell(B^+) = B^+$ necessarily $(h_{<>} \cup g^+) \cup \text{ id}
_{N_{<0>}}$ is an elementary embedding.  But this mapping has domain and range
including $N_{<0>} \cup N^2_{<>}$ and included in $N^2_{<0>}$, but the latter
is $\aleph_\epsilon$-primary and $\aleph_\epsilon$-minimal over the former.
Hence $(h_{<>} \cup g^+) \cup \text{ id}_{N_{<>}}$ can be extended to an
automorphism of $N^2_{<0>}$ which we call $h_{<0>}$.

Now we define by induction on $n \in [2,\omega)$ for every $\eta \in I$ 
of length $n$, an automorphism $h_\eta$ of $N^2_\eta$ extending $h_{\eta^-} 
\cup \text{ id}_{N_\eta}$, which exists as $N^2_\eta$ is 
$\aleph_\epsilon$-primary over $N^2_{\eta^-} \cup N_\eta$ (and
$\nonfork{N^2_{\eta^-}}{N_\eta}_{N_{\eta^-}}$).  Now 
$\dsize \bigcup_{\eta \in I} h_\eta$ is an elementary mapping (as 
$\langle N^2_\eta:\eta \in I \rangle$ is a non-forking tree;
i.e. \scite{1.5B}(10)), with domain and range 
$\dsize \bigcup_{\eta \in I} N^2_\eta$ hence can be extended to an 
automorphism $h^*$ of $M'$, (we can demand 
$h^* \restriction M^- = \text{ id}_{M^-}$ but not necessarily).  So as
$h^*$ extends $g$, the conclusion follows. 
\hfill$\square_{\scite{1.8}}$
\enddemo
\bigskip

\proclaim{\stag{1.9} Claim}  1) For every $\Upsilon = \text{ tp}_\delta
\bigl[\binom BA,M\bigr]$, and $\bold{\bar a},\bold{\bar b}$ listing $A,B$
respectively there is $\psi = \psi(\bar x_A,\bar x_B)  \in  
{\Cal L}_{\infty,\aleph_\epsilon}(q.d.)$ of depth $\delta$  such that:

$$
\text{tp}_\delta \bigl[\binom BA,M \bigr]  = \Upsilon  
\Leftrightarrow  M \models \psi [\bold{\bar a},\bold{\bar b}].
$$
\medskip

\noindent
2) Assume $\bigotimes_{M_1,M_2}$ of \scite{1.1B} holds as exemplified by 
${\Cal F}$ and $\binom BA \in \Gamma(M_1)$ and \newline
$g \in {\Cal F},\text{Dom}(g) = B$; and $\alpha$ an ordinal then

$$
\text{tp}_\alpha \left(\binom BA,M \right) = \text{ tp}_\alpha \left(
\binom{g(B)}{g(A)},M_2 \right).
$$
\medskip

\noindent
3) Similarly for tp$_\alpha([A],M),\text{tp}_\alpha[M]$.
\endproclaim
\bigskip

\demo{Proof}  Straightforward (remember 
we assume every first order formula is equivalent to a predicate). 
\hfill$\square_{\scite{1.9}}$
\enddemo
\bigskip

\demo{\stag{1.10} Proof of Theorem 1.2}  [The proof does not require that the
$M^\ell$ are $\aleph_\epsilon$-saturated, but only that \scite{1.8}, 
\scite{1.9} hold except in constructing $g_{\alpha(*)}$ (see 
$\otimes_{14},\otimes_{15}$ in \scite{1.11}(E), we could instead use NOTOP]. 

So suppose 
\medskip
\roster
\item "{$(*)_0$}"  $M^1 \equiv_{{\Cal L}_{\infty,\aleph_\epsilon}(d.q.)}M^2$ 
or (at least) $\bigotimes_{M^1,M^2}$ from \scite{1.1B} holds.
\endroster
\medskip 

\noindent
We shall prove  $M^1 \cong M^2$.  By \scite{1.9} (i.e. by \scite{1.9}(1) if
the first possibility in $(*)_0$ holds and by \scite{1.9}(2) if the second
possibility in $(*)_0$ holds) 
\medskip
\roster
\item "{$(*)_1$}"   tp$_\infty[M^1] = \text{ tp}_\infty[M^2]$.
\endroster
\enddemo
\bigskip

\noindent
So it suffices to prove:
\proclaim{\stag{1.11} Claim}  If $M^1,M^2$ are 
$\aleph_\epsilon$-saturated models (of $T,T$ as in \scite{1.2}), then:
\medskip

$(*)_1 \Rightarrow  M^1 \cong M^2$.
\endproclaim
\bigskip

\demo{Proof}  Let $\langle W_k,W'_k:k < \omega \rangle$  
be a partition of $\omega$ to infinite sets (so pairwise disjoint).
\medskip

\noindent
\underbar{\stag{1.11a} Explanation}:  (If seems opaque, 
the reader may return to it after reading part of the proof). 

We shall now define an approximation to a decomposition.  
We are approximating a non-forking trees $\langle N^\ell_\eta,a^\ell_\eta:
\eta \in I^* \rangle$  of countable elementary submodels of $M^\ell$ and
$\langle f^*_\eta:\eta \in I^\ast \rangle$ such that $f^*_\eta$ an 
isomorphism from
$N^1_\eta$ onto $N^2_\eta$ increasing with $\eta$ such that $M^\ell$ is 
$\aleph_\epsilon$-prime over  $\dsize \bigcup_{\eta \in I^*} N^\ell_\eta$. 
\medskip
\roster
\item "{$(\alpha)$}"  In the approximation we have: $I$ approximating 
$I^*$, \newline
[it will not be  $I^* \cap {}^{n \ge}$ordinal and we may ``discover" more 
immediate successor to each  $\eta \in I$; as the approximation to  
$N_\eta$ improves we have more regular types, but some member of $I$ will
be later will drop] 
\item "{$(\beta)$}"  $A^\ell_\eta$ approximates $N^\ell_\eta$ and is 
$\epsilon$-finite
\item "{$(\gamma)$}"  $a^\ell_\eta$ is the $a^\ell_\eta$ (if $\eta$ survives,
i.e. will not be dropped)
\item "{$(\delta)$}"  $B^\ell_\eta,b^\ell_{\eta,m}$ expresses commitments 
on constructing $A^\ell_\eta$: we ``promise" \newline
$B^\ell_\eta \subseteq N^\ell_\eta$ and $B^\ell_\eta$ is countable; 
$b^\ell_{\eta,m}$ for $m < \omega$ list
$B^\ell_\eta$ (so in the choice $B^\ell_\eta \subseteq M^\ell$ there is 
some arbitrariness)
\item "{$(\varepsilon)$}"  $f_\eta$ approximate $f^*_\eta$ 
\item "{$(\zeta)$}"  $p^\ell_{\eta,m}$ also expresses commitments on the 
construction.
\endroster
\medskip  

Since there are infinitely many commitments that we must meet in a 
construction of length $\omega$ and we want many chances to meet each of 
them, the sets  $W_k,W'_k$  are introduced as a further bookkeeping device.
At stage  $n$ in the construction we will deal e.g. with the $b^\ell_{\eta,m}$
for  $\eta$  that are appropriate and for  $m \in W_k$ for some $k < n$  and 
analogously for  $p^\ell_{\eta,m}$ and the $W'_k$. 

Note that while the  $A^\ell_\eta$ satisfy the independence properties of a 
decomposition, the  $B^\ell_\eta$ do not and may well intersect 
non-trivially.  Nevertheless, a conflict arises if an  
$a^\ell_{\eta \char 94 <i>}$ falls into  $B^\ell_\eta$ since the  
$a^\ell_{\eta \char 94 <i>}$ are supposed to represent independent elements
realizing regular types over the model approximated by  $A^\ell_\eta$ but 
now $a^\ell_{\eta \char 94 <i>}$ is in that model.  This problem is 
addressed by pruning  $\eta \char 94 <i>$  from the tree $I$.  
\enddemo
\bigskip

\definition{\stag{1.11A} Definition}  An approximation $Y$ to an isomorphism 
consist of: 
\medskip
\roster
\item "{$(a)$}" \underbar{natural numbers $n,k^*$ and index set}:
$I \subseteq {}^{n \ge}$ ordinals   \newline
(and  $n$  minimal)
\item "{$(b)$}"  $\langle A^\ell_\eta,B^\ell_\eta,a^\ell_\eta,b^\ell_{\eta,m}:
\eta \in I$ and $m \in \dsize \bigcup_{k < k^*} W_k \rangle$ for $\ell = 1,2$
(this is an approximated decomposition)
\item "{$(c)$}"  $\langle f_\eta:\eta \in I\rangle$ 
\item "{$(d)$}"  $\langle p^\ell_{\eta,m}:\eta \in I \text{ and }
m \in \dsize \bigcup_{k < k^*} W'_k \rangle$
\endroster
\medskip

\noindent
such that:  
\medskip
\roster
\item  $I$ closed under initial segments  
\item  $<> \in  I$  
\item  $A^\ell_\eta \subseteq B^\ell_\eta \subseteq M^\ell,A^\ell_\eta$ is 
$\epsilon$-finite, $ac\ell(A^\ell_\eta) = A^\ell_\eta,B^\ell_\eta$ is
countable, \newline
$B^\ell_\eta = \{ b^\ell_{\eta,m}:m \in \dsize 
\bigcup_{k < k^*} W_k \}$
\item  $A^\ell_\nu \subseteq A^\ell_\eta$ if $\nu \triangleleft \eta \in I$  
\item   if $\eta \in I \backslash \{<>\}$, then ${\frac{a^\ell_\eta}
{A^\ell_{(\eta^-)}}}$ is a (stationary) regular type and $a^\ell_\eta \in 
A^\ell_\eta$; if in addition $\ell g(\eta) > 1$ then 
${\frac{a^\ell_\eta}{A^\ell_{(\eta^-)}}} \perp 
A^\ell_{(\eta^{--})}$
(note $a^\ell_{<>}$ may be not defined or $\in A^\ell_{<>}$)
\item  ${\frac{A^\ell_\eta}{A^\ell_{\eta^-} + a_\eta}} \perp_a 
A^\ell_{\eta^-}$ if $\eta \in I,\ell g(\eta) > 0$  
\item   for  $\eta \in I$, not $\triangleleft$-maximal in $I$ the set  
$\{a^\ell_\nu:\nu \in I \text{ and } \nu^- = \eta\}$ is a maximal family 
of elements realizing 
over  $A^\ell_\eta$ regular types  $\perp A^\ell_{(\eta^-)}$ (when  
$\eta^-$ is defined), independent over $\left( A^\ell_\eta,B^\ell_\eta
\right)$, (and we can add: if \newline
$\nu^-_1 = \nu^-_2 = \eta$ and
$\frac{a^\ell_{\nu_1}}{A_\eta} \pm \frac{a^\ell_{\nu_2}}{A_\eta}$ then
$a^\ell_{\nu_1}/A_\eta = a^\ell_{\nu_2}/A_\eta$)  
\item   $f_\eta$ an elementary map from  $A^1_\eta$ onto $A^2_\eta$
\item   $f_{(\eta^-)} \subseteq f_\eta$ when $\eta \in I,\ell g(\eta) > 0$
\item   $f_\eta(a^1_\eta) = a^2_\eta$
\item $(\alpha) \quad f_\eta \left( \text{tp}_\infty \left[ 
\binom{A^1_\eta}{A^1_{(\eta^-)}},M^1 \right] \right) = 
\text{ tp}_\infty \left[ \binom{A^2_\eta}{A^2_{(\eta^-)}},M^2 \right]$ when
$\eta  \in  I\backslash \{<>\}$ \newline
\smallskip
\noindent
$\quad(\beta) \qquad 
f_{<>} \left( \text{tp}_\infty \left[ A^1_{<>},
M^1 \right] \right)  = \text{ tp}_\infty \left[ A^2_{<>},M^2 \right]$
\item   $B^\ell_\eta \prec M^\ell$   
\item   $\langle p^\ell_{\eta,m}:m \in \dsize \bigcup_{k < k^*} W'_k
\rangle$  is a sequence of types over  $A^\ell_\eta$.
\endroster
\enddefinition
\bigskip

\demo{\stag{1.11B} Notation}  We write 
$n = n_Y = n[Y],I = I_Y = I[Y],A^\ell_\eta = A^\ell_\eta[Y]$, \newline
$B^\ell_\eta = B^\ell_\eta[Y],f_\eta = f^Y_\eta 
= f_\eta[Y],a^\ell_\eta = a^\ell_\eta[Y],b^\ell_\eta = b^\ell_\eta[Y],
k^\ast = k^\ast_Y = k^\ast[Y]$ and \newline
$p^\ell_{\eta,m}
= p^\ell_{\eta,k}[Y]$.
\enddemo
\bigskip

\remark{Remark}  We may want to demand: each 
${\frac{a^\ell_{\eta \char94 <i>}}{A_\eta}}$ is strongly regular; also: 
if two such types are not orthogonal 
then they are equal (or at least have same witness  $\varphi$  for  
$\left( \varphi,{\frac{a_{\eta \char94 <i>}}{A_\eta}} \right)$ regular).  Why 
can we?  As the models are $\aleph_\epsilon$-saturated (so take $p' \pm  p$,
$rk(p')$ minimal).
\endremark
\bigskip

\demo{\stag{1.11C} Observation}  $(*)_1$ implies that 
there is an approximation.
\enddemo
\bigskip

\demo{Proof}  Let $I = \{<>\},A^\ell_{<>} = ac\ell(\emptyset),k^* = 1$ and
then choose countable $B^\ell_{<>}$ to satisfy condition (12) and then 
choose $f_\eta,p^\ell_k,b^\ell_{\eta,m}(k \in W_0)(m \in W_0)$ as required.
\enddemo
\bigskip

\demo{\stag{1.11D} Main Fact}  For any approximation  $Y,i \in \dsize 
\bigcup_{k < k^*_Y} (W_k \cup W'_k)$  and $m \le n_Y$ and $\ell(*) \in 
\{1,2\}$ we can find an approximation  $Z$  such that: 
\medskip
\roster
\item "{$(\bigotimes)(\alpha)$}"   $n_Z = \text{ Max}\{m+1,n_Y\},I_Z \cap
{}^{m \ge} \text{ord } = I_Y \cap {}^{m \ge} \text{ord}$, \newline
(we mean $m$ not $n_Y$) and $k^*_Z = k^*_Y + 1$ 
\item "{${}(\beta)$}" $\quad (a) \quad$ 
if $\eta \in I_Y,\ell g(\eta) < m$  then        
$$
A^\ell_\eta[Z] = A^\ell_\eta[Z],
$$

$$
a^\ell_\eta[Z] = a^\ell_\eta[Z]
$$

$$
B^\ell_\eta[Z] = B^\ell_\eta[Y]
$$
     
$\,\,(b) \quad$ 
if  $\eta \in I_Y \cap I_Z,k < k^*_Y$ and $j \in W'_k$ then  \newline

$\quad \qquad
\,\,p^\ell_{\eta,j}[Z] = p^\ell_{\eta,j}[Y]$ \newline    

$\,\,(c) \quad$ 
if $\eta \in I_Y \cap I_Z,k < k^*_Y$ and $j \in W_k$ then \newline 

$\quad \qquad
\,\,b^\ell_{\eta,j}[Z] = b^\ell_{\eta,j}[Y]$  
\medskip
\noindent
\item "{$(\gamma)^1$}"  \underbar{if} $\eta \in I_Y,\ell g(\eta) = m,k <
k^*_Y$ and $i \in W_k$ and the element $b \in M^{\ell(*)}$ satisfies 
clauses $(a),(b)$ below \underbar{then} for some such
$b:A^{\ell(*)}_\eta[Z] = ac\ell(A^{\ell(*)}_\eta[Y] \cup \{b\})$; where
\medskip
\noindent
{\roster
\itemitem{ $(a)$ }  $b^{\ell(*)}_{\eta,i}[Y] \notin A^{\ell(*)}_\eta[Y]$ and
${\frac{b}{A^{\ell(*)}_\eta[Y]}} \perp_a A^{\ell(*)}_{\eta^-}[Y]$
\itemitem{ $(b)$ }  one of the conditions (i),(ii) listed below holds 
for $b$ \newline

$\qquad (i) \quad
b = b^{\ell(*)}_{\eta,i}[Y]$ and \newline

$\qquad \qquad \ell g(\eta) > 0 \Rightarrow
{\frac{b}{A^{\ell(*)}_\eta[Y]}} \perp_a A^{\ell(*)}_{\eta^-}[Y]$
\underbar{or} \newline
\medskip
\noindent
$\qquad (ii) \quad$
for no $b$ is $(i)$ satisfied (so $\ell g(\eta) > 0$) and 
$b \in M^{\ell(*)}$, \newline

$\qquad \qquad {\frac{b}{A^{\ell(*)}_\eta[Y]}} \perp_a 
A^{\ell(*)}_{\eta^-}[Y]$ and $b^\ell_{\eta,i}
\underset{A^{\ell(*)}_\eta[Y]} {}\to \biguplus b$ \underbar{or} \newline
\endroster}
\item "{$(\gamma)^2$}"  assume $\eta \in I_Y,\ell g(\eta) = m,k <
k^*_Y$ and $i \in W'_k$ then we have:      
{\roster
\itemitem{ $(a)$ }  \underbar{if} $p^{\ell(*)}_{\eta,i}$ is realized by some
$b \in M^{\ell(*)}$ such that \newline     
\smallskip
\noindent
$\qquad R \left( {\frac{b}{A^{\ell(*)}_\eta[Y]}},L,\infty \right) = 
R \left( p^{\ell(*)}_{\eta,i},L,\infty \right)$ and \newline  
\smallskip
\noindent
$\qquad \left[ \ell g(\eta) > 0 \Rightarrow {\frac{b}{A^{\ell(*)}_\eta[Y]}} 
\perp_a A^{\ell(*)}_{\eta^-}[Y] \right]$ \newline
\smallskip
\noindent
$\qquad$ \underbar{then} for some such  $b$ we have \newline
\smallskip
\noindent
$\qquad A^{\ell(*)}_\eta[Z] = 
ac\ell \left( A^{\ell(*)}_\eta[Y] \cup \{b\}\right)$
\medskip
\noindent
\itemitem{ $(b)$ }  \underbar{if} the assumption of clause $(a)$ fails and
$p^{\ell(*)}_{\eta,i}$ is realized \newline
\smallskip
\noindent
$\qquad$ by some 
$b \in  M^{\ell(*)} \backslash A^{\ell(*)}_\eta$ such that \newline
\smallskip
\noindent
$\qquad \left[ \ell g(\eta) > 0 \Rightarrow
{\frac{b}{A^{\ell(*)}_\eta[Y]}} \perp_a A^{\ell(*)}_{\eta^-}[Y] \right]$
\newline
\smallskip
\noindent
$\qquad$ \underbar{then} for some such  $b$ we have \newline
\smallskip
\noindent
$\qquad A^{\ell(*)}_\eta[Z] = ac\ell \left[ A^{\ell(*)}_\eta [Y] \cup \{b\}
\right]$
\endroster}
\item "{$(\delta)$}"  If $\eta \in I_Y$ and 
$\ell g(\eta) = m$, \underbar{then}  
$B^\ell_\eta[Z] = \{ b^\ell_{\eta,j}[Y]:j \in \cup \{W_k:k < k^*_Z \}\}$
\newline
is a countable subset of $M^\ell$,  containing  
$\{B^\ell_\nu[Z]:\nu \trianglelefteq \eta \text{ and } \nu \in Y\} \cup  
B^\ell_\eta[Y]$, with $B^\ell_\eta[Z] \prec M^\ell$ moreover
$B^\ell_\eta[Z] \subseteq_{na} M^\ell$ i.e. if $\bar a \subseteq B^\ell_\eta
[Z],\varphi(x,\bar y)$ is first order and
$(\exists x \in M^\ell \backslash ac\ell(\bar a))\,\varphi(x,\bar a)$ then
$(\exists x \in B^\ell_\eta[Z] \backslash ac\ell(\bar a))\varphi(\bar x,
\bar a))$ and 
$\{ a^\ell_{\eta \char 94 < \alpha >}[Y]:\eta \char 94 \langle \alpha
\rangle \in I_Y \text{ and } a^\ell_{\eta \char 94 < \alpha >}[Y] \notin
B^\ell_\eta[Z] \}$ is independent over
$(B^\ell_\eta[Z],A^\ell_\eta[Y])$  
\item "{$(\epsilon)$}"  if $\eta \in I_Y,\ell g(\eta) > m$, then $\eta \in  
I_Z \Leftrightarrow a^\ell_{\eta \restriction (m+1)}[Y] \notin  
B^\ell_{\eta \restriction m}[Z]$  
\item "{$(\zeta)$}"  if $\eta \in I_Y \cap I_Z,\ell g(\eta) > m$  then  
$A^\ell_\eta[Z] = ac\ell(A^\ell_\eta[Y] \cup A^\ell_{\eta \restriction m}[Z])$
and \newline
$B^\ell_\eta[Z] = B^\ell_\eta[Y]$  
\item "{$(\eta)$}"  if  $\eta \in I_Z \backslash I_Y$ then $\eta^- \in I_Y$
and $\ell g(\eta) = m+1$
\item "{$(\theta)$}"  $\{ p^\ell_{\eta,i}[Z]:i \in 
W'_{k^*_Z-1}\}$ is ``rich enough", e.g. include all finite types over
$A^\ell_\eta$
\item "{$(\iota)$}"  $\{b^\ell_{\eta,i}:i \in W_{k^*_z-1}\}$ list
$B^\ell_\eta[Z]$, each appearing infinitely often.
\endroster
\enddemo
\bigskip

\demo{Proof}  \underbar{First}
we choose $A^{\ell(*)}_\eta[Z]$  for  $\eta  \in  I$  of 
length $m$ according to condition $(\gamma) = (\gamma)^1 + (\gamma)^2$.  
(Note: one of the clauses $(\gamma)^1,(\gamma)^2$ necessarily holds as 
trivially
$\dsize \bigcup_k W_k \cap \dsize \bigcup_k W'_k = \emptyset$). \newline
\medskip

\noindent
\underbar{Second}, we choose (for such $\eta)$ an elementary mapping 
$f^Z_\eta$ 
extending $f^Y_\eta$ and a set  $A^{3-\ell (\ast )}_\eta [Z] \subseteq  
M^{3-\ell(*)}$ satisfying ``$f^Z_\eta$ is from $A^1_\eta[Z]$ onto  
$A^{3-\ell(*)}_\eta[Z]$" such that 
\medskip
\roster
\item "{$(*)_2$}"  if $m > 0$, then \newline
$f^Z_\eta \bigl(\text{tp}_\infty\bigl(\binom{A^1_\eta[Z]}{A^1_{\eta^-}[Y]}
,M_1 \bigr)\bigr) = \text{ tp}_\infty\bigl(\binom{A^2_\eta[Z]}
{A^2_{\eta^-}[Y]},M_2\bigr)$ 
\medskip
\noindent
\item "{$(*)_3$}"  if  $m = 0$, then  
$f^Z_\eta \left( \text{tp}_\infty \left( A^1_\eta[Z],M_1 \right) \right)  = 
\text{ tp}_\infty \left( A^2_\eta[Z],M_2 \right)$.
\endroster
\medskip

\noindent
[Why possible?  If we ask just the equality of  $\text{tp}_\alpha$ for 
an ordinal $\alpha$,  this follows by the first component of
$\text{tp}_{\alpha +1}$.  But (overshooting) for $\alpha \ge 
\left[ (\|M_1\| + \|M_2\|)^{\aleph_0} \right]^+$, equality of 
$\text{tp}_\alpha$ implies equality of $\text{tp}_\infty$].
\medskip

\noindent
\underbar{Third}, we choose $B^\ell_\eta[Z]$ for $\eta \in I_Y,\ell g(\eta)
 = m$ according to condition $(\delta)$ (here we use the countability of the 
language, you can do it by extending it  $\omega$  times) in both sides, i.e.
for $\ell = 1,2$.
\medskip

\noindent
\underbar{Fourth}, let $I' = \{\eta \in I:\text{ if } \ell g(\eta) > m$ then
$a^\ell_{\eta \restriction (m+1)}[Y] \notin B^\ell_{\eta \restriction m}[Z]\}$
(this will be $I_Y \cap I_Z$).
\medskip

\noindent
\underbar{Fifth}, we choose $A^\ell_\eta[Z]$ for $\eta \in I'$: if $\ell 
g(\eta) < m$, let $A^\ell_\eta[Z] = A^\ell_\eta[Y]$, if $\ell g(\eta) = m$
this was done, lastly if  $\ell g(\eta) > m$, let $A^\ell_\eta[Z] = 
ac\ell \left( A^\ell_\eta[Y] \cup A^\ell_{\eta \restriction m}[Z]\right)$.
\medskip

\noindent
\underbar{Sixth}, by induction on $k \le n$ we choose 
$f^Z_\eta$ for $\eta \in I'$ of length $k$: if $\ell g(\eta) < 
m$, let $f^Z_\eta = f^Y_\eta$, if $\ell g(\eta) = m$  this was done, lastly if
$\ell g(\eta) > m$ choose an elementary mapping from  $A^1_\eta $ onto  
$A^2_\eta$ extending $f^Y_\eta \cup f^Z_{\eta^-}$ (possible as
$f^Y_\eta \cup f^2_{\eta^-}$ is an elementary mapping and 
Dom$(f^Y_\eta) \cap \text{ Dom}(f^Z_{\eta^-}) =
A^{\ell(*)}_{\eta^-},\nonfork{\text{ Dom}(f^Y_\eta)}
{\text{Dom}(f^Z_{\eta^-})}_{A^{\ell(*)}_{\eta^-}}$ \newline
and $A^{\ell(*)}_{\eta^-} = ac \ell(A^{\ell(*)}_{\eta^-})$.)  
Now $f^Z_\eta$ satisfies clause (11) of Definition \scite{1.11A} when 
$\ell g(\eta) > m$  by applying \scite{1.8}(2).
\medskip

\noindent
\underbar{Seventh}, for $\eta \in I'$,  of length $< n_Z$, let
$v_\eta =:\{ \alpha:\eta \char 94 \langle \alpha \rangle \in I\}$, and we
choose $\{a^1_{\eta \char 94 < \alpha >}[Z]:\alpha \in u_\eta \},[\alpha \in  
u_\eta \Rightarrow \eta \char 94 \langle \alpha \rangle \notin I]$, a 
set of elements of $M^1$ realizing (stationary) regular types over 
$A^1_\eta[Z]$, orthogonal to $A_{\eta^-}[Y]$ when $\ell g(\eta) > 0$, 
such that it is independent over 
$\left( \cup \{a^1_{\eta \char 94 < \alpha >}[Y]:\eta \char 94 \langle \alpha
\rangle \in  I'\} \cup B^1_\eta[Z],A^1_\eta[Z] \right)$ and maximal under 
those restrictions.
Without loss of generality $\text{ sup}(v_\eta) < \text{ min}(u_\eta)$ and
for $\alpha_1 \in v_\eta \cup u_\eta$ and $\alpha_2 \in u_\eta$ we have: 
\medskip
\roster
\item "{$(*)_1$}"  if (for the given $\alpha_2$ and $\eta$) $\alpha_1$ is
minimal such that \newline
${\frac{a^1_{\eta \char 94 < \alpha_1 >}[Z]}
{A^1_\eta[Z]}} \pm {\frac{a^1_{\eta \char 94 < \alpha_2 >}[Z]}
{A^1_\eta[Z]}}$ then  
${\frac{a^1_{\eta \char 94 < \alpha_1 >}[Z]}{A^1_\eta[Z]}} = 
{\frac{a^1_{\eta \char 94 < \alpha_2 >}[Z]}{A^1_\eta[Z]}}$
\medskip
\noindent
\item "{$(*)_2$}"  \underbar{if} $\alpha_1 < \alpha_{2_1}$ and 
$a^1_{\eta \char 94 \langle \alpha_1 \rangle}[Z]/A^1_\eta[Z] =
a^1_{\eta \char 94 \langle \alpha_2 \rangle}[Z]/A^1_\eta[Z]$ and for some
$b \in M^1$  realizing 
${\frac{a^1_{\eta \char 94 < \alpha_2 >}[Z]}{A^1_\eta[Z]}}$ we have
\newline
$b \underset{A^1_{\eta^\ell}[Z]} {}\to \biguplus 
a^1_{\eta \char 94 < \alpha_2>}$ and  
\newline
$\text{tp}_\infty \left[ \binom{b}{A^1_{\eta \char 94 < \alpha_2 >}},M\right]
= \text{ tp}_\infty \left[ \left( 
{\frac{a^1_{\eta \char 94 < \alpha_1 >}}{A^1_\eta \char 94 < \alpha_2 >}}
\right),M \right]$ \newline
and $\alpha_1$ is minimal (for the given $\alpha_2$ and $\eta$)
\underbar{then} \newline
$\text{tp}_\infty \left[ \binom{a^1_{\eta \char 94 < \alpha_2 >}}
{A^1_{\eta \char 94 < \alpha_2 >}},M \right] = \text{ tp}_\infty \left[ \binom
{a^1_{\eta \char 94 < \alpha_1 >}}{A^1_{\eta \char 94 < \alpha_1 >}},M
\right]$.
\endroster
\medskip

\noindent
Easily (as in \cite[X]{Sh:c}) if $\alpha \in u_\eta,\eta 
\char 94 \langle \beta \rangle \in I'$ then
${\frac{a^1_{\eta \char 94 < \alpha >}[Z]}{A^1_\eta[Z]}}
\perp {\frac{a^1_{\eta \char 94 < \beta >}[Y]}{A^1_\eta[Y]}}$. 

For $\alpha \in u_\eta$ let $A^1_{\eta \char 94 < \alpha >}[Z] = 
ac\ell \left( A^1_\eta[Y] \cup \{a^1_{\eta \char 94 < \alpha >}[Z]\}\right)$.
\medskip

\noindent
\underbar{Eighth}, by the second component in the definition of 
$\text{tp}_{\alpha +1}$ (see Definition \scite{1.4}) we can choose 
(for $\alpha \in u_\eta)\,a^2_{\eta \char 94 < \alpha >}[Z],
A^2_{\eta \char 94 < \alpha >}[Z]$ and then 
$f^Z_{\eta \char 94 < \alpha >}$ as required (see (7) of 
Definition \scite{1.11A}).
\medskip

\noindent
\underbar{Ninth} and lastly, we let $I_Z = I' \cup \{\eta \char 94 < \alpha >
:\eta \in I',\ell g(\eta) < n_Z \text{ and } \alpha \in u_\eta\}$ and we 
choose  $B^\ell_\eta$ for $\eta \in I_Z \backslash I_Y$ and the $p^\ell_{\eta,
i},b^\ell_{\eta,j}$ as required (also in other case left). \newline
${}{}$ \hfill$\square_{\scite{1.11D}}$
\enddemo
\bigskip

\demo{\stag{1.11E} Finishing the Proof of 1.11}  We define by induction on
$n < \omega$ an approximation $Y_n = Y(n)$.  Let $Y_0$ be the trivial one 
(as in the proof of \scite{1.11}(C)). 

$Y_{n+1}$ is gotten from $Y_n$ as in \scite{1.11D} for 
$m_n,i_n \le n,\ell_n(*) \in \{1,2\}$  defined by reasonable bookkeeping 
(so $i_n \in 
\dsize \bigcup_{k < k^*_{Y(n)}} (W_k \cup W'_k)$ such that any triples 
appear infinitely often); 
without loss of generality: if $n_1 < n_2 \and \eta \in I^\ell_{n_1} \cap  
I^\ell_{n_2}$ then $\eta \in \dsize \bigcap^{n_2}_{n=n_1} I_n$.
\medskip

\noindent
Let $I^* = \text{ lim}(I^{Y(n)}_\ell) =: \{\eta:\text{ for every large 
enough } n,\eta \in I_n\}$ \newline
\medskip
for $\eta \in I^*$ let: $A^\ell_\eta[*] = \dsize \bigcup_{n < \omega}
A^\ell_\eta[Y_n],f^\ell_\eta[*] = \dsize \bigcup_{n < \omega} f^{Y(n)}_\eta$
and \newline

$B^\ell_\eta[*] = \dsize \bigcup_{n < \omega} B^\ell_\eta[Y_n]$.
\medskip

\noindent
Easily 
\medskip
\roster
\item "{$\bigoplus_1$}"   for  $\eta \in I^*,\langle B^\ell_\eta[Y_n]:n < 
\omega \text{ and } \eta \in I[Y_n] \rangle$  is an increasing sequence of 
$\subseteq_{na}$-elementary submodels of $M^\ell$
\endroster
\medskip

\noindent
hence 
\medskip
\roster
\item "{$\bigoplus_2$}"  for $\eta \in I^*,B^\ell_\eta[*] \subseteq_{na} 
M^\ell$
\item "{$\bigoplus_3$}"   $\nu \triangleleft \eta \in I^* \Rightarrow 
B^\ell_\nu[*] \subseteq B^\ell_\eta[*]$.
\endroster
\medskip

\noindent
[Why?  Because for infinitely many $n,m_n = \ell g(\eta)$  and clause 
$(\delta)$ of Main Fact \scite{1.11D}]. 
\medskip
\roster
\item "{$\bigoplus_4$}"  if $\eta \in I[Y_{n_1}] \cap I^*,\eta^- = \nu$ and
$n_1 \le n_2$ \underbar{then}

$$
A^\ell_\eta \left[ Y_{n_1} \right] \nonfork{}{}_{A^\ell_\nu[Y_{n_1}]} 
A^\ell_\nu[Y_{n_2}].
$$
\endroster
\medskip

\noindent
[Why?  Prove by induction on $n_2$ (using the non-forking calculus), for
$n_2 = n_1$ this is trivial, so assume $n_2 > n_1$.  If $m_{(n_2-1)} >
\ell g(\nu)$ we have $A^\ell_\nu[Y_{n_2}] = A^\ell_\nu[Y_{n_2-1}]$ (see
\scite{1.11D} and we have nothing to prove.  If $m_{(n_2-1)} < \ell g(\nu)$
then we note that $A^\ell_\nu[Y_{n_2}] = \text{ acl}(A^\ell_\nu[Y_{n_2-1}] 
\cup A^\ell_{\nu \restriction m_{(n_2-1)}}[Y_{n_2}])$ and
$\nonfork{A^\ell_\nu[Y_{n_2-1}]}{A^\ell_{\nu \restriction m_{(n_2-1)}}
[Y_{n_2}]}_{A^\ell_{\nu \restriction m_{(n_2-1)}}}$
(as $\nu \in I[Y_{n_2}]$, by \scite{1.11D} clause $(\delta)$ last phrase) and
now use clauses (5), (6) of Definition \scite{1.11D}.  Lastly if
$m_{(n_2-1)} = \ell g(\nu)$ again use $\nu \in I[Y_{n_2}]$ by \scite{1.11D},
clause $(\delta)$, last phrase]. 
\medskip
\roster
\item "{$\bigoplus_5$}"   if $\eta \in I[Y_{n_1}] \cap I^*,\eta^- = \nu$
and $n_1 \le n_2$ \underbar{then}

$$
{\frac{A^\ell_\eta[*]}{A^\ell_\nu[*] + a^\ell_\eta[*]}} \perp_a 
A^\ell_\nu.
$$
\endroster
\medskip

\noindent
[Why?  By clause (6) of Definition \scite{1.11A}, and 
orthogonality calculus]. 
\medskip
\roster
\item "{$\bigoplus_6$}"  if $\eta \in I^*$, then $A^\ell_\eta[*] \subseteq
B^\ell_\eta[*] \prec M^\ell$ moreover 
\item "{$\bigotimes_7$}"  $A^\ell_{\eta(*)}[*] \subseteq_{na} B^\ell_\eta[*]
\subseteq_{na} M^\ell$.
\endroster
\medskip

\noindent
[Why?  The second relation holds by $\bigotimes_2$.  The first relation we 
prove by induction on $\ell g(\eta)$;  clearly $A^\ell_\eta[*] = ac\ell 
(A^\ell_\eta[*])$ by clause (3) of Definition \scite{1.11A}; 
we prove this by induction on $m = \ell g(\eta)$, so 
suppose this is true for every $m' < m,m = \ell g(\eta),\eta \in I^*$, 
let $\varphi(x)$ be a formula with parameters in $A^\ell_\eta
[*]$  realized in $M^\ell$ as above say by $b \in M^\ell$ for some $n,
b \nonfork{}{}_{A^\ell_\eta[Y_n]} A^\ell_\eta[*]$.

So $\{\varphi(x)\} = p^\ell_{\eta,i}$ for some $i$ and for some $n' > n$
defining $Y(n+1)$ we have used \scite{1.11D} with $(\ell(*),i,m)$ there being
$(\ell,i,\ell g(\eta))$ here, so we consider clause $(\gamma)^2$ of
\scite{1.11D}.  So the case left is when the assumption of both clauses (a) 
and (b) of $(\gamma)^2$ fail, so we have $\ell g(\eta) > 0$ and

$$
b' \notin A^\ell_\eta[Y_{n'}],b' \in M^\ell \models \varphi[b'] \Rightarrow
{\frac{b'}{A^\ell_\eta[Y_{n'}]}} \pm A^\ell_{\eta^-}[Y_{n'}].
$$
\medskip

\noindent
We can now use the induction hypothesis (and \cite[5.3,p.292]{BeSh:307}).] 
\medskip
\roster
\item "{$\bigoplus_8$}"  if $\eta \in I^*$ and $\ell =1,2$, \underbar{then}
\newline
$\{ a^\ell_{\eta \char 94 < \alpha >}[*]:\eta \char 94 \langle \alpha 
\rangle \in I^*\}$ is a maximal subset of 
$$
\{ c \in M_\ell:\,{\frac{c}{A^\ell_\eta[*]}} \text{ regular, }
c \nonfork{}{}_{A^\ell_\eta[*]} B^\ell_\eta[*] \text{ and }
\ell g(\eta) > 0 \Rightarrow {\frac{c}{A^\ell_\eta[*]}} \perp  
A^\ell_{\eta^-}[*] \}
$$   
\medskip

\noindent
independent over $(A^\ell_\eta[*],B^\ell_\eta[*])$.
\endroster
\medskip

\noindent
[Why?  Note clause (7) of Definition \scite{1.11A} and clause 
$(\delta)$ of Main Fact \scite{1.11D}]. 
\medskip
\roster
\item "{$\bigotimes_9$}"  $A^\ell_{<>}[*] = B^\ell_{<>}[*]$
\endroster
\medskip

\noindent
[Why?  By the bookkeeping every  $b \in B^\ell_{<>}[*]$  is considered 
for addition to $A^\ell_{<>}[*]$ see \scite{1.11D}, clause $(\gamma)^1$,
subclause (b)(i) and for $\langle \rangle$ there is nothing to stop us]. 
\medskip
\roster
\item "{$\bigotimes_{10}$}"  if $\eta \in I^*$ and $p \in S(A^\ell_\eta[*])$
is orthogonal to $A^\ell_{\eta^-}[*]$ then 
${\frac{B^\ell_\eta[*]}{A^\ell_\eta[*]}} \perp p$ \newline
[why? if not, as $A^\ell_\eta[*] \subseteq B^\ell_\eta[*]$ by
\cite[Th.B,p.277]{BeSh:307} there is \newline
$c \in B^\ell_\eta[*] \backslash
A^\ell_\eta[*]$ such that: $\frac{c}{A^\ell_\eta[*]}$ is $p$.  As $c \in
B^\ell_\eta[*] = \dsize \bigcup_{n < \omega} B^\ell_\eta[Y_n]$, for every
$n < \omega$ large enough $c \in B^\ell_\eta[*]$, and $p$ does not fork
over $A^\ell_\eta[Y_n]$.  So for some such $n$ the triple $(i_n,\ell_n,m_n)$
is such that $\ell_n = \ell,m_n = \ell g(\eta)$ and $b^\ell_{\eta,i_n} = c$,
so by clause $(\gamma)^1$(b)(ii) of \scite{1.11D} we have $c \in A^\ell_\eta
[Y_n] \subseteq A^\ell_\eta[*]$.]
\item "{$\bigotimes_{11}$}"  if  $\eta \in I^*,\ell \in \{1,2\}$
\underbar{then} $\{a_{\eta \char 94 < \alpha >}:\eta \char 94 \langle \alpha
\rangle \in I^* \}$  is a maximal subset of
$\{ c \in M^\ell:\,{\frac{c}{A^\ell_\gamma[*]}} \text{ regular, } \perp
A^\ell_{\eta^-}[*] \text{ when meaningful}\}$ independent over 
$A^\ell_\eta[*]$.
\endroster
\medskip

\noindent
[Why?  If not, then for some $c \in M, 
\{a^\ell_{\eta \char 94 \langle \alpha \rangle}:\eta \char 94 \langle 
\alpha \rangle \in I^*\} \cup \{c\}$ is independent over $A^\ell_\eta[*]$ and
tp$(c,A^\ell_\eta[*])$ is regular (and stationary).  Hence by
$\otimes_{10}$ we have $\{a^\ell_\eta[Y_n]:\eta \char 94 \langle 
\alpha \rangle \in I^*\} \cup \{c\}$ is independent over $(A^\ell_\eta[*],
B^\ell_\eta[*])$.  Now for large enough $n$ we have 
$c \nonfork{}{}_{A^\ell_\eta[Y_n]} A^\ell_\eta[*]$ and by $\otimes_{10}$ 
we have
$\nonfork{c}{B^\ell_\eta[*]}_{A^\ell_\eta[Y_n]}$, hence 
$\nonfork{c}{B^\ell_\eta[Y_n]}_{A^\ell_\eta[*]}$, hence
$\{c\} \cup \{a^\ell_{\eta \char 94 \langle \alpha \rangle}[Y_n]:\eta
\char 94 \langle \alpha \rangle \in I[Y_n]\}$ is not independent over
$(A^\ell_\eta[Y_n],B^\ell_\eta[Y_n])$, but $\{a^\ell_{\eta \char 94
\langle \alpha \rangle}[Y_n]:\eta \char 94 \langle \alpha \rangle \in 
I[Y_n]\}$ is independent over $(A^\ell_\eta[Y_n],B^\ell_\eta[Y_n])$.  So
there is a finite set $w$ of ordinals such that $\alpha \in w \Rightarrow 
\eta \char 94 \langle \alpha \rangle \in I[Y_n]$ and
$\{c\} \cup \{a^\ell_{\eta \char 94 \langle \alpha \rangle}[Y_n]:\alpha \in
w\}$ is not independent over $(A^\ell_\eta[Y_n],B^\ell_\eta[Y_n])$, and
without loss of generality $w$ is minimal.
Now $\{a^\ell_{\eta \char 94 \langle \alpha \rangle}[*]:\eta \char 94
\langle \alpha \rangle \in I^*\} \cup B^\ell_\eta[*]$ include
$\{a^\ell_{\eta \char 94 \langle \alpha \rangle}[Y_n]:\eta \char 94 \langle
\alpha \rangle \in I[Y_n]\}$ hence it includes $\{a^\ell_{\eta \char 94 
\langle \alpha \rangle}[Y_n]:\alpha \in \omega\}$, easy contradiction.]
\medskip
\roster
\item "{$\bigoplus_{12}$}"  $f^*_\eta = \dsize \bigcup_{m < \omega}
f_\eta[Y_n]$ (for $\eta \in I^*$)  is an elementary map from  
$A^1_\eta[*]$ onto $A^2_\eta[*]$.
\endroster
\medskip

\noindent
[Easy]. 
\medskip
\roster
\item "{$\bigoplus_{13}$}"  $\dsize \bigcup_{\eta \in I^*} f^*_\eta$ is 
an elementary mapping from  $\dsize \bigcup_{\eta \in I^*} A^1_\eta[*]$ 
onto  $\dsize \bigcup_{\eta \in I^*} A^2_\eta[*]$.
\endroster
\medskip

\noindent
[Clear using by $\otimes_5 + \otimes_6$ and non-forking calculus].
\medskip
\roster
\item "{$\bigoplus_{14}$}"  We can find  $\langle d^\ell_\alpha:\alpha < 
\alpha(*) \rangle$  such that:  $d^\ell_\alpha \in M^\ell$, \newline  
$\text{tp}\bigl( d^\ell_\alpha,\dsize \bigcup_{\eta \in I[*]} A^\ell_\eta
[*] \cup \{d^\ell_\beta:\beta < \alpha \}\bigr)$ is $\aleph_\epsilon$
-isolated and $\bold F^\ell_{\aleph_0}$-isolated, and $g_\alpha = \dsize
\bigcup_{\eta \in I^*} f^\ast_\eta \cup \{\langle (d^1_\alpha,d^2_\alpha):
\alpha < \alpha(*)\rangle \}$  is an elementary mapping, $\alpha(*)$ is 
maximal.
\endroster
\medskip

\noindent
[Why?  We can choose by induction on $\alpha$, a member $d^1_\alpha$ of
$M^1 \backslash \dsize \bigcup_{\eta \in I[*]} \cup \{d^1_\beta:\beta <
\alpha\}$ such that tp$(d^1_\alpha,\dsize \bigcup_{\eta \in I[*]}
A^\ell_\eta[*] \cup \{d^1_\beta:\beta < \alpha\})$ is
$\aleph_\varepsilon$-isolated and $\bold F^\ell_{\aleph_0}$-isolated.  So for
some $\alpha(*)$, $d^1_\alpha$ is well defined iff $\alpha < \alpha(*)$ (as
$\beta < \alpha \Rightarrow d^1_\beta \ne d^1_\alpha \in M^1$).  Now choose
by induction on $\alpha < \alpha(*),d^2_\alpha \in M^2$ as required above,
possible by ``$M^2_i$ being $\aleph_\varepsilon$-saturated (see
\cite[XII,2.1,p.591]{Sh:c}, \cite[IV,3.10,p.179]{Sh:c}.]
\medskip
\roster 
\item "{$\bigotimes_{15}$}"  Dom$(g_{\alpha(*)})$, Rang$(g_{\alpha(*)})$ are 
universes of elementary submodels of $M^1,M^2$ respectively called 
$M'_1,M'_2$ respectively.
\endroster
\medskip

\noindent
[Why?  See \cite[XII,1.2(2),p.591]{Sh:c} and the proof of $\otimes_{14}$].
\medskip
\roster
\item "{$\bigotimes_{16}$}"  $M^\ell \ne M'_\ell$ then for some $d \in M_\ell 
\backslash M'_\ell,{\frac{d}{M'_\ell}}$ is regular.
\endroster
\medskip

\noindent
[Why?  By \cite[Th.5.9,p.298]{BeSh:307} as $N^\ell_\eta \subseteq_{na} M^\ell$
by $\bigotimes_7$]. 
\medskip
\roster
\item "{$\bigotimes_{17}$}"  if $M^\ell \ne M'_\ell$ then for some $\eta \in  
I^*$,  there is $d \in  M^\ell \backslash M'_\ell,{\frac{d}{A^\ell_\eta[*]}}$
is regular, $d \nonfork{}{}_{A^\ell_\eta[*]} M'_\ell,
\left[ \ell g(\eta) > 0 \Rightarrow  {\frac{d}{A^\ell_\eta[*]}} \perp 
A^\ell_{\eta^-}[*] \right]$.
\endroster
\medskip

\noindent
[Why?  By \cite[XII,1.4,p.529]{Sh:c} every non-algebraic $p \in S(M'_\ell)$
is not orthogonal to some $A^\ell_\eta[*]$ so by $\otimes_{16}$ we
can choose $\eta \in I^*$ and $d \in  M^\ell \backslash M'_\ell$ such that
${\frac{d}{M'_\ell}}$ regular $\pm A^\ell_\eta[*]$; without loss of 
generality  $\ell g(\eta)$  is minimal,
now $A^\ell_\eta[*] \subseteq_{\text{na}} M^\ell$ and by 
\cite[4.5,p.290]{BeSh:307} without loss of generality  
$d \nonfork{}{}_{A^\ell_\eta[*]} M'_\ell$; the last clause is by 
``$\ell g(\eta)$ minimal"]. 
\medskip
\roster
\item "{$\bigoplus_{18}$}"   $M_\ell  = M'_\ell$.
\endroster
\medskip

\noindent
[Why?  By $\bigoplus_{11} + \bigoplus_{17}$].
\medskip
\roster
\item "{$\bigoplus_{19}$}"  there is an isomorphism from  $M_1$ onto 
$M_2$ extending \newline
$\dsize \bigcup_{\eta \in I^*} f^*_\eta$.
\endroster
\medskip

\noindent
[Why?  By $\bigoplus_{14} + \bigotimes_{15}$ we have $M'_1 \cong M'_2$, so
by $\otimes_{18}$ we are done].  \hfill$\square_{\scite{1.11E}}
\,\,\square_{\scite{1.11}}$
\enddemo
\bigskip

\proclaim{\stag{1.12} Lemma}   Assume $B \nonfork{}{}_{A} C,A = ac\ell(A) = 
B \cap C$ and $A,B,C$ are $\epsilon$-finite, $A \cup B \cup C \subseteq M,M$ 
an $\aleph_\epsilon$-saturated model of $T$.  For notational simplicity make
$A$ a set of individual constants.
\newline
\underbar{Then}  $\text{tp}_{{\Cal L}_{\infty,\aleph_\epsilon}(d.q.)}
(B + C;M) = \text{ tp}_{{\Cal L}_{\infty,\aleph_\epsilon}(d.q.)}(B;M) + 
\text{ tp}_{{\Cal L}_{\infty,\aleph_\epsilon}(d.q.)}[C;M]$ \newline
where
\endproclaim
\bigskip

\definition{\stag{1.12A} Definition}  1) For any logic ${\Cal L}$ let, 
$\bar b$ a sequence from a model $M$,

$$
\align
\text{tp}_{\Cal L}(\bar b;M) = \biggl\{ \varphi(\bar x):&M \models 
\varphi[B],\varphi \text{ a formula in the vocabulary of }  M, \\
  &\text{from the logic }  {\Cal L} \text{ (with free variables from} \\
  &\bar x, \text{ where } \bar x = \langle x_i:i < \ell g(\bar b)\rangle)
\biggr\}.
\endalign
$$
\medskip

\noindent
2)  Replacing  $\bar b$  by a set  $B$  means we use the variables  
$\langle x_b:b \in  B \rangle$. \newline
3)  Saying  $p_1 = p_2 + p_3$ in \scite{1.12} 
means that we can compute  $p_1$ from  $p_2$ and $p_3$ 
(and the knowledge how the variables fit and the knowledge
of $T$, of course).
\enddefinition
\bigskip

\demo{Proof of the Lemma \scite{1.12}} \newline
It is enough to prove: 
\enddemo
\bigskip

\proclaim{\stag{1.12B} Claim}  Assume    
\medskip
\roster
\item "{$(a)$}"  $M^1,M^2$ are $\aleph_\epsilon$-saturated and
\item "{$(b)$}"  $A^\ell_1 \nonfork{}{}_{A^\ell_0} A^\ell_2$
\item "{$(c)$}"  $A^\ell_0 = ac\ell(A^\ell_0)$ and $A^\ell_i$ 
is $\epsilon$-finite for $\ell =1,2$ and $i < 3$
\item "{$(d)$}"  for $m = 0,1,2$ we have 
$f_m:A^1_m \overset\text{onto}\to\longrightarrow A^2_m$ is an elementary
mapping preserving $\text{tp}_\infty$ (in $M^1,M^2$ respectively) and
\item "{$(e)$}"  $f_0 \subseteq f_1,f_2$.
\endroster
\medskip

\noindent
\underbar{Then} there is an isomorphism from $M^1$ 
onto $M^2$ extending $f_1 \cup f_2$.
\endproclaim
\bigskip

\demo{Proof of \scite{1.12B}}  Repeat the proof of \scite{1.2}, 
but starting with $Y_0$ such that \newline
$A^\ell_\eta[Y_0] = A^\ell_0,A^\ell_\eta[Y_0] = A^\ell_1,A^\ell_{<1>}
[Y_0] = A^\ell_2,f^{Y_0}_{<>} = f_0,f^{Y_0}_{<0>} = f_1,f^{Y_0}_{<1>}
= f_2$ and that  $\langle \rangle,\langle 0 \rangle,\langle 1 \rangle$ belongs
to all $I[Y_0]$.  During the construction we preserve $\langle 0 \rangle,
\langle 1 \rangle \in I[Y_n]$ and for this $B^\ell_{<>}[Y_n]
\nonfork{}{}_{A^\ell_0} A^\ell_1 \cup A^\ell_2$ [you can use 
$\aleph_\epsilon$-saturation!]. \hfill$\square_{\scite{1.12}}$
\enddemo
\newpage

\head {\S2 Finer Types} \endhead
\resetall
\bigskip

We shall use here alternative types showing us probably a finer way to 
manipulate tp.
\bigskip

\demo{\stag{2.1} Convention}  $T$  is superstable,  NDOP 

$M,N$ are $\aleph_\epsilon$-saturated $\prec {\frak C}^{\text{eq}}$.
\enddemo
\bigskip

\definition{\stag{2.2} Definition}  $\Gamma_p = 
\biggl\{ \binom{\bar b}{\bar a}:\bar a \subseteq \bar b \text{ are } 
\epsilon$-finite$\biggr\}$

$$
\Gamma_r = \biggl\{ \binom{p}{\bar a}:\bar a \text{ is } 
\epsilon\text{-finite, }
p \in S(\bar a) \text{ is regular (so stationary)}\biggr\}
$$

$$ 
\align
\Gamma_{pr} = \biggl\{ \binom{p,r}{\bar a}:&\,\bar a \text{ is } \epsilon
\text{-finite, }  p  \text{ is a regular type of depth } > 0, \\
  &\,p \pm \bar a \text{ (really only the equivalence class } p / \pm
\text{ matters),} \\
  &\,r = r(x,\bar y) \in S(\bar a) \text{ is such that for } (c,\bar b)
\text{ realizing }  r, \\
  &\,c/(\bar a + \bar b) \text{ is regular } \pm p, \text{ and }
{\frac{\bar b}{\bar a}} = (r \restriction \bar y) \perp  p \biggr\}.
\endalign
$$
\medskip

\noindent
We may add (to $\Gamma_x$) superscripts: 
\medskip
\roster
\item "{$(\alpha)$}"  $f$ if $\bar a$ (or $\bar a,\bar b$) is finite
\item "{$(\beta)$}"   $s:\text{ for } \Gamma_p$ if ${\frac{\bar b}{\bar a}}$
is stationary, for $\Gamma_r$ if $p$ is stationary and for 
$\Gamma_{pr}$ if $r$ is stationary and every automorphism of ${\frak C}$ 
over $\bar a$ fix $p/\pm$ 
\item "{$(\gamma)$}"  $c$ if $\bar a$ (or $\bar a,\bar b$) are algebraically 
closed.
\endroster
\enddefinition
\bigskip

\proclaim{\stag{2.3} Claim}  If $p$ is regular of depth $> 0$ and 
$p \pm \bar a$ and $\bar a$ is $\epsilon$-finite \underbar{then} 
for some $\bar a',\bar a \subseteq \bar a' \subseteq ac\ell(\bar a)$ 
and for some $q$ \underbar{we have}
$\binom{p,q}{\bar a} \in \Gamma^s_{pr}$. \newline
[Why?  Use e.g. \cite[V,4.11,p.272]{Sh:c} assume ${\frac{\bar b}{\bar a}} \pm
p$;  we can define inductively equivalence relations  $E_n$,  with 
parameters from  $ac\ell(\bar a^\ell),\bar a^\ell = \bar a \char 94(\bar b
/E_0) \char 94 \cdots \char 94 (\bar b/E_{n-1})$, 
such that  $\text{tp}(\bar b/E_n,ac\ell(\bar a^n))$ is semi-regular. By 
superstability this stop for some $n$ hence $\bar b \subseteq 
ac\ell(\bar a^n)$.  For some first  $m \,\text{ tp}(\bar b/E_m,
ac\ell(\bar a^n))$ is $\pm p$, by \cite[X,7.3(5),p.552]{Sh:c} the type 
is regular (as because $p$ is trivial having depth  $> 0$; see 
\cite[X,7.2,p.551]{Sh:c}.
\endproclaim
\bigskip

\definition{\stag{2.4} Definition}  We define by induction on an 
ordinal $\alpha$ 
the following (simultaneously): note --- if a definition of something 
depends on another which is not well defined, neither is the something) 

$$
\text{tp}^1_\alpha\bigl[\binom{p}{\bar a},M \bigr] \qquad \text{ for }
\qquad \binom{p}{\bar a} \in \Gamma_r,\bar a \subseteq M
$$

$$
\text{tp}^2_\alpha\bigl[\binom{p,r}{\bar a},M\bigr] \qquad \text{ for }
\qquad \binom{p,r}{\bar a} \in \Gamma_r,\bar a \subseteq M
$$

$$
\text{tp}^3_\alpha\bigl[\binom{\bar b}{\bar a},M\bigr] \qquad \text{ for }
\qquad \binom{\bar b}{\bar a} \in \Gamma^c_p,
\bar a \subseteq \bar b \subseteq M
$$
\enddefinition
\bigskip

\noindent
\underbar{Case A} \,\,\underbar{$\alpha = 0$}: $\text{tp}^1_\alpha \left[ 
\binom{p}{\bar a},M \right]$ is  $\text{tp}((c,\bar a),\emptyset)$ 
for any  $c$  realizing  $p$.

$$
tp^2_\alpha \bigl[\binom{p,r}{\bar a},M \bigr] \text{ is tp}((c,\bar b,
\bar a),\emptyset) \text{ for any } (c,\bar b) \text{ realizing } r
$$

$$
\text{tp}^3_\alpha \bigl[\binom{\bar b}{\bar a},M \bigr] \text{ is tp}
((\bar b,\bar a),\emptyset).
$$
\bigskip

\noindent
\underbar{Case B} \,\,\underbar{$\alpha = \beta + 1$}: \newline
\smallskip
\noindent
(1)  $\text{tp}^1_\alpha \left[ \binom{p}{\bar a},M \right]$ is:  
\medskip
\noindent
\underbar{Subcase B1}  \,\,\underbar{if} $p$ has depth zero, 
it is $w_p(M/\bar a)$ 
(the $p$-weight, equivalently, the dimension)  
\medskip
\noindent
\underbar{Subcase B2}  \,\,\underbar{if} $p$  has depth  $> 0$  
(hence is trivial),  then it is  $\{\langle y,\lambda^y_{\bar a,p} \rangle:
y\}$ where    

$$
\align
\lambda^y_{\bar a,p} = &\text{ dim}({\Cal I}^y_{\bar a,p}[M],a) \text{ where }
{\Cal I}^y_{\bar a,p}[M] = \\        
  &\,\biggl\{ c \in M:c \text{ realize }  p \text{ and } y = \text{ tp}^3
_\beta \left[ \binom{ac\ell(\bar a + c)}{ac\ell(\bar a)},M \right] \biggr\}
\endalign
$$
\medskip

\noindent
an alternative probably more transparent and simpler in use:

$$
\align
\lambda^y_{\bar a,p} = \text{ dim}\biggl\{ c \in M:&\text{ realizes }
p \text{ and} \\
  &\,y = \{ \text{tp}^3_\beta \left[ \binom{ac\ell(\bar a + c')}
{ac\ell(\bar a)},M \right]:c' \in p(M) \text{ and } c' \nonfork{}{}_{\bar a}
c \}.
\endalign
$$  
\medskip

\noindent
(2) $\text{tp}^2_\alpha \left[ \binom{p,r}{\bar a},M \right]$ is: \newline    
\medskip
\roster
\item "{{}}"  $\text{tp}^1_\alpha \left[ \binom{c/\bar b^+}{\bar b^+},M
\right]$  for any  $(c,\bar b)$  realizing  $r,\bar b^+ = ac\ell(\bar a + 
\bar b)$  (so not well defined if we get at least two different 
cases; so remember  $c/b^+ \in  S(\bar b^+)$).
\endroster
\medskip

\noindent
(3) $\text{tp}^3_\alpha \left[ \binom{\bar b}{\bar a},M \right]$  is  
$\biggl\{ \left< p,\text{tp}^2_\alpha \left[ \binom{p,r}{\bar b},M \right]
\right>:\binom{p,r}{\bar b} \in \Gamma^s_{pr} \text{ and}: p \pm \bar a
\biggr\}$.
\bigskip

\noindent
\underbar{Case C} \,\,\underbar{$\alpha$ limit}:  For any $\ell \in \{1,2,3\}$
and suitable object OB:

$$
\text{tp}^\ell_\alpha[OB,M] = \langle \text{tp}^\ell_\beta[OB,M]:\beta < 
\alpha \rangle.
$$
\bigskip

\definition{\stag{2.5} Definition}  1) For $\binom{p}{\bar a} \in \Gamma_r$ 
where $\bar a \in M$, let (remembering \scite{1.6}(8)):

$$
\align
{\Cal P}^M_{\binom{p}{\bar a}} = \biggl\{ q \in M:&\,q \text{ regular and}: 
q \pm p \text{ or for some} \\
  &\,c \in p(M),q \in {\Cal P}^M_{\binom{c}{\bar a}} \biggr\}.
\endalign
$$
\medskip

\noindent
2)  For  $\binom{p,r}{\bar a} \in \Gamma_{pr}$ let

$$
\align
{\Cal P}^M_{\binom{p,r}{\bar a}} = \biggl\{ q \in S(M):&\,q 
\text{ regular and}: q \pm p \text{ or for some} \\
  &\,(c,\bar b) \in r(M),q \in {\Cal P}^M_{\binom{c}
{\bar a + \bar b}} \biggr\}.
\endalign
$$
\medskip

\noindent
3)  For a set  ${\Cal P}$  is stationary regular types not orthogonal to  
$M_1$,  let  $M_1 \le_{\Cal P} M_2$ means  $M_1 \prec  M_2$ and for every  $p
\in {\Cal P}$  and  $\bar c \in  M_2,{\frac{\bar c}{M_1}} \perp p$. \newline
4)  If (in (3)), ${\Cal P} = {\Cal P}^{M_1}_{\binom{p}{\bar a}}$ we may 
write  $\binom{p}{\bar a}$ instead ${\Cal P}$, similarly if ${\Cal P} = 
{\Cal P}^{M_1}_{\binom{p,r}{\bar a}}$ we may write $\binom{p,r}{\bar a}$.
\enddefinition
\bigskip

\proclaim{\stag{2.6} Claim}: \newline
1)  From  $\text{tp}^1_\alpha \bigl[ \binom{p}{\bar a},M \bigr]$  we 
can compute $\text{tp}^1_\infty \bigl[ \binom{p}{\bar a},M \bigr]$ if  
$\text{Dp}(p) < \alpha$. \newline
\smallskip
\noindent
2)  From $\text{tp}^2_\alpha \bigl[ \binom{p,q}{\bar a},M \bigr]$  we 
can compute $\text{tp}^2_\infty \bigl[ \binom{p,q}{\bar a},M \bigr]$ if  
$\text{Dp}(p) < \alpha$. \newline
\smallskip
\noindent
3)  From $\text{tp}^3_\alpha \bigl[ \binom{\bar b}{\bar a},M \bigr]$  we 
can compute $\text{tp}^3_\infty \bigl[ \binom{b}{\bar a},M \bigr]$ if  
$\text{Dp}(\bar b/\bar a) < \alpha$.
\endproclaim
\bigskip

\demo{Proof}  By the definition.
\enddemo
\bigskip

\demo{\stag{2.7} Observation}  From $\text{tp}^\ell_\alpha(OB,M)$ we can 
compute $\text{tp}^\ell_\beta[OB,M]$ if $\beta \le \alpha$ and the former 
is well defined.
\enddemo
\bigskip

\proclaim{\stag{2.8} Lemma}  For every ordinal $\alpha$ the following holds:
\medskip
\roster
\item   $\text{tp}^1_\alpha$ is well defined \footnote{i.e. in all the 
cases we have tried to define it in Definition \scite{2.9}}
\item   $\text{tp}^2_\alpha$ is well defined.
\item   $\text{tp}^3_\alpha$ is well defined.
\item   If  $\bar a \in M_1,\binom{p}{a} \in \Gamma_r,M_1 
\le_{\binom{p}{\bar a}} M_2$ \underbar{then} \newline
$\text{tp}^1_\alpha \left[ \binom{p}{\bar a},M_1 \right] = 
\text{tp}^1_\alpha \left[ \binom{p}{\bar a},M_2 \right]$.
\item   If  $\bar a \in M_1,\binom{p,r}{\bar a} \in \Gamma^s_{pr},M_1 
\le_{\binom{p}{\bar a}} M_2$ \underbar{then} \newline
$\text{tp}^2_\alpha \left[ \binom{p,r}{\bar a},M_1 \right] = 
\text{tp}^2_\alpha \left[ \binom{p,r}{\bar a},M_2 \right]$.
\item  If $\bar a \subseteq \bar b \subseteq M_1,\binom{\bar b}{\bar a} 
\in \Gamma^c_p,M_1 \le_{\binom{\bar b}{\bar a}} M_2$ \underbar{then} \newline
$\text{tp}^3_\alpha \left[ \binom{\bar b}{\bar a},M_1 \right] = 
\text{tp}^3_\alpha \left[ \binom{\bar b}{\bar a},M_2 \right]$.
\endroster
\endproclaim
\bigskip

\demo{Proof}  We prove it, by induction on $\alpha$, simultaneously (for all 
clauses and parameters). 

If $\alpha$ is zero, they hold trivially by the definition. 

If  $\alpha$ is limit, they hold trivially by the definition and induction 
hypothesis.

\noindent
So for the rest of the proof let $\alpha = \beta + 1$.
\enddemo
\bigskip

\demo{Proof of $(1)_\alpha$}  If $p$ has depth zero --- check directly.    

If $p$ has depth $> 0$ - by $(3)_\beta$    

(i.e. induction hypothesis) no problem.
\enddemo
\bigskip

\demo{Proof of $(2)_\alpha$}  Like \scite{1.8} (and $(4)_\alpha$).
\enddemo
\bigskip

\demo{Proof of $(3)_\alpha$}  Like $(2)_\alpha$.
\enddemo
\bigskip

\demo{Proof of $(4)_\alpha$}  Like \scite{1.7} (and $(3)_\beta,(6)_\beta$).
\enddemo
\bigskip

\demo{Proof of $(5)_\alpha$}  By $(2)_\alpha$ we can look only at  
$(c,\bar b^+)$ in $M_1$,  then use $(4)_\alpha$.
\enddemo
\bigskip

\demo{Proof of $(6)_\alpha$}  By $(5)_\alpha$.  \hfill$\square_{\scite{2.8}}$
\enddemo
\bigskip

\proclaim{\stag{2.9} Lemma}  For an ordinal $\alpha$ restricting ourselves to 
the cases (the types $p,p_1$ being) of depth $< \alpha$:
\medskip
\roster
\item "{$(A1)$}"  Assume $\binom{p}{\bar a} \in \Gamma_r,\bar a \subseteq  
\bar a_1 \subseteq M,\bar a_1$ is $\epsilon$-finite, ${\frac{\bar a_1}
{\bar a}} \perp p$ and $p_1$ is the stationarization of $p$  
over $\bar a_1$.  \newline
\underbar{Then} from  $\text{tp}^1_\alpha \left[ \binom{p}{\bar a},M \right]$
we can compute $\text{tp}^1_\alpha \bigl[ \binom{p_1}{\bar a_1},M \bigr]$.
\item "{$(A2)$}" Under the assumption of $(A1)$ all the inverse 
computations are O.K.
\item "{$(A3)$}"  Assume  $\binom{p_\ell}{\bar a} \in \Gamma_r$ for $\ell = 
1,2,\bar a \subseteq M$  and  $p_1 \pm p_2$. \newline
\underbar{Then} from  $\text{tp}^1_\alpha \bigl[ \binom{p_1}{\bar a_1},M
\bigr]$  we can compute  $\text{tp}^1_\alpha \left[ \binom{p_2}{\bar a},M
\right]$.
\item "{$(B1)$}"  Assume  $\binom{p_\ell,r_\ell}{\bar a} \in \Gamma^{sc}_{pr}$
for $\ell = 1,2,\bar a \in M$ and  $p_1 \pm  p_2$. \newline
\underbar{Then} (from the first order information on $\bar a,p_1,p_2,r_1,r_2$,
of course, and) $\text{tp}^2_\alpha \left[\binom{p,r_2}{\bar a},M \right]$ 
we can compute  $\text{tp}^2_\alpha \left[ \binom{p,r_2}{\bar a},M \right]$.
\item "{$(B2)$}"  Assume  $\bar a \subseteq \bar a_1 \subseteq M,
{\frac{\bar a_1}{\bar a}} \pm p,\binom{p,r}{\bar a} \in \Gamma^{cs}_{pr},r 
\subseteq r_1 \in S(\bar a_1),r_1$ does not fork over $\bar a,
\bigl(\text{ so } \binom{p,r_1}{\bar a_1} \in \Gamma_{pr}\bigr)$. \newline
\underbar{Then} from  $\text{tp}^2_\alpha \left[ \binom{p,r_1}{\bar a},M
\right]$  we can compute  $\text{tp}^2_\alpha \left[ \binom{p,r_2}{\bar a}
,M \right]$.
\item "{$(B3)$}"  Under the assumption of $(B2)$, the inverse 
computation are O.K.
\item "{$(C1)$}"  Assume  $\binom{\bar b}{\bar a} \in \Gamma^c,\bar a 
\subseteq \bar b \subseteq M,\bar a \subseteq \bar a_1,\bar b 
\nonfork{}{}_{\bar a} \bar a_1,\bar b_1 = ac\ell(\bar a_1 + \bar b)$.
\newline
\underbar{Then} from  $\text{tp}^3_\alpha \bigl[\binom{\bar b}{\bar a},M
\bigr]$  we can compute $\text{tp}^3_\alpha \bigl[\binom{\bar b_1}
{\bar a_1},M \bigr]$.
\item "{$(C2)$}"  Under the assumptions of $(C1)$ the inverse computation 
is 0.K.
\item "{$(C3)$}"  Assume $\binom{\bar b}{\bar a} \in \Gamma_p,\bar b 
\subseteq b^*,{\frac{\bar b^*}{\bar b}} \perp_a \bar a,
\bar b^* = ac\ell(\bar b^*)$.  Then from $\text{tp}^3_\alpha
\bigl[\binom{\bar b}{\bar a},M \bigr]$ we can compute \newline
$\biggl\{ \text{tp}^3_\alpha \left[ \binom{\bar b'}{\bar a},M \right]:
\bar b \subseteq b' \subseteq M \text{ and }
{\frac{\bar b'}{\bar b}} = {\frac{\bar b^*}{\bar b}} \biggr\}$.
\endroster
\endproclaim
\bigskip

\demo{Proof}  We prove it, simultaneously, for all clauses and parameters, by 
induction on  $\alpha$ and the order of the clauses. 

For  $\alpha = 0$:  easy. 

For  $\alpha$ limit: very easy. 

So assume $\alpha = \beta + 1$.
\enddemo
\bigskip

\demo{Proof of $(A1)_\alpha$}  As $p$ is stationary $\perp 
{\frac{\bar a_1}{\bar a}}$, for every  $c \in p(M),{\frac{c}{\bar a}} \vdash
{\frac{c}{\bar a_1}}$,  which necessarily is $p_1$, hence $p(M) = p_1(M)$.  
Also the dependency relation on $p(M)$ is the same over $\bar a_1$, hence 
dimension.  So it suffices to show: 
\medskip
\roster
\item "{$(*)$}"  for  $c \in p(M)$, from $\text{tp}^3_\beta \bigl[ 
\binom{ac\ell(\bar a + c)}{ac\ell \, \bar a},M \bigr]$ we can compute  
\newline
$\text{tp}^3_\beta \bigl[\binom{ac\ell(\bar a_1 + c)}{ac\ell \, \bar a_1},M
\bigr]$.
\endroster
\medskip

\noindent
But this holds by $(C1)_\beta$.
\enddemo
\bigskip

\demo{Proof of $(A2)_\alpha$}   Similar using $(C2)_\beta$.
\enddemo
\bigskip

\demo{Proof of $(A3)_\alpha$}   If $p_1$ (equivalently $p_2$) has 
depth zero --- the dimensions are equal.  Assume they have depth  
$> 0$  hence are trivial and dependency over  
$\bar a$  is an equivalence relation on  $p_1(M) \cup  p_2(M)$. 

Now for $c_1 \in  p_1(M)$,  from  $\text{tp}^3_\beta 
\bigl[\binom{ac\ell(a + c_1)}{ac\ell(\bar a)},M \bigr]$ we can compute 
for every complete type over  
$ac\ell(\bar a + c_1)$ not forking over $\bar a$, and $\bar d$ realizing  
$r,\text{tp}^3_\beta \bigl[\binom{ac\ell(\bar a + \bar d + c_1)}{ac\ell
(\bar a + \bar d)},M \bigr]$ - by $(C1)_\beta$, then we can compute for 
each such $r,\bar d$,

$$
\align
\biggl\{ \text{tp}^3_\beta \left[ \binom{ac\ell(\bar a + \bar d + c_2)}
{ac\ell(\bar a + \bar d)},M \right]:&\,c_2 \in p_2(M) \text{ and }  
{\frac{c_2}{ac\ell(\bar a + \bar d + c_1)}} \perp_a (\bar a + \bar d) \\
  &\text{ (necessarily }  c_2 \nonfork{}{}_{\bar a} \bar d)\biggr\}
\endalign
$$
\medskip

\noindent
(this by $(C3)_\beta$).
\enddemo
\bigskip

\demo{Proof of $(B1)_\alpha$}   As in earlier cases we can restrict 
ourselves to the case  $\text{Dp}(p_\ell) > 0$.  We can find $(c_\ell,
\bar b_\ell) \in r_\ell(M),\bar b_1 \nonfork{}{}_{\bar a} \bar b_2,c_1
\bar b_1 \nonfork{}{}_{\bar a} \bar b_2$ (by \cite[X,7.3(6)]{Sh:c}].  
By \scite{2.8}(2) (and the definition) from  $\text{tp}^2_\alpha \bigl[ 
\binom{p_1,r_1}{\bar a},M \bigr]$  we can compute that it is equal to  
\newline
$\text{tp}^1_\alpha \bigl[\binom{c_1/ac\ell(\bar a + \bar b_1)}{ac\ell(\bar a
+ b_1)},M \bigr]$. \newline
\medskip

\noindent
By $(A1)_\alpha$ we can compute 
$\text{tp}^1_\alpha \bigl[\binom{c_1/ac\ell(\bar a + \bar b_1 + \bar b_2)}
{ac\ell(\bar a + b_1 + \bar b_2)},M \bigr]$ hence by $(A3)_\alpha$ we
can compute $\text{tp}^1_\alpha 
\bigl[\binom{c_2/ac\ell(\bar a + \bar b_1 + \bar b_2)}{ac\ell(\bar a
+ \bar b_1 + \bar b_2)},M\bigr]$. \newline
\medskip

\noindent
Now use $(A2)_\alpha$ to compute  
$\text{tp}^1_\alpha \bigl[\binom{c_2/ac\ell(\bar a + \bar b_2)}
{ac\ell(\bar a + \bar b_2)},M \bigr]$ and by \scite{2.8}(2), \scite{2.4}(2) 
it is equal to $\text{tp}^2_\alpha \bigl[\binom{p,r}{\bar a},M \bigr]$.
\enddemo
\bigskip

\demo{Proof of $(B2)_\alpha$}   Choose $(c,\bar b) \in  r(M),c \bar b 
\nonfork{}{}_{\bar a} \bar a_1$. 

\relax From  $\text{tp}^2_\alpha \bigl[\binom{p,r_1}{\bar a},M \bigr]$ we can 
compute $\text{tp}^1_\alpha \bigl[\binom{c/(\bar a + \bar b)}
{\bar a + \bar b},M \bigr]$ (just --- see \scite{2.8}(2) and 
Definition \scite{2.4}), 
from it we can compute $\text{tp}^1_\alpha \bigl[\binom{c/(\bar a + \bar b
+ \bar a_1)}{(\bar a + \bar b + \bar a_1)},M \bigr]$ (by $(A1)_\alpha)$,
from it we can compute $\text{tp}^2_\alpha \bigl[\binom{p,r_2}{\bar a_2},
M \bigr]$ (see \scite{2.8}(2) and Definition \scite{2.4}).
\enddemo
\bigskip

\demo{Proof of $(B3)_\alpha$}   Let  $\binom{p,r}{\bar b_1} \in \Gamma^s_r,p 
\perp \bar a_1$ be given.  So necessarily  ${\frac{\bar a_1}{\bar a}} \pm p$
(this to enable us to use $(B2,3)$.  It suffices to compute  
$\text{tp}^2_\alpha \bigl[\binom{p,r}{\bar b_1},M \bigr]$ and we can 
discard the case  $\text{Dp}(p) = 0$. 

So $p$ is regular $\pm \bar b_1,\perp \bar a_1$, hence $p \pm \bar b,p \perp
\bar a$,  and as  $\bar a \subseteq \bar b,\bar b = ac\ell(\bar b)$ we can 
find  $r,\binom{p,r_1}{\bar b} \in \Gamma_{pr}$, (see \scite{2.3}) 
and we know  
$\text{tp}^2_\alpha \left[ \binom{p,r_1}{\bar b},M \right]$, and we can find 
$r_2$, a complete type over  $\bar b_1$ extending  $r_1$ which does not 
fork over  $\bar b_1$.  From  $\text{tp}^2_\alpha \left[ \binom{p,r_1}{\bar b}
,M \right]$ we can compute $\text{tp}^2_\alpha \bigl[\binom{p,r_2}
{\bar b_1},M \bigr]$ by $(B2)_\alpha$, and from it   
$\text{tp}^2_\alpha \bigl[\binom{p,r}{\bar b_1},M \bigr]$ by $(B1)_\alpha$.
\enddemo
\bigskip

\demo{Proof of $(C2)_\alpha$}  Similar, use $(B3)_\alpha$ instead of 
$(B2)_\alpha$.
\enddemo
\bigskip

\demo{Proof of $(C3)_\alpha$}  Without loss of generality ${\frac{\bar b^*}
{\bar b}}$ is semi regular, let $p^*$ be a regular type not orthogonal to 
it and without loss of generality  $\text{Dp}(p^*) > 0 \Rightarrow
{\frac{\bar b^*}{\bar b}}$  regular (as in \scite{2.3}). 

If $p^*$ has depth zero, then the only  $p$  appearing in the definition  
$\text{tp}^3_\alpha \bigl( \bigl[ {\frac{\bar b}{\bar a}} \bigr],M \bigr)$
is $p^*$ (up to $\pm$)  and this is easy.  Then  
$\text{tp}^2_\alpha$ is just the dimension and we have no problem. 

So assume $p^*$ has depth  $> 0$.  We can by $(B1)_\alpha,(B2)_\alpha$ 
compute $\text{tp}^2_\alpha \bigl[\binom{p',q'}{\bar b^*},M \bigr]$ when
$p' \pm \bar b,p' \pm  p^*$ (regardless of the choice of $\bar b^*$).  
Next assume  $p' \pm  p^*$; by $(B1)_\alpha$ without loss of generality $q'$
does not fork over $\bar b$.  As  $\text{Dp}(p^*) > 0$, it is trivial (and 
we assume $w_p(\bar b^*,\bar b) = 1)$ hence $\bar b^* /\bar b$ is regular so 
in $\text{tp}^2_\alpha \bigl[\binom{p,q'}{\bar b^*},M \bigr]$ we just lose 
a weight 1 for one specific  $\text{tp}^3_\beta$ type: the one $\bar b^*$ 
realizes concerning which we have a free choice.  We are left with the 
cases  $p' \pm  \bar b,p' \pm  p^*$; well we know $\text{tp}^3_\beta$ but 
we have to add  $\text{tp}^3_\alpha$?   Use Claim \scite{2.6}(3) 
(and $(A1)_\alpha$ as we add a parameter). \hfill$\square_{\scite{2.9}}$
\enddemo
\bigskip

\proclaim{\stag{2.10} Claim}  $\text{tp}^3_\gamma 
\bigl[\binom{\bar b}{\bar a},M\bigr],
\text{tp}^3_\gamma[\bar a,M],\text{tp}^3_\gamma[M]$  are 
expressible by formulas in \newline
${\Cal L}^\gamma_{\infty,\aleph_\epsilon}(d.q.)$.
\endproclaim
\bigskip

\noindent
By \scite{2.9} we have
\demo{\stag{2.11} Conclusion}  If  $\text{Dp}(T) < \infty$ then: \newline
1)  From  $\text{tp}^3_\infty \bigl[\binom BA,M \bigr]$  we can compute  
$\text{tp}_\infty \bigl[ \binom BA,M \bigr]$ (type --- from \S1).
\newline
2)  Similarly from $\text{tp}^3_\infty[A,M]$ we can compute $\text{tp}_\infty
[(A),M]$.
\enddemo
\bigskip

\noindent
\relax From \scite{2.6}, \scite{2.10}, \scite{2.11} and \scite{1.11} we get
\proclaim{\stag{2.12} Corollary}  If  $\gamma = \text{Dp}(T)$ and $M,N$ are 
$\aleph_\epsilon$-saturated, then

$$
M \cong  N \Leftrightarrow  \text{ tp}^3_\gamma[M] = \text{ tp}^3_\gamma[N] 
\Leftrightarrow  M \cong_{{\Cal L}^\gamma_{\infty,\aleph_\epsilon}(d.q.)}N.
$$
\endproclaim
\newpage

\head {Appendix} \endhead
\bigskip

The following clarifies several issues raised by Baldwin.  A consequence of 
\medskip
\roster
\item "{$\bigotimes$}"  the existence of nice invariants for 
characterization up to isomorphism (or characterization of the models up to 
isomorphism by their ${\Cal L}$-theory for suitable logic ${\Cal L})$
\endroster
\medskip

\noindent
naturally give absoluteness, e.g. extending the universe say by nice forcing 
preserve non-isomorphism.  So negative results for 
\medskip
\roster
\item "{$(*)$}"  is non-isomorphism (of models of $T$) preserved by forcing 
by ``nice forcing notions"?
\endroster
\medskip

\noindent
implies that we cannot characterize models up to isomorphism by their 
${\Cal L}$-theory when the logic  ${\Cal L}$  is ``nice", i.e. when  
$Th_{\Cal L}(M)$  preserved by nice forcing notions. 
So coding a stationary set by the isomorphism type can be interpreted 
as strong evidence of ``no nice invariants",  see \cite{Sh:220}.  
Baldwin, Laskowski, Shelah \cite{BLSh:464} show that not only for 
every unsuperstable; but also for some quite 
trivial superstable (with NDOP, NOTOP) countable  $T$,  there are 
non-isomorphic models which can be made isomorphic by some ccc (even 
$\sigma$-centered) forcing notion.  This shows that the lack of a really 
finite characterization is serious. \newline
Can we still get from the characterization in this paper an absoluteness 
result?
Note that for preserving $\aleph_\epsilon$-saturation (for simplicity for 
models of countable  $T$)  we need to add no reals \footnote{(the set of
$\{ac\ell(\bar a):\bar a \in {}^{\omega >}M\}$ is absolute 
but the set of their enumeration and of the  
$\{f \restriction(ac\ell(\bar a)):f \in \text{ AUT}({\frak C}),f(\bar a) = 
\bar a\}$ is not).}, and in order not to erase distinction of dimensions 
we want not to collapse cardinals, so the following questions is natural, 
for a first order (countable) complete  $T$: 
\medskip
\roster
\item "{$(*)^1_T$}"  assume  $V_1 \subseteq V_2$ are transitive models of ZFC 
with the same cardinals and reals, the theory $T \in V_1$.  If the models
$M_1,M_2$ are from $V_1$ and they are models of $T$ not isomorphic in $V_1$;
must they still be not isomorphic in $V_2$ \footnote{Note we did not say 
they have the same $\omega$-sequences of ordinals, 
e.g. if $V_2 = V^P_1,P$ Prikry forcing, then the 
assumption of $(*)_T$ holds though a new $\omega$-sequence of ordinals was
added.  So for  $V_1 \subseteq V_2$ as in $(*)_T$, the 
${\Cal L}_{\infty,\aleph_1}$-theory is not necessarily preserved.}
\item "{$(*)^2_T$}"  like $(*)^1_T$ we assume in addition 
${\Cal P}(|T_1|)^{V_1} = {\Cal P}(|T_1|)^{V_2}$.
\endroster
\medskip

Of course, for countable  $T$  the answer is negative even for 
$\aleph_\epsilon$-saturated models except for superstable, NDOP, NOTOP 
theories so we restrict ourselves to these.  It should be quite transparent 
that ${\Cal L}_{\infty,\aleph_\epsilon}(q.d.)$-theory is preserved from  
$V_1$ to $V_2$ (as well as the set of sentences in the logic), hence for 
the class of $\aleph_\epsilon$-saturated models (of superstable NDOP, NOTOP
theory $T$) the answer to $(*)^2_T$ is yes.
\newpage

\newpage
    
REFERENCES.  
\bibliographystyle{lit-plain}
\bibliography{lista,listb,listx,listf}

\def\germ{\frak} \def\scr{\cal}
  \ifx\documentclass\undefinedcs\def\rm{\fam0\tenrm}\fi
  \def\defaultdefine#1#2{\expandafter\ifx\csname#1\endcsname\relax
  \expandafter\def\csname#1\endcsname{#2}\fi} \defaultdefine{Bbb}{\bf}
  \defaultdefine{frak}{\bf} \defaultdefine{mathbb}{\bf}
  \defaultdefine{beth}{BETH} \def\bbfI{{\Bbb I}} \def\mbox{\hbox}
  \def\text{\hbox} \def\om{\omega} \def\Cal#1{{\bf #1}} \def\pcf{pcf}
  \defaultdefine{cf}{cf} \defaultdefine{reals}{{\Bbb R}}
  \defaultdefine{real}{{\Bbb R}} \def\restriction{{|}} \def\club{CLUB}
  \def\w{\omega} \def\exist{\exists} \def\se{{\germ se}} \def\bb{{\bf b}}
  \def\equivalence{\equiv} \let\lt< \let\gt> \def\cite#1{[#1]}
  \def\implies{\Rightarrow}
\begin{thebibliography}{BLSh 464}
\makeatletter \renewcommand{\@biblabel}[1]{[#1]} \makeatother

\bibitem[BLSh 464]{BLSh:464}John~T. Baldwin, Michael~C. Laskowski, and Saharon
  Shelah.
\newblock {Forcing Isomorphism}.
\newblock {\em {Journal of Symbolic Logic}}, {\bf 58}:1291--1301, 1993.

\bibitem[BeSh 307]{BeSh:307}Steven Buechler and Saharon Shelah.
\newblock {On the existence of regular types}.
\newblock {\em {Annals of Pure and Applied Logic}}, {\bf 45}:277--308, 1989.

\bibitem[Sh 220]{Sh:220}Saharon Shelah.
\newblock {Existence of many $L_ {\infty,\lambda}$-equivalent, nonisomorphic
  models of $T$ of power $\lambda$}.
\newblock {\em {Annals of Pure and Applied Logic}}, {\bf 34}:291--310, 1987.
\newblock Proceedings of the Model Theory Conference, Trento, June 1986.

\bibitem[Sh 225]{Sh:225}Saharon Shelah.
\newblock {On the number of strongly $\aleph_ \epsilon$-saturated models of
  power $\lambda$}.
\newblock {\em {Annals of Pure and Applied Logic}}, {\bf 36}:279--287, 1987.
\newblock See also [Sh:225a].

\bibitem[Sh:c]{Sh:c}Saharon Shelah.
\newblock {\em {Classification theory and the number of nonisomorphic models}},
  volume~92 of {\em {Studies in Logic and the Foundations of Mathematics}}.
\newblock {North-Holland Publishing Co., Amsterdam, xxxiv+705 pp}, 1990.

\end{thebibliography}

\shlhetal
\enddocument

\bye